\documentclass[12pt]{article}
\def\Arxived{extended}
\newif\ifarxived
\ifx\Arxived\Kaba\else\arxivedtrue\fi
\ifarxived
\ifarxived\else
\documentclass[12pt]{article}
\fi
\usepackage[bookmarks=true,bookmarksnumbered=true,
  bookmarkstype=toc,
  pdftitle={Strong downward Löwenheim-Skolem theorems for stationary logics, II --- 
  reflection down to the continuum},%
  pdfauthor={Sakaé Fuchino, André Ottenbreit Maschio 
  Rodrigues and Hiroshi Sakai},%
  pdfsubject={},%
  pdfkeywords={Strong Downward Löwenheim Skolem Theorem, stationary logic, generically 
  large cardinals, mixed support iteration, Laver 
  function, superhuge cardinal, continuum problem}]{hyperref}
\else
\usepackage[dvipdfmx,bookmarks=true,bookmarksnumbered=true,
  bookmarkstype=toc,
  pdftitle={Strong downward Löwenheim-Skolem theorems for stationary logics, II --- 
  reflection down to the continuum},%
  pdfauthor={Sakaé Fuchino (渕野 昌), André Ottenbreit Maschio 
  Rodrigues and Hiroshi Sakai (酒井 拓史)},%
  pdfsubject={},%
  pdfkeywords={Strong Downward Löwenheim Skolem Theorem, stationary logic, generically 
  large cardinals, mixed support iteration, Laver 
  function, superhuge cardinal, continuum problem}]{hyperref}
\fi
\ifarxived
\else
\usepackage[dvipdfmx]{pxjahyper}
\fi
\hypersetup{
  colorlinks=true,
  linkcolor=blue,
  filecolor=blue,      
  citecolor=cyan,
  urlcolor=cyan
}
\ifarxived
\usepackage{graphicx, color}
\else
\usepackage[dvipdfmx]{graphicx, color}
\fi
\usepackage{amsmath, amssymb}
\usepackage{bbm}
\usepackage{dsfont}
\usepackage{bbold}
\usepackage{mathtools}
\usepackage[normalem]{ulem}
\usepackage{setspace}
\newif\ifextended
\ifarxived\extendedtrue\else
\ifx\Extended\Kaba\else\extendedtrue\fi
\fi
\newif\ifJapanese
\newif\iftesting
\ifextended\testingtrue\fi
\makeatletter
\ifextended\input{SDLS-x.aux}\else\input{SDLS.aux}\fi
\makeatother

\usepackage{theorem}
\newcommand{\bbd}[1]{{\mathbb{#1}}}
\theorembodyfont{\ifJapanese\rm\else\it\fi}
\newcount\minute	
\newcount\hour		
\newcount\hourMins  
\ifJapanese
\def\today%
{
  \the\year 年\,\zeroPadTwo{\the\month}月\,\zeroPadTwo{\the\day}日%
}
\fi
\def\now%
{
  \minute=\time    
  \hour=\time \divide \hour by 60 
  \hourMins=\hour \multiply\hourMins by 60
  \advance\minute by -\hourMins 
  \zeroPadTwo{\the\hour}:\zeroPadTwo{\the\minute}%
}
\def\zeroPadTwo#1%
{
  \ifnum #1<10 0\fi    
  #1
}
\title{Strong downward L\"owenheim-Skolem theorems for stationary logics, II\\--- 
  reflection down to the continuum}
\author{Saka\'e Fuchino (\ifarxived\hspace{4em}\else 渕野 昌\fi), Andr\'e Ottenbreit Maschio 
  Rodrigues\\and Hiroshi Sakai (\ifarxived\hspace{5em}\else 酒井 拓史\fi)${}^\ast$}
\date{}

\renewcommand{\baselinestretch}{1.2}
\renewcommand{\thefootnote}{(\arabic{footnote})\,}
\iftesting
\fi
\iftesting
\newcommand{\Label}[1]{\label{#1}\marginpar{{\renewcommand{\baselinestretch}{0.4}\tiny 
		  #1}}}
\else
\newcommand{\Label}[1]{\label{#1}}
\fi

\def\memo#1{\iftesting\marginpar{{\normalsize\renewcommand{\baselinestretch}{0.4}\tiny%
			#1\par}}\else\fi}%
\newcounter{frml}[section]
\newcounter{frmla}[section]
\def\thefrml{{\arabic{section}.\arabic{frml}}}
\def\thefrmla{{$\aleph$\arabic{section}.\arabic{frmla}}}
\def\frmlabel#1{\refstepcounter{frml}{\def\baka{#1}\ifx\baka\empty\else\label{#1}\fi}%
{\rm({\thefrml})\hfill\hfill\hfill}}
\def\frmlabela#1{\refstepcounter{frmla}{\def\baka{#1}\ifx\baka\empty\else\label{#1}\fi}%
{\rm({\thefrmla})\hfill\hfill\hfill}}
\def\xitem[#1]{\item[\frmlabel{#1}]\mbox{}%
	\iftesting\marginpar{{\renewcommand{%
				\baselinestretch}{0.6}\tiny#1}}\fi\ignorespaces}
\def\xitemq[#1]{\item[\frmlabel{#1}]\mbox{}%
	\ignorespaces}
\def\xitemd[#1]#2{\item[(\ref{#1})$#2$\hfill\hfill\hfill]}
\def\xitemA[#1]{\item[\frmlabela{#1}]\mbox{}%
	\iftesting\marginpar{{\renewcommand{%
				\baselinestretch}{0.6}\tiny#1}}\fi\ignorespaces}
\def\xitemsub[#1]#2{\item[\frmlabel{#1}$_{#2}$]\mbox{}%
	\iftesting\marginpar{{\renewcommand{%
				\baselinestretch}{0.6}\tiny#1}}\fi\ignorespaces}
\def\xxitem[#1][#2]{\item[(\ref{#1}{\makebox[1.4ex][c]{#2}})]\mbox{}%
	\iftesting\marginpar{{\renewcommand{%
				\baselinestretch}{0.6}\tiny\{#1\}\{#2\}}}\fi\ignorespaces}
\def\xitemof#1{{\rm({\ref{#1}})}}
\def\xitemAof#1{{\rm({\ref{#1}})}}

\newenvironment{xitemize}{\begin{list}{}{\parsep=0.5\smallskipamount%
			\itemindent=-0.4ex%
			\itemsep=0.5\smallskipamount\leftmargin=4em\labelwidth=3em\labelsep=0.7em}}%
							 {\end{list}}
\def\assert#1{\noindent\makebox[4.8ex][r]{\rm(\makebox[2.2ex][c]{#1})}\ \ \ignorespaces}
\def\wassertof#1{\makebox[4.8ex][r]{\rm(\makebox[2.2ex][c]{#1})}}%
\def\assertof#1{(#1)}%

\def\daimaru#1{\makebox[1zw][c]{\mbox{\leavevmode\lower.08zh\hbox{%
        \rlap{\hbox to 0.76zw{\hfil\mbox{}\hfill{}\raisebox{0.03zh}{\scalebox{1.2}{○}}\hfil}}%
        \raise0.19zh\hbox to 1zw{\hfil{\hspace{0.16zw}\footnotesize#1}\hfil}}}}\,}

\newtheorem{Thm}{\ifJapanese{\bf 定理}\else {\bf Theorem}\fi}[section]

\newtheorem{ThmA}{\ifJapanese{\bf 定理\,A\!}\else{\bf Theorem\,A\!}\fi}[section]

\newtheorem{Prop}[Thm]{\ifJapanese{\bf 命題}\else{\bf Proposition}\fi}

\newtheorem{Lemma}[Thm]{\ifJapanese{\bf 補題}\else{\bf Lemma}\fi}
\newtheorem{LemmaA}[ThmA]{\ifJapanese{\bf 補題\,A\!}\else{\bf Lemma\,A\!}\fi}
\newtheorem{Cor}[Thm]{\ifJapanese{\bf 系}\else{\bf Corollary}\fi}

\newtheorem{Claim}{{\bf Claim}}[Thm]
\newtheorem{ClaimA}{{\bf Claim\,A\!}}[ThmA]

\newcommand{\prf}{\ifJapanese{\bf 証明．\ }\ignorespaces\else{\bf 
		Proof.\ \ }\ignorespaces\fi}

\newcommand{\prfofClaim}{\raisebox{-.4ex}{\Large $\vdash$\ \ }}

\newcommand{\prfofClaimA}{\prfofClaim}
\newcommand{\prfof}[1]{\ifJapanese{\bf #1 の証明．\ \ }%
	\ignorespaces\else{\bf Proof of #1:}\ \ \ignorespaces\fi}
\newcommand{\Thmof}[1]{\ifJapanese{定理\,\ref{#1}}\else{Theorem~\ref{#1}}\fi}

\newcommand{\Lemmaof}[1]{\ifJapanese{補題\,\ref{#1}}\else{Lemma~\ref{#1}}\fi}

\newcommand{\LemmaAof}[1]{\ifJapanese{補題\,A\,\ref{#1}}\else{Lemma\,A\,\ref{#1}}\fi}

\newcommand{\Propof}[1]{\ifJapanese{命題\,\ref{#1}}\else{Proposition~\ref{#1}}\fi}

\newcommand{\Corof}[1]{\ifJapanese{系\,\ref{#1}}\else{Corollary~\ref{#1}}\fi}

\newcommand{\Claimof}[1]{{Claim \ref{#1}}}
\newcommand{\ClaimAof}[1]{{Claim\,A\,\ref{#1}}}

\newcommand{\Thmabove}{{\ifJapanese 定理\else Theorem\fi\ \number\theThm}}

\newcommand{\Claimabove}{{Claim \number\theClaim}}

\newsavebox{\qedbox}\sbox{\qedbox}{
{\unitlength=0.05mm \begin{picture}(40,60)
\put(0,0){\framebox(30,44)[cc]{}}
\put(30,-7){\rule{7\unitlength}{44\unitlength}}
\put(10,-7){\rule{27\unitlength}{7\unitlength}}
\end{picture}}}
\newcommand{\qed}{\mbox{}\hfill\usebox{\qedbox}}
\newcommand{\smallqed}%
{\mbox{}\smallskip\hfill\raisebox{-.4ex}{\Large $\dashv$}}
\newcommand{\qedof}[1]%
{\mbox{} \hspace*{\fill}{\usebox{\qedbox}{\tiny~(#1)}}}
\newcommand{\Qedof}[1]%
{\mbox{} \hspace*{\fill}{\usebox{\qedbox}%
{\tiny~(#1~\number\theThm)}}}
\newcommand{\QedAof}[1]%
{\mbox{} \hspace*{\fill}{\usebox{\qedbox}%
{\tiny~(#1~\number\theThmA)}}}
\newcommand{\qedofThm}{\Qedof{\ifJapanese 定理\else Theorem\fi}}

\newcommand{\qedofCor}{\Qedof{\ifJapanese 系\else Corollary\fi}}
\newcommand{\qedofProp}{\Qedof{\ifJapanese 命題\else Proposition\fi}}
\newcommand{\qedofLemma}{\Qedof{\ifJapanese 補題\else Lemma\fi}}
\newcommand{\qedofLemmaA}{\QedAof{\ifJapanese 補題\,A\!\else Lemma\,A\!\fi}}

\newcommand{\qedskip}{\medskip}

\newcommand{\qedofClaim}%
{\mbox{}\hfill\raisebox{-.4ex}{\Large $\dashv$ }\nolinebreak%
\mbox{\tiny~(Claim~\number\theClaim)}}
\newcommand{\qedofClaimA}%
{\mbox{}\hfill\raisebox{-.4ex}{\Large $\dashv$ }\nolinebreak%
\mbox{\tiny~(Claim~A\,\number\theClaimA)}}
\newcommand{\qedofClaimAof}[1]%
{\mbox{}\hfill\raisebox{-.4ex}{\Large $\dashv$ }\nolinebreak%
\mbox{\tiny~(Claim~A\,\ref{#1})}}
\newcommand{\qedofSubclaim}%
{\mbox{}\hfill\raisebox{-.4ex}{\Large $\dashv$ }\nolinebreak%
\mbox{\tiny~(Subclaim~\number\theSubclaim)}}
\newcommand{\cardof}[1]{\mathopen{|\,}#1\mathclose{\,|}}
\newcommand{\ccardof}[1]{\mathopen{\,\|}#1\mathclose{\|\,}}
\newcommand{\ulsetof}[1]{\mathopen{\,|}#1\mathclose{|\,}}
\newcommand{\Card}{{\it Card\/}}
\newcommand{\setof}[2]{\{#1\,:\,#2\}}
\newcommand{\ssetof}[1]{\{#1\}}

\newcommand{\subseteqand}[1]{\mathrel{\mathop{\subseteq}%
		\limits_{\scriptscriptstyle\hbox to 14pt{$\scriptscriptstyle #1$\hss}}}}

\newcommand{\mapping}[3]{#1:#2\rightarrow #3}

\newcommand{\elembed}[3]{#1:#2\stackrel{\preccurlyeq\hspace{0.8ex}}{\rightarrow}#3}

\newcommand{\fnsp}[2]{\mbox{}^{{#1}\hspace{-0.02em}}#2}
\newcommand{\imageof}{{}^{\,{\prime}{\prime}}}
\newcommand{\ro}{\mathop{\mathrm r\mathrm o}}
\newcommand{\seqof}[2]{\langle#1\,:\,#2\rangle}
\newcommand{\pairof}[1]{\langle#1\rangle}
\newcommand{\psof}[1]{{\mathcal P}\/(#1)}
\newcommand{\forces}[2]{\,\|\hspace{-.35ex}\mbox{\sf--}_{\,#1\,}%
\mbox{\rm``}\,#2\,\mbox{\rm''}}
\newcommand{\notforces}[2]{\rlap{\rm\ 
 /}\|\hspace{-.35ex}\mbox{\sf--}_{\,#1\,}%
 \mbox{\rm``}\,#2\,\mbox{\rm''}}
\newcommand{\CHECK}[2]{{\surd\!}_{#1}(#2)}
\newcommand{\modelof}[1]{\models\!\mbox{\rm``\,}#1\mbox{\rm''}}

\newcommand{\crit}{\mbox{\it crit\/}}
\newcommand{\bbone}{{\mathord{\mathbb{1}}}}
\newcommand{\bbzero}{{\mathord{\mathbb{0}}}}
\newcommand{\circleq}{\mathrel{{\leqslant}%
		\hspace{-0.86ex}{\lower-0.53ex\hbox{$\scriptscriptstyle\circ$}}}}
\newcommand{\symb}[1]{{\mathord{\hspace{0.08em}\underbracket[0.6pt][2pt]{#1}}\hspace{0.08em}}}

\newcommand{\variables}[2]{{#1}_0\ctentenc {#1}_{#2}}

\newcommand{\restr}{\restriction}

\newcommand{\cf}{\mathop{cf\/}}

\newcommand{\Col}{{\rm Col}}

\newcommand{\Fn}{{\rm Fn}}

\newcommand{\trcl}{\mathop{\mbox{\it trcl\/}}}
\newcommand{\dom}{\mathop{\it dom}}

\newcommand{\poC}{\bbd{C}}
\newcommand{\poP}{\bbd{P}}
\newcommand{\poQ}{\bbd{Q}}

\newcommand{\On}{{\rm On}}

\newcommand{\genG}{\mathbb{G}}
\newcommand{\utildea}{\utilde{a}}
\newcommand{\utildegenG}{\utilde{\mathbb{G}}}
\newcommand{\utildecalG}{\utilde{\cal{G}}}
\newcommand{\utildepoP}{\utilde{\mathbb{P}}}

\newcommand{\utildepoQ}{\utilde{\mathbb{Q}}}
\newcommand{\utpoQ}{\utilde{\mathbb{Q}}}
\newcommand{\utpoR}{\utilde{\mathbb{R}}}

\newcommand{\genH}{\mathbb{H}}

\newcommand{\condp}{\mathbbm{p}}
\newcommand{\condq}{\mathbbm{q}}
\newcommand{\condr}{\mathbbm{r}}

\newcommand{\LT}{{<}\,}

\newcommand{\ctenten}{,\mbox{}\hspace{0.08ex}{.}{.}{.}\hspace{0.1ex}}
\newcommand{\ctentenc}{,{}\linebreak[0]\hspace{0.04ex}{{.}{.}{.}\hspace{0.1ex},\,}\linebreak[0]}

\newcommand{\xmbox}[1]{ $\relax{\rm #1}\relax$ }
\newcommand{\veca}{{\vec a}}

\newcommand{\gmA}{\mathfrak{A}}
\newcommand{\gmB}{\mathfrak{B}}

\newcommand{\continuum}{2^{\aleph_0}}

\newcommand{\calB}{{\mathcal B}}
\newcommand{\calC}{{\mathcal C}}
\newcommand{\calD}{{\mathcal D}}

\newcommand{\calF}{{\mathcal F}}
\newcommand{\calG}{{\mathcal G}}
\newcommand{\calH}{{\mathcal H}}

\newcommand{\calL}{{\mathcal L}}

\newcommand{\calP}{{\mathcal P}}

\newcommand{\calS}{{\mathcal S}}
\newcommand{\calT}{{\mathcal T}}
\newcommand{\calU}{{\mathcal U}}

\newcommand{\calY}{{\mathcal Y}}

\newcommand{\uta}{\utilde{a}}

\newcommand{\utC}{\utilde{C}}
\newcommand{\utS}{\utilde{S}}
\newcommand{\utildeg}{\utilde{g}}
\newcommand{\utcalU}{\utilde{\calU}}

\newcommand{\varin}{\mathrel{\varepsilon}}
\newcommand{\notvarin}{\mathrel{{\not}\,\varepsilon\,}}
\newcommand{\ZF}{{\sf ZF}}
\newcommand{\ZFC}{{\sf ZFC}}
\newcommand{\CH}{{\sf CH}}

\newcommand{\SCH}{{\sf SCH}}

\newcommand{\MA}{{\sf MA}}

\newcommand{\PFA}{{\sf PFA}}

\newcommand{\IC}{{\sf IC}}
\newcommand{\IS}{{\sf IS}}
\newcommand{\IU}{{\sf IU}}

\newcommand{\RP}{{\sf RP}}

\newcommand{\DRP}{{\sf DRP}}

\newcommand{\SDLS}{{\sf SDLS}}
\newcommand{\intnl}{int}

\newcommand{\st}{such that}
\newcommand{\wrt}{with respect to}
\newcommand{\Wolog}{Without loss of generality}
\newcommand{\wolog}{without loss of generality}

\newcommand{\tfae}{the following are equivalent}

\newcommand{\po}{poset}
\newcommand{\pos}{posets}
\newcommand{\uniV}{{\sf V}}
\newcommand{\Ba}{Boolean algebra}

\newcommand{\cBa}{complete Boolean algebra}

\newcommand{\Pkl}[2]{\ifx\bakakaba#1\bakakaba\ifx\bakakaba#2\bakakaba{\mathcal 
    P}_\kappa(\lambda)\else{\mathcal P}_\kappa(#2)\fi\else{\mathcal P}_{#1}(#2)\fi}

\newcommand{\utildeT}[1]{%
  \hbox to 0pt{\smash{$\mathop{\textstyle #1}\limits_{%
			\raisebox{0.4ex}[0pt]{$\scriptstyle\sim$}}$}\hss}%
  \relax\phantom{\mathord{{#1}_{\rule[-0.6ex]{0pt}{1pt}}}}}
\newcommand{\utildeS}[1]{%
	\hbox to 0pt{\smash{$\mathop{\scriptstyle #1}\limits_{%
				\raisebox{0.6ex}[0pt]{$\scriptscriptstyle\sim$}}$}\hss}%
	\relax\phantom{\mathord{{#1}_{\rule[-0.6ex]{0pt}{1pt}}}}}
\newcommand{\utildeSS}[1]{%
	\hbox to 0pt{$\mathop{\scriptscriptstyle #1}%
		\limits_{\scriptscriptstyle\sim}$\hss}%
		\relax\phantom{\underline{#1}}}
\newcommand{\utilde}[1]{%
	\mathchoice{\utildeT{#1}}{\utildeT{#1}}{\utildeS{#1}}{\utildeSS{#1}}}

{\end{minipage}\end{trivlist}}
\newcommand{\baA}{\bbd{A}}
\newcommand{\baB}{\bbd{B}}

\begin{document}
\maketitle
\renewcommand{\thefootnote}{$\ast$\ }
  \footnotetext{Graduate School of System Informatics, Kobe University \\Rokko-dai 1-1, Nada, Kobe 657-8501 Japan
   \\
    \scalebox{0.95}[1]{\tt fuchino@diamond.kobe-u.ac.jp,\,andreomr@gmail.com,\,hsakai@people.kobe-u.ac.jp}}
\ifextended

\begin{quotation}
	\footnotesize
	\noindent
	\centerline{\normalsize\tt Contents\hspace{6em}\mbox{}}\mbox{}\\
       {\mbox{}\hspace{-1.6em}\tt\makebox[3.4ex][l]{\ref{Einf}.}%
         Introduction}\ \ \dotfill\ \ \pageref{Einf}\\ 
       {\mbox{}\hspace{-1.6em}\tt\makebox[3.4ex][l]{\ref{lt-conti}.}%
         Reflection down to $\LT\continuum$}\ \ \dotfill\ \ \pageref{lt-conti}\\ 
       {\mbox{}\hspace{-1.6em}\tt\makebox[3.4ex][l]{\ref{internal}.}%
         Internal interpretation of stationary logic}\ \ \dotfill\ \ \pageref{internal}\\ 
       {\mbox{}\hspace{-1.6em}\tt\makebox[3.4ex][l]{\ref{PKL}.}%
         Stationarity quantifier in PKL logic}\ \ \dotfill\ \ \pageref{PKL}\\ 
       {\mbox{}\hspace{-1.6em}\tt\makebox[3.4ex][l]{\ref{laver}.}%
         Laver-generic large cardinals}\ \ \dotfill\ \ \pageref{laver}\\ 



       \noindent
       {\mbox{}\hspace{-1.6em}\tt
         References}\ \ \dotfill\ \ \pageref{ref}\\ 

\end{quotation}

\fi
\renewcommand{\thefootnote}{}
\footnotetext{{\it Date:} February 9, 2018
  \qquad {\it Last update:} 
  \today\ (\now\ JST)\vspace{-1\smallskipamount}
}
\footnotetext{\vspace{-1\smallskipamount}}
\footnotetext{The first author would like to thank Joel Hamkins for several 
  valuable comments and suggestions. }

\ifextended
\footnotetext{\tt This is an extended version of the paper with the same title.  
\\
  All additional 
  details not to be contained in the submitted version of the paper are either typeset in 
  typewriter font (the font this paragraph is typeset) or put in separate appendices. 
  The numbering of the assertions is kept identical with the submitted version. The most 
  up-to-date file of this extended version \ifarxived with correct typesetting of 
  East Asian fonts \fi can be downloaded as:\quad \href{https://fuchino.ddo.jp/papers/SDLS-II-x.pdf}{\tt https://fuchino.ddo.jp/papers/SDLS-II-x.pdf}}
\else
\footnotetext{An updated and extended version of this paper with more details and 
  proofs is downloadable as:\qquad \href{https://fuchino.ddo.jp/papers/SDLS-II-x.pdf}{\tt https://fuchino.ddo.jp/papers/SDLS-II-x.pdf}}
\fi
\begin{abstract}
  Continuing \cite{I}, we study the Strong Downward L\"owenheim-Skolem Theorems (SDLSs)
  of the stationary logic and their variations. 

  In \cite{I}, it has been shown that the SDLS for the ordinary stationary logic 
  with weak second-order parameters $\SDLS(\calL^{\aleph_0}_{stat},\LT\aleph_2)$ 
  down to $\LT\aleph_2$ is equivalent to the conjunction of \CH\ and Cox's 
  Diagonal Reflection Principle for internally clubness. 

  We show that the SDLS for the stationary logic without weak second-order 
  parameters $\SDLS^-(\calL^{\aleph_0}_{stat},\LT\continuum)$ down to 
  $\LT\continuum$ implies that the size of the continuum is $\aleph_2$.
  In contrast, an internal interpretation of the stationary logic  
  can satisfy the SDLS down 
  to $\LT\continuum$ under the continuum being of size $>\aleph_2$. This SDLS is shown to be 
  equivalent to an internal version of the Diagonal Reflection  
  Principle down to an internally stationary set of size $\LT2^{\aleph_0}$. 

  We also consider a $\Pkl{}{}$ version of the stationary logic and show that the 
  SDLS for this logic in internal interpretation $\SDLS^{\intnl}_+(\calL^{PKL}_{stat},\LT\continuum)$
  for reflection down to $\LT2^{\aleph_0}$ is consistent under the 
  assumption of the consistency of \ZFC\ $+$ ``the existence of a supercompact 
  cardinal'' and this SDLS implies that the continuum is (at least) weakly Mahlo.  

  These three ``axioms'' in terms of SDLS are consequences of three 
  instances of a 
  strengthening of generic supercompactness which we call Laver-generic
  supercompactness. Existence of a Laver-generic supercompact cardinal in each of these 
  three instances also fixes the cardinality of the continuum to be 
  $\aleph_1$ or $\aleph_2$ or very large respectively. 
  We also show that the existence of one of these generic large 
  cardinals implies the ``$++$'' version of the corresponding forcing axiom.  
\end{abstract}
\renewcommand{\thefootnote}{\arabic{footnote})\,}
\section{Introduction}
\Label{Einf}
We use the notation and conventions set up in \cite{I}: We assume that all structures and 
languages have at most countable signature. 
$\calL^{\aleph_0}$ is the weak second 
order logic extending the usual first-order logic with monadic second-order 
variables with the interpretation that they  
run over countable subsets of the underlying set of the structure and also with 
the built-in binary predicate symbol $\varin$ with which  
we can build a new type of atomic formulas of the form $x\varin X$ where $x$ is 
a first order and $X$ a weak second-order variables. The interpretation of the 
atomic formula $x\varin X$ for a structure $\gmA=\pairof{A\ctenten}$ with
$a\in A$ and $U\in[A]^{\aleph_0}$ is, as expected,  
\begin{xitemize}
\xitem[Einf-0] 
  $\gmA\models x\varin X\,(a,U)$\ \ $\Leftrightarrow$ $a\in U$.
\end{xitemize}

$\calL^{\aleph_0,II}$ is then the logic obtained from $\calL^{\aleph_0}$ by adding the weak second order 
existential quantifier $\exists X$  (and its dual $\forall X$) where for a formula
$\varphi(x_0\ctentenc X_0\ctentenc X)$ in $\calL^{\aleph_0,II}$ and a structure
$\gmA=\pairof{A\ctenten}$ with $a_0\ctenten\in A$ and
$U_0\ctenten\in[A]^{\aleph_0}$, 
\begin{xitemize}
\xitem[Einf-1] 
  $\gmA\models\exists X\,\varphi(a_0\ctentenc U_0\ctentenc X)$\\ 
  $\Leftrightarrow$\ \
  there is $U\in[A]^{\aleph_0}$ \st\
  $\gmA\models\varphi(a_0\ctentenc U_0\ctentenc U)$. 
\end{xitemize}
$\calL^{\aleph_0}_{stat}$ is the logic obtained form $\calL^{\aleph_0}$ by adding 
the new weak second-order quantifier $stat\,X$ (and its dual $aa\,X$) where 
for a formula
$\varphi(x_0\ctentenc X_0\ctentenc X)$ in $\calL^{\aleph_0}_{stat}$ and a structure
$\gmA=\pairof{A\ctenten}$ with $a_0\ctenten\in A$ and
$U_0\ctenten\in[A]^{\aleph_0}$, 
\begin{xitemize}
\xitem[Einf-2] 
  $\gmA\models stat\,X\,\varphi(a_0\ctentenc U_0\ctentenc X)$\\ 
  $\Leftrightarrow$\ \
  $\setof{U\in[A]^{\aleph_0}}{\gmA\models\varphi(a_0\ctentenc U_0\ctentenc U)}$ 
  is stationary in $[A]^{\aleph_0}$. 
\end{xitemize}

Finally $\calL^{\aleph_0,II}_{stat}$ is the logic obtained 
from $\calL^{\aleph_0}$ by adding both types of the weak second-order quantifiers.

For one of the logics $\calL$ introduced above, 
and structures $\gmA$, $\gmB$ of the same signature 
with $\gmB\subseteq\gmA$. We say $\gmB$ is $\calL$-elementary submodel of $\gmA$ 
(notation: $\gmB\prec_\calL\gmA$)
if, for any formula $\varphi(x_0\ctentenc X_0\ctenten)$ in $\calL$ of the signature
where $x_0\ctenten$ are first order and $X_0\ctenten$ weak second-order 
variables, for any $b_0\ctenten\in\ulsetof{\gmB}$ and for any countable subsets
$B_0\ctenten$  of $\ulsetof{\gmB}$, we have 
\begin{xitemize}
\xitem[Einf-3] 
  $\gmB\models\varphi(b_0\ctenten,B_0\ctenten)$ holds if and only if $\gmA\models\varphi(b_0\ctenten,B_0\ctenten)$.
\end{xitemize}

$\gmB$ is a weakly $\calL$-elementary submodel of $\gmA$ (notation:
$\gmB\prec^-_{\calL}\gmA$), 
if  
\begin{xitemize}
\xitem[Einf-4] 
  $\gmB\models\varphi(\variables{b}{n-1})$ holds if and only if
  $\gmA\models\varphi(\variables{b}{n-1})$ holds 
\end{xitemize}
for all formulas $\varphi=\varphi(x_0\ctenten)$ in $\calL$ without free weak second-order variables,
and for all $\variables{b}{n-1}\in\ulsetof{B}$. 

The 
Strong Downward L\"owenheim-Skolem Theorem \footnote{The  
  adjective ``strong'' is added to indicate that $\gmB$ in the statement of the 
  property is not merely elementarily equivalent to but also elementary submodel of
  $\gmA$.} (abbreviated in the following as SDLS) for (elementary substructures \wrt\ the formulas of a language) $\calL$ down to 
$\LT\kappa$ is the assertion defined by

\begin{xitemize}
\item[$\SDLS(\calL,\LT\kappa)$:] {\it For any 
  structure $\gmA$ of 
  countable signature there is $\gmB\prec_\calL\gmA$ of cardinality $\LT\kappa$.}
\end{xitemize}

We also consider the SDLS \wrt\ the weak
$\calL$-elementary submodel relation: 

\begin{xitemize}
\item[$\SDLS^-(\calL,\LT\kappa)$:] {\it For any 
  structure $\gmA$ of 
  countable signature, there is $\gmB\prec^-_{\calL}\gmA$ of cardinality $\LT\kappa$.}
\end{xitemize}

We shall call the cardinal $\kappa$ as above 
the 
{\it reflection cardinal} or {\it L\"owenheim-Skolem cardinal}
of the respective SDLS.

The SDLSs formulated as above have the 
following natural strengthenings:

\begin{xitemize}
\item[$\SDLS_+(\calL,\LT\kappa)$:] {\it For any structure $\gmA=\pairof{A\ctenten}$ of 
  countable signature with $\cardof{A}\geq\kappa$, there are {\it stationarily many} 
  $M\in[A]^{\LT\kappa}$ \st\ $\gmA\restr M\prec_\calL\gmA$. }
\item[$\SDLS^-_+(\calL,\LT\kappa)$:] {\it For any structure $\gmA=\pairof{A\ctenten}$ of 
  countable signature with $\cardof{A}\geq\kappa$, there are {\it stationarily many} 
  $M\in[A]^{\LT\kappa}$ \st\ $\gmA\restr M\prec^-_\calL\gmA$. }
\end{xitemize}

The SDLS theorems for the logics introduced above 
can be characterized by some (combinations of) known principles:

\begin{Thm}{\rm (\Thmof{main-thm},\,\assertof{1}, \Propof{P-sdls-1}, 
    \Lemmaof{P-DRP-3},\,\assertof{1}, \Corof{P-DRP-4} in 
    \cite{I})} \Label{P-Einf-0}
  Suppose that $\kappa$ is a 
  regular cardinal $\geq\aleph_2$. \smallskip
  
  \assert{1} $\SDLS^-(\calL^{\aleph_0},\LT\kappa)$ is a theorem in \ZFC.\smallskip

  \assert{2} Each of $\SDLS_+(\calL^{\aleph_0},\LT\kappa)$,
  $\SDLS^-_+(\calL^{\aleph_0,II},\LT\kappa)$ and 
  $\SDLS_+(\calL^{\aleph_0,II},\LT\kappa)$ is equivalent 
  to $\mu^{\aleph_0}<\kappa$ for all $\mu<\kappa$.\smallskip
  
  \assert{3} $\SDLS^-_+(\calL^{\aleph_0}_{stat},\LT\kappa)$ is equivalent to
  $\DRP(\LT\kappa, \IC_{\aleph_0})$.\smallskip

  \assert{4} Each of 
  $2^{\aleph_0}<\kappa$ $+$ $\SDLS^-_+(\calL^{\aleph_0}_{stat},\LT\kappa)$,
  $\SDLS^-_+(\calL^{\aleph_0,II}_{stat},\LT\kappa)$,
  $\SDLS_+(\calL^{\aleph_0}_{stat},\LT\kappa)$, 
  $\SDLS_+(\calL^{\aleph_0,II}_{stat},\LT\kappa)$ and $2^{\aleph_0}<\kappa$ $+$
  $\DRP(\LT\kappa,\IC_{\aleph_0})$ is equivalent to $\mu^{\aleph_0}<\kappa$ for all $\mu<\kappa$ $+$
  $\DRP(\LT\kappa,\IC_{\aleph_0})$. \qed
\end{Thm}

Note that the parameter ``$\LT\kappa$'' in the SDLS 
statements for the logics $\calL$ as above (except for 
$\SDLS^-(\calL^{\aleph_0},\LT\kappa)$) is impossible for $\kappa\leq\aleph_1$. 
In case of $\SDLS^-(\calL^{\aleph_0}_{stat},\kappa)$, this can be seen in the fact 
that the first order quantifier $Qx\,\varphi$ which states ``there are uncountably 
many $x$ with $\varphi$'' is expressible with the stationarity quantifier as
$stat\,X\,\exists x(x\notvarin X\land\varphi)$. 

For $\kappa=\aleph_2$, the $+$-version of the strong downward L\"owenheim-Skolem 
statements are equivalent with the corresponding statements without $+$:

\begin{Lemma}{\rm (\Lemmaof{P-sdls-0} in \cite{I})} \Label{P-Einf-1}
  Suppose that $\calL$ is one of the four logics as above. 
  Then \smallskip

  \assert{1} $\SDLS_+(\calL,\LT\aleph_2)$ and 
    $\SDLS(\calL,\LT\aleph_2)$ are equivalent and\smallskip

  \assert{2}
    $\SDLS^-_+(\calL,\LT\aleph_2)$ and 
    $\SDLS^-(\calL,\LT\aleph_2)$ are equivalent.\qed
\end{Lemma}

The following is immediate from \Thmof{P-Einf-0} and \Lemmaof{P-Einf-1}.
\begin{Cor}{\rm (\Thmof{main-thm} in \cite{I})}
  \Label{P-Einf-2}
  \wassertof{1} $\SDLS^-(\calL^{\aleph_0},\LT\aleph_2)$ is a theorem in \ZFC.\smallskip

  \assert{2} Each of $\SDLS(\calL^{\aleph_0},\LT\aleph_2)$,
  $\SDLS^-(\calL^{\aleph_0,II},\LT\aleph_2)$ and 
  $\SDLS(\calL^{\aleph_0,II},\LT\aleph_2)$ is equivalent to \CH.\smallskip

  \assert{3} $\SDLS^-(\calL^{\aleph_0}_{stat},\LT\aleph_2)$ is equivalent to
  $\DRP(\IC_{\aleph_0})$.\smallskip

  \assert{4} Each of 
  $\CH$ $+$ $\SDLS^-(\calL^{\aleph_0}_{stat},\LT\aleph_2)$,
  $\SDLS^-(\calL^{\aleph_0,II}_{stat},\LT\aleph_2)$,
  $\SDLS(\calL^{\aleph_0}_{stat},\LT\aleph_2)$ and
  $\SDLS(\calL^{\aleph_0,II}_{stat},\LT\aleph_2)$ is equivalent to $\CH$ $+$
  $\DRP(\IC_{\aleph_0})$. \qed
\end{Cor}

$\DRP(\LT\aleph_2, \IC_{\aleph_0})$ is the Diagonal Reflection Principle for 
internal clubness introduced in Cox \cite{cox:} and 
$\DRP(\LT\kappa, \IC_{\aleph_0})$ is a generalization of it introduced in 
\cite{I}. We will not repeat the definition of $\DRP(\LT\kappa, \IC_{\aleph_0})$.  
Instead, we cite the following combinatorial characterization of this 
principle given in \cite{I}:

For a cardinal $\mu$, let
\begin{xitemize}
\xitem[Einf-5] $\IU_\mu=\setof{X}{[X]^\mu\cap X\mbox{ is cofinal in }[X]^\mu}$;
\xitem[Einf-6] $\IS_\mu=\setof{X}{[X]^\mu\cap X\mbox{ is stationary }[X]^\mu}$;
\xitem[Einf-7] $\IC_\mu=\setof{X}{[X]^\mu\cap X\mbox{ contains a subset which is club in }[X]^\mu}$.
\end{xitemize}
Elements of $\IU_\mu$, $\IS_\mu$, $\IC_\mu$ are said to be {\it internally 
  unbounded}, {\it internally stationary} and {\it internally club} (\wrt\ subsets 
of cardinality $\mu$) respectively. 

Let $\calC$ be one of  $\IU_{\aleph_0}$, $\IS_{\aleph_0}$,
$\IC_{\aleph_0}$. For $\calC$ and regular cardinals $\kappa$, $\lambda$ with
$\kappa\leq\lambda$, 
\begin{xitemize}
\item[$(\ast)^{\calC}_{\LT\kappa,\lambda}$:] {\it For any countable expansion 
$\tilde{\gmA}$ of $\pairof{\calH(\lambda),\in}$ and sequence 
$\seqof{S_a}{a\in\calH(\lambda)}$ \st\ $S_a$ is a stationary subset of 
$[\calH(\lambda)]^{\aleph_0}$ for all $a\in\calH(\lambda)$, there is an
  $M\in[\calH(\lambda)]^{\LT\kappa}$ \st\ 
  \begin{xitemize}
  \item[\wassertof{1}] $M\in\calC$;
  \item[\wassertof{2}] $\tilde{\gmA}\restr M\prec\tilde{\gmA}$ and 
  \item[\wassertof{3}] $S_a\cap[M]^{\aleph_0}$ is stationary in $[M]^{\aleph_0}$ 
    for all $a\in M$.
  \end{xitemize}}
\end{xitemize}

Similarly to $\SDLS_+(\cdots)$ and $\SDLS^-_+(\cdots)$ we can also define the 
following strengthening of the principle $(\ast)^{\calC}_{\LT\kappa,\lambda}$ above 
as
\begin{xitemize}
\item[$(\ast)^{+\,\calC}_{\LT\kappa,\lambda}$:] {\it For any countable expansion 
$\tilde{\gmA}$ of $\pairof{\calH(\lambda),\in}$ and sequence 
$\seqof{S_a}{a\in\calH(\lambda)}$ \st\ $S_a$ is a stationary subset of 
$[\calH(\lambda)]^{\aleph_0}$ for all $a\in\calH(\lambda)$, there are 
  {\it stationarily many} $M\in[\calH(\lambda)]^{\LT\kappa}$ \st\ 
  \begin{xitemize}
  \item[\wassertof{1}\phantom{$'$}] $M\in\calC$;
  \item[\wassertof{2}\phantom{$'$}] $\tilde{\gmA}\restr M\prec\tilde{\gmA}$ and 
  \item[\wassertof{3}\phantom{$'$}] $S_a\cap [M]^{\aleph_0}$ is stationary in $[M]^{\aleph_0}$ 
    for all $a\in M$.\memo{(3) corrected.}
  \end{xitemize}}
\end{xitemize}

We also have the equivalence of this strengthening with the original
$(\ast)^{\calC}_{\LT\kappa,\lambda}$ in case $\kappa=\aleph_2$ (\Lemmaof{P-DRP-1} in \cite{I}). The 
generalizations of Cox's Diagonal Reflection Principle is characterized as 
the global version of $(\ast)^{+\,\calC}_{\LT\kappa,\lambda}$:

\begin{Thm}{\rm (\Lemmaof{L-sdls-0-1} in \cite{I})}
  \Label{P-Einf-3}
  Suppose that $\calC$ is one of\/ $\IU_{\omega_0}$, $\IS_{\omega_0}$,
  $\IC_{\omega_0}$. For a regular cardinal
  $\kappa>\aleph_1$,    $\DRP(\LT\kappa,\calC)$ holds if and only if  
  \begin{xitemize}
  \xitem[Einf-8] 
    $(\ast)^{+\,\calC}_{\LT\kappa,\lambda}$ holds for all 
    cardinal $\lambda\geq\kappa$. \qed
  \end{xitemize}  

\end{Thm}

Below, we shall review some basic facts about the forcing constructions which are 
used in later sections.  

For a \po\ $\poP$ we denote with 
$\ro(\poP)$ the \cBa\ $\baA$ \st\ the separative quotient of $\poP$ can be 
densely embedded into $\baA^+$. 
For \pos\ $\poP$ and $\poQ$ we write $\poP\sim\poQ$ if $\ro(\poP)\cong\ro(\poQ)$. 
For a cardinal $\kappa$ and a set $X$, $\Col(\kappa,S)$ denotes the \po\ as is 
defined in Kanamori \cite{higher-inf:}. 

For \pos\ $\poP$ and $\poQ$, we write $\poP\circleq\poQ$ if 
$\ro(\poP)$ can be completely embeddable into $\ro(\poQ)$.

The following \Thmof{Th-col-0} is a generalization of Proposition 10.20 in Kanamori \cite{higher-inf:}.
It can be proved similarly to the Proposition. 
\begin{Thm}
  \Label{Th-col-0}
  Suppose that $\kappa$ is regular and $\kappa<\lambda$. If $\poP$ is a 
  separative \po\ \st\ $\cardof{\poP}=\lambda$, $\poP$ is $\kappa$-closed and
  \begin{xitemize}
  \xitem[col-1] 
    $\forces{\poP}{\mbox{there is a surjection }\kappa\rightarrow\lambda}$, 
  \end{xitemize}
  then 
  $\ro(\poP)\cong\ro(\Col(\kappa,\ssetof{\lambda}))$.
  \ifextended\else\qed\fi
\end{Thm}
\ifextended{\tt
\prf 
Note that \xitemof{col-1} implies that $\poP$ is atomless. 
\Wolog, we assume that $\poP$ is a dense sub-ordering of $\ro(\poP)^+$. 
\begin{Claim}
  \Label{cl-col-0}
  For any $r\in\ro(\poP)^+$ there is a pairwise incompatible
  $D\subseteq\poP\downarrow r$ of cardinality $\lambda$. 
\end{Claim}
\prfofClaim By the assumption \xitemof{col-1}, $\ro(\poP)^+\downarrow r$ does not 
have the $\lambda$-cc. Hence there is a pairwise disjoint
$D'\subseteq\ro(\poP)^+\downarrow r$ of cardinality $\lambda$. Since $\poP\downarrow r$ is dense 
in $\ro(\poP)^+\downarrow r$, $D'$ has a refinement $D\subseteq\poP\downarrow r$ of 
cardinality $\geq\lambda$. On the other hand, we have $\cardof{D}\leq\lambda$ since $\cardof{\poP}=\lambda$. 
\qedofClaim\qedskip

Let $\utildeg$ be a $\poP$-name \st\
\begin{xitemize}
\xitemA[col-1-0] 
  $\forces{\poP}{\mapping{\utildeg}{\kappa}{\utildegenG}\mbox{ is a bijection}}$ 
\end{xitemize}
where 
$\utildegenG$ is the standard $\poP$-name of a $(\uniV,\poP)$-generic set. 
Note that, for a $(\uniV,\poP)$-generic $\genG$, we have $\kappa\leq\cardof{\genG}$ 
  by $\kappa$-closedness of $\poP$ and $\cardof{\genG}\leq\cardof{\lambda}=\kappa$ by 
$\cardof{\poP}=\lambda$ in $\uniV$. Hence we have $\cardof{\genG}=\kappa$ in $\uniV[\calG]$.

Let
\begin{xitemize}
\xitemA[] 
  $D=\setof{p\in\Col(\kappa,\ssetof{\lambda})}{\dom(p)=\ssetof{\lambda}\times\alpha
  \mbox{ for some }\alpha<\kappa}$. 
\end{xitemize}
Clearly $D$ is a dense subset of $\poP$. It is enough to show that there is an 
order and incompatibility preserving $\mapping{e}{D}{\ro(\poP)^+}$ \st\
$e\imageof D$ is dense in $\ro(\poP)^+$. 

We define such $e$ by induction on $\ell(p)<\kappa$ for $p\in D$ where $\ell(p)$ 
is defined to be the ordinal with 
$\dom(p)=\ssetof{\lambda}\times\ell(p)$. 

Let $e(\emptyset)=\bbone_\poP$.

Having defined $e(p)$ for $p\in D$ with $\ell(p)=\alpha$, let 
$\setof{a^p_\xi}{\xi<\lambda}$ be a maximal antichain in $\poP\downarrow e(p)$ \st, 
for each $\xi<\lambda$, $a^p_\xi\forces{\poP}{\utildeg(\alpha)=\check{r}}$ for some
$r\in\poP$. Note that the construction of such $a^p_\xi$'s is possible by 
\Claimof{cl-col-0}.
Set $e(p\cup\ssetof{\pairof{\pairof{\lambda,\alpha},\xi}})=a^p_\xi$. 
For a limit $\delta<\kappa$ and $p\in D$ with $\ell(p)=\delta$, suppose that
$e(p\restr\beta)$, $\beta<\delta$ has been defined \st\
$\seqof{e(p\restr\beta)}{\beta<\delta}$ is a decreasing sequence. 
Let $e(p)=\prod_{\beta<\delta}e(p\restr \beta)$. Note that 
$e(p)>\bbzero_{\ro(\poP)}$ by $\kappa$-closedness of $\poP$ and since the elements of 
$\poP$ are cofinal in the sequence $\seqof{e(p\restr\beta)}{\beta<\delta}$. 

By the induction on $\alpha<\kappa$ we can show that 
$A_\alpha=\setof{e(p)}{p\in D,\,\ell(p)=\alpha+1}$ is a maximal antichain. 
It is also clear that $e$ is order and 
incompati\-bility preserving. Thus the following claim implies the 
forcing equivalence of $\poP$ and $\Col(\kappa,\ssetof{\lambda})$. 
\begin{Claim}
  $e\imageof D$ is dense in $\ro(\poP)^+$. 
\end{Claim}
\prfofClaim
For an arbitrary $r\in\poP$ we have $r\forces{\poP}{\check{r}\in\utildecalG}$. Hence 
there is $r'\leq_\poP r$ and $\alpha<\lambda$ \st\
$r'\forces{\poP}{\utildeg(\check{\alpha})=\check{r}}$.  Since $A_\alpha$ is a 
maximal antichain in $\poP$, there is a $p\in D$ with $\ell(p)=\alpha+1$ \st\ 
$e(p)=a^{p\restr\alpha}_{p(\alpha)}$ is compatible with $r'$. Since
$e(p)=a^{p\restr\alpha}_{p(\alpha)}$ decides $\utildeg(\check{\alpha})$, we have
$e(p)\forces{\poP}{\utildeg(\check{\alpha})=\check{r}}$. Thus 
$e(p)\forces{\poP}{\check{r}\in\utildecalG}$ by \xitemof{col-1-0}. 
Since $\poP$ is separative 
it follows that $e(p)\leq_\poP r$. 
\qedofClaim\\
\qedofThm
}\fi

\begin{Cor}
  \Label{Cor-col-0}
  \wassertof{1}
  Suppose that $\lambda^{<\kappa}=\lambda$. Then we have
  \begin{xitemize}
  \xitem[col-2] 
    $\Col(\kappa,\ssetof{\lambda})\sim\Col(\kappa,S)\sim\Fn(\kappa,\lambda,<\kappa)$ 
  \end{xitemize}
  for all $S\subseteq\lambda^+$ with $\lambda+1\leq\sup S<\lambda^+$. \smallskip
  
  \assert{2} For any $\kappa$-closed $\poP$ with $\cardof{\poP}\leq\mu=\mu^{\LT\kappa}<\lambda$ we have 
  $\Col(\kappa,\lambda)\sim\poP\times\Col(\kappa,\lambda)\sim\poP*\Col(\kappa,\lambda)^{\uniV^\poP}$. 
  In particular, we have $\poP\circleq\Col(\kappa,\lambda)$. 
\end{Cor}
\prf
\assertof{1}: 
Let $\poP=\Col(\kappa,S)$ for $S$ as above. By the assumption we have 
$\cardof{\poP}=\lambda$, $\poP$ is separative and $\kappa$-closed. Also $\poP$ 
adds a surjection from $\kappa$ to $\lambda$. 
Hence 
\Thmof{Th-col-0} implies $\Col(\kappa,\ssetof{\lambda})\sim\Col(\kappa,S)$. 
$\Col(\kappa,\ssetof{\lambda})\sim\Fn(\kappa,\lambda,<\kappa)$ also can be shown 
similarly. \smallskip

\assertof{2}: Let $\mu<\lambda$ be \st\ $\cardof{\poP}\leq\mu$ 
and $\mu^{\LT\kappa}=\mu$. Then 

\begin{xitemize}
\xitem[] 
  $
  \begin{array}[t]{ll}
    \poP*\Col(\kappa,\lambda)^\poP\\[\jot]
    {}\sim \poP\times\Col(\kappa,\lambda)
    	&\mbox{; by $\kappa$-closedness of $\poP$}\\[\jot]
    {}\sim (\poP\times \Col(\kappa,\mu+1))\times\Col(\kappa,\lambda\setminus(\mu+1))\\[\jot]
    {}\sim 
    \Col(\kappa,\ssetof{\mu})\times\Col(\kappa,\lambda\setminus(\mu+1))
    	&\mbox{; by \Thmof{Th-col-0}}\\[\jot]
    {}\sim 
    \Col(\kappa,\mu+1)\times\Col(\kappa,\lambda\setminus(\mu+1))
     &\mbox{; by \assertof{1}}\\[\jot]
    {}\sim 
    \Col(\kappa,\lambda).
  \end{array}
  $\\\qedofCor
\end{xitemize}

\section{\texorpdfstring{Reflection down to $\LT\continuum$}{Reflection down 
    to < continuum}}
\Label{lt-conti}
$\SDLS^-(\calL^{\aleph_0}_{stat},\LT\continuum)$ is consistent under
$\continuum=\aleph_2$ (e.g.\ under the assumption of the existence of 
a supercompact cardinal):\ \ $\MA^{+\omega_1}(\sigma\mbox{-closed})$ implies
$\SDLS^-(\calL^{\aleph_0}_{stat},\LT\aleph_2)$ and 
$\MA^{+\omega_1}(\sigma\mbox{-closed})$ is consistent (under the large cardinal assumption) with $\continuum=\aleph_2$.
This is no more the case if the continuum is larger than $\aleph_2$:

\begin{Prop}
  \Label{P-lt-conti-0} 
  $\SDLS^-(\calL^{\aleph_0}_{stat},<\kappa)$ for $\kappa>\aleph_2$ implies
  $\kappa>2^{\aleph_0}$. 
\end{Prop}

The proof of the Proposition uses the following \Thmof{P-lt-conti-1} which is often the main 
ingredient of a proof showing that a certain principle implies
$\continuum\leq\aleph_2$ (see e.g.\ the proof of Theorem 37.18 in 
\cite{millennium-book:}). 

\begin{Thm}{\rm (Theorem 3.2 \assertof{a} in Baumgartner and Taylor \cite{baumgartner-taylor:})}
  \Label{P-lt-conti-1} If $\calC$ is a club subset of $[\omega_2]^{\aleph_0}$, then 
  there is a countable set $A\subseteq\omega_2$ \st\ $\cardof{[A]^{\aleph_0}\cap C}=2^{\aleph_0}$.\\\qed
\end{Thm}

\noindent
\prfof{\Propof{P-lt-conti-0}} $\SDLS^-(\calL^{\aleph_0}_{stat},\LT\aleph_2)$ implies
$2^{\aleph_0}\leq\aleph_2$: it is easy to see that 
$\SDLS^-(\calL^{\aleph_0}_{stat},\LT\aleph_2)$ implies the reflection principle
$\RP(\omega_2)$ in \cite{millennium-book:}. $\RP(\omega_2)$ implies $2^{\aleph_0}\leq\aleph_2$
(a result by Todor\v{c}evi\'c, see Theorem 37.18 in \cite{millennium-book:}). 
We have 
$\kappa>\aleph_2\geq2^{\aleph_0}$.  

Thus we may assume that $\SDLS^-(\calL^{\aleph_0}_{stat},\LT\aleph_2)$ does not 
hold. Hence there is a structure $\gmA$ \st, for any  
$\gmB\prec^-_{(\calL^{\aleph_0}_{stat})}\gmA$, we 
have $\ccardof{\gmB}\geq\aleph_2$. Let $\lambda=\ccardof{\gmA}$. \Wolog, we may 
assume that $\ulsetof{\gmA}=\lambda$. Let 
\begin{xitemize}
\xitem[lt-conti-0] 
  $\gmA^*=\pairof{\calH(\lambda^+),\underbrace{\lambda\ctenten}_{=\gmA},\in}$.
\end{xitemize}
Note that we have
\begin{xitemize}
\xitem[lt-conti-1] 
  $\gmA^*\models
  \underbrace{aa\,X\,\exists x\forall y\,(y\varin X\leftrightarrow y\in x)}_{=\varphi}$
\end{xitemize}
where ``$aa\,X$'' is the dual quantifier to ``$stat\,X$''. That is, we treat 
``$aa\,X\,\psi$'' just as an abbreviation of ``$\neg stat\,X\neg\psi$''. 
Note that
$\gmA\models aa\,X\,\varphi(X\ctenten)$ if and only if there are club many 
$U\in[A]^{\aleph_0}$ with $\gmA\models\varphi(U\ctenten)$.

By $\SDLS^-(\calL^{\aleph_0}_{stat},<\kappa)$, there is 
$M\in[\calH(\lambda^+)]^{<\kappa}$ \st\
$\gmA^*\restr M\prec^-_{(\calL^{\aleph_0}_{stat})}\gmA^*$.
It follows that
$\gmA\restr (\lambda\cap M)\prec^-_{(\calL^{\aleph_0}_{stat})}\gmA$.  
By the choice of $\gmA$, we have $\cardof{M}\geq\cardof{\lambda\cap M}\geq\aleph_2$. 

Since
$\gmA^*\restr M\models\varphi$ 
by \xitemof{lt-conti-1}, and by elementarity, there is
$C\subseteq[M]^{\aleph_0}\cap M$ which is a club in
$[M]^{\aleph_0}$. By \Thmof{P-lt-conti-1}, it 
follows that $\kappa>\cardof{M}\geq\cardof{C}\geq 2^{\aleph_0}$. 
\qedof{\Propof{P-lt-conti-0}}

\begin{Cor}
  \Label{P-lt-conti-2} $\SDLS^-(\calL^{\aleph_0}_{stat},\LT\continuum)$ implies
  $\continuum=\aleph_2$. 
\end{Cor}
\prf $\SDLS^-(\calL^{\aleph_0}_{stat},\LT\continuum)$ implies
$\continuum\leq\aleph_2$ by \Propof{P-lt-conti-0}. $\continuum=\aleph_1$ is 
impossible under $\SDLS^-(\calL^{\aleph_0}_{stat},\LT\continuum)$ since, as we 
have seen already, 
$\SDLS^-(\calL^{\aleph_0}_{stat},\LT\aleph_1)$ does not hold (in \ZFC).\qedofCor\qedskip

Note that $\SDLS^-(\calL^{\aleph_0}_{stat},\LT\continuum)$ holds under
$\MA^{+\omega_1}(\sigma\mbox{-closed})$ $+$ $2^{\aleph_0}=\aleph_2$. However,

\begin{Cor}
  \label{P-lt-conti-2-0}
  $\SDLS(\calL^{\aleph_0}_{stat},\LT\continuum)$ is inconsistent.
\end{Cor}
\prf Assume that $\SDLS(\calL^{\aleph_0}_{stat},\LT\continuum)$ holds. Then
$\SDLS^-(\calL^{\aleph_0}_{stat},\LT\continuum)$ holds and hence
$\continuum=\aleph_2$ by \Corof{P-lt-conti-2}. Thus,
$\SDLS(\calL^{\aleph_0}_{stat},\LT\aleph_2)$ holds. By 
\Lemmaof{P-Einf-1},\,\assertof{1} and \Thmof{P-Einf-0},\,\assertof{4}, it follows 
that $\continuum=\aleph_1$. This is a contradiction.\qedofCor\qedskip

Translated in terms of diagonal reflection, \Propof{P-lt-conti-0} can also be 
reformulated as:
\begin{xitemize}
\xitem[lt-conti-2] 
  If $\kappa>\aleph_2$, $(\ast)^{+\IC_{\aleph_0}}_{\LT\kappa,\lambda}$ for all 
  cardinal $\lambda\geq\kappa$ implies $\continuum<\kappa$. 
\end{xitemize}

However, the following internal variation of $(\ast)^{+\IS_{\aleph_0}}_{\LT\kappa,\lambda}$ is 
compatible with $\continuum\geq\kappa$ (see \Thmof{P-lt-conti-3}). 
For regular cardinals $\kappa$, $\lambda$ with
$\kappa\leq\lambda$, let 
\begin{xitemize}
\item[$(\ast)^{\intnl+}_{\LT\kappa,\lambda}$:] {\it For any countable expansion 
$\tilde{\gmA}$ of $\pairof{\calH(\lambda),\in}$ and sequence 
$\seqof{S_a}{a\in\calH(\lambda)}$ \st\ $S_a$ is a stationary subset of 
$[\calH(\lambda)]^{\aleph_0}$ for all $a\in\calH(\lambda)$, there are 
  stationarily many 
  $M\in[\calH(\lambda)]^{\LT\kappa}$ \st\ 
  \begin{xitemize}
  \item[\wassertof{2}] $\tilde{\gmA}\restr M\prec\tilde{\gmA}$; and 
  \item[\wassertof{3$'$}] $S_a\cap M$ is stationary in $[M]^{\aleph_0}$ 
    for all $a\in M$.
  \end{xitemize}}
\end{xitemize}
Note that \assertof{2} implies that $a\subseteq M$ for all
$a\in [M]^{\aleph_0}\cap M$. 

By the definition of $(\ast)^{\intnl+}_{\LT\kappa,\lambda}$ it is easy to see that we have
\begin{xitemize}
\xitem[] 
  $(\ast)^{+\IC_{\aleph_o}}_{\LT\kappa,\lambda}$\ \ $\Rightarrow$\ \ 
  $(\ast)^{\intnl+}_{\LT\kappa,\lambda}$\ \ $\Rightarrow$\ \ 
  $(\ast)^{+\IS_{\aleph_0}}_{\LT\kappa,\lambda}$
\end{xitemize}
for any regular $\kappa$, $\lambda$ with $\kappa\leq\lambda$.

For a class $\calP$ of \pos, a cardinal $\kappa$ is said to be {\it generically 
  supercompact by $\calP$}
if, for any regular $\lambda\geq\kappa$, there is a \po\ $\poP\in\calP$ \st, for 
a $(\uniV,\poP)$-generic $\genG$, there are transitive $M\subseteq\uniV[\genG]$ and 
$j\subseteq\uniV[\genG]$ \st\ 
\begin{xitemize}
\xitem[lt-conti-2-0] 
  $\elembed{j}{\uniV}{M}$,
\xitem[lt-conti-2-1] $\crit(j)=\kappa$, $j(\kappa)>\lambda$,
\xitem[lt-conti-2-2] $j\imageof\lambda\in M$. 
\end{xitemize}

Similarly, a cardinal $\kappa$ is said to be {\it generically measurable by a forcing
$\poP$} if, for a $(\uniV,\poP)$-generic $\genG$, there are transitive 
$M\subseteq\uniV[\genG]$  and $\elembed{j}{\uniV}{M}$ \st\ $\crit(j)=\kappa$ and 
$j\imageof\kappa\in M$. 

Note that the condition \xitemof{lt-conti-2-2} is weaker than the closure property of 
$M$ assumed in the usual definition of supercompactness. On the other hand, this 
condition is often enough to prove strong reflection properties of $\kappa$ (see 
the following \Lemmaof{L-lt-conti-0}). 

By definition, if $\kappa$ is generically supercompact by $\calP$ then $\kappa$ 
is generically measurable by some $\poP\in\calP$. 

The following is easy to see:
\begin{Lemma}
  \Label{L-lt-conti-0}
  Suppose that $\genG$ is a $(\uniV,\poP)$-generic filter for a
  \po\/ $\poP\in\uniV$ and $\elembed{j}{\uniV}{M\subseteq\uniV[\genG]}$ \st, for 
  cardinals $\kappa$, $\lambda$ in $\uniV$ 
  with $\kappa\leq\lambda$, $\crit(j)=\kappa$ and $j\imageof\lambda\in M$. 

  \assert{1} For any set $A\in\uniV$  
  with $\uniV\models\cardof{A}\leq\lambda$, we have $j\imageof A\in M$. \smallskip

  \assert{2} $j\restr\lambda$, $j\restr\lambda^2\in M$.\smallskip

  \assert{3} For any $A\in\uniV$ with $A\subseteq\lambda$ or $A\subseteq\lambda^2$ 
  we have $A\in M$.\smallskip

  \assert{4} $(\lambda^+)^M\geq(\lambda^+)^\uniV$, Thus, if
  $(\lambda^+)^\uniV=(\lambda^+)^{\uniV[\genG]}$,  then
  $(\lambda^+)^M=(\lambda^+)^\uniV$. \smallskip

  \assert{5} $\calH(\lambda^+)^\uniV\subseteq M$.\smallskip

  \assert{6} $j\restr A\in M$ for all $A\in\calH(\lambda^+)^\uniV$. 
\end{Lemma}
\prf \assertof{1}: In $\uniV$, let $\mapping{f}{\lambda}{A}$ be a surjection. 

For each $a\in A$ with $a=f(\alpha)$, we have
\begin{xitemize}
\xitem[pr-L-lt-conti-0] 
  $j(a)=j(f(\alpha))=j(f)(j(\alpha))$
\end{xitemize}
by elementarity. 
Thus $j\imageof A=j(f)\imageof(j\imageof\lambda)$. Since $j(f)$,
$j\imageof\lambda\in M$, it follows that $j\imageof A\in M$.\smallskip

\assertof{2}: Since $j\imageof\lambda\in M$ and $(j\restr\lambda)(\xi)$ for 
$\xi\in\lambda$ is the $\xi$th element of $j\imageof\lambda$, $j\restr\lambda$ is 
definable subset of $\lambda\times j\imageof\lambda$ in $M$ and hence is an 
element of $M$. Similarly, $j\restr\lambda^2\in M$. \smallskip

\assertof{3}: Suppose that $A\in\uniV$ and $A\subseteq\lambda$ (the case of
$A\subseteq\lambda^2$ can be treated similarly). Then $j\imageof A\in M$ by 
\assertof{1}. Thus, by \assertof{2},
$A=(j\restr\lambda)^{-1}\imageof(j\imageof A)\in M$. \smallskip

\assertof{4}: Suppose that $\mu<(\lambda^+)^\uniV$. Then there is $A\in\uniV$ 
with $A\subseteq\lambda^2$ \st\ $A$ codes the order type of $\mu$. $A\in M$ by 
\assertof{3}. Thus $M\modelof{\cardof{\mu}\leq\lambda}$.

If $(\lambda^+)^\uniV=(\lambda^+)^{\uniV[\genG]}$, we have
\begin{xitemize}
\xitem[pr-L-lt-conti-1] 
  $(\lambda^+)^\uniV=(\lambda^+)^{\uniV[\genG]}\geq(\lambda^+)^M\geq(\lambda^+)^\uniV$.
\end{xitemize}

\assertof{5}: For $A\in\calH(\lambda^+)^\uniV$, let $U\in\uniV$ be \st\
$\trcl(A)\subseteq U$ and $\uniV\modelof{\cardof{U}=\lambda}$. Let 
$c_A\subseteq\lambda^2$ and $d_A$, $e_A\subseteq\lambda$ be \st\ $c_A$, $d_A$, $e_A\in\uniV$ 
and 
\begin{xitemize}
\xitem[pr-L-lt-conti-2] 
  $\pairof{\lambda, c_A, d_A, e_A}\cong\pairof{U,\in\restr U^2, \trcl(A), A}$. 
\end{xitemize}
By \assertof{3}, $c_A$, $d_A$, $e_A\in M$ and hence 
$\pairof{\lambda, c_A, d_A, e_A}\in M$. Since $\trcl(A)$ and 
then $A$ can be recovered 
from this quadruplet in $M$, it follows that $A\in M$. 
\smallskip

\assertof{6}: Suppose that $A\in\calH(\lambda^+)^\uniV$. Since $A\in M$ by \assertof{5}, 
it is enough to show 
that $j\restr\trcl(A)\in M$. 

We have $\trcl(A)\in\calH(\lambda^+)^\uniV$ and hence $\trcl(A)\in M$ by 
\assertof{5}. Thus $j\imageof\trcl(A)$, $j\imageof(\in\restr\trcl(A))\in M$ by \assertof{1}. 
But then the mapping $(j\restr\trcl(A))^{-1}$ is the Mostowski collapse of $j\imageof\trcl(A)$. 
Thus $j\restr\trcl(A)\in M$. 
\qedofLemma
\qedskip

The following Lemmas \ref{L-lt-conti-1-0} and \ref{L-lt-conti-1-1} should be well-known and easy to check. 
\Lemmaof{L-lt-conti-1-1} can be proved similarly to 
\Propof{T-laver-1-1},\,\assertof{2}. 

\begin{Lemma}
  \Label{L-lt-conti-1-0}
  If $\kappa$ is generically measurable for some \po\/ $\poP$, then $\kappa$ is regular.
  \ifextended\else\qed\fi
\end{Lemma}
\ifextended{\tt
\prf Suppose that $\kappa=\sup\setof{\kappa_\xi}{\xi<\mu}$ for some $\mu<\kappa$ 
and $\kappa_\xi<\kappa$ for $\xi<\mu$. Let $S=\setof{\kappa_\xi}{\xi<\mu}$ and 
let $\elembed{j}{\uniV}{M\subseteq\uniV[\genG]}$ be as in the definition of 
generic measurability. Then $j(S)=S$ by $\crit(j)=\kappa$. By elementarity it 
follows that $M\modelof{j(\kappa)=\sup S=\kappa}$. This is a contradiction to
$\crit(j)=\kappa$. 
\qedofLemma
\qedskip
}\fi

\begin{Lemma}
  \Label{L-lt-conti-1-1} 
  \wassertof{1} Suppose that $\kappa$ is generically measurable for a 
  \po\/ $\poP$ and $j$, $M\subseteq\uniV[\genG]$ for a $(\uniV,\poP)$-generic
  $\genG$ \st\ $M$ is an inner model of\/ $\uniV[\genG]$ $\elembed{j}{V}{M}$,
  $\crit(j)=\kappa$. Then, in $V[\genG]$, 
  \begin{xitemize}
  \xitem[lt-conti-2-2-0] 
    $F=\setof{a\in(\psof{\kappa})^\uniV}{\kappa\in j(a)}$
  \end{xitemize}
  is a $\uniV$-normal ultrafilter on (the \Ba) $(\psof{\kappa})^V$.\smallskip

  \assert{2} If $\mu<\kappa$ and $\kappa$ is generically measurable for 
  a $\mu$-cc \po\/ $\poP$ then there is a $\mu$-saturated normal ideal 
  over $\kappa$ (in $\uniV$). In particular, $\kappa$ is $\kappa$-weakly Mahlo.
\end{Lemma}
\ifextended{\tt
\prf \assertof{1}: See the proof of \Claimof{Cl-laver-0}.\smallskip

\assertof{2}:  Suppose that $\poP$ is a $\mu$-cc \po\ and $\genG$ 
a $(\uniV,\poP)$-generic filter with $j$, $M\subseteq\uniV[\genG]$ \st\ $M$ is an 
inner model of $\uniV[\genG]$, $\elembed{j}{\uniV}{M}$ and $\crit(j)=\kappa$. Let 
$F$ be defined as in \xitemof{lt-conti-2-2-0} and let $\utilde{F}$ be 
a $\poP$-name of $F$ \st\ all the properties of $F$ we need below are forced for
$\utilde{F}$ by $\bbone_\poP$. In $\uniV$, let
\begin{xitemize}
\xitemA[lt-conti-2-2-1]
  $I_0=\setof{u\in\psof{\kappa}}{\forces{\poP}{\check{u}\notvarin\utilde{F}}}$. 
\end{xitemize}
It is easy to see that $I_0$ is an ideal on $\psof{\kappa}$. 
We show that $I_0$ is as desired. 

Let $F_0$ be the 
dual filter of $I_0$.
\begin{Claim}
  \Label{Cl-lt-conti-0}
  $I_0$ is a normal ideal.
\end{Claim}
\prfofClaim
It is enough to show that $F_0$ is a normal filter. Suppose that
$u_\alpha\in F_0$, for all $\alpha<\kappa$. This means that we have 
$\forces{\poP}{\check{u}_\alpha\varin\utilde{F}}$ for all $\alpha<\kappa$. Since 
$F$ is a $\uniV$-normal ultrafilter, it follows that\\
$\forces{\poP}{\CHECK{\poP}{\bigtriangleup_{\alpha<\kappa}u_\kappa}\varin\utilde{F}}$.
Thus, $\bigtriangleup_{\alpha<\kappa}u_\kappa\in F_0$.
\qedofClaim

\begin{Claim}
  \Label{Cl-lt-conti-1} For any $a\in\psof{\kappa}$, $a$ is $I_0$-stationary if 
  and only if there is $\condp\in\poP$ 
  \st\ $\condp\forces{\poP}{a\in\utilde{F}}$. 
\end{Claim}
\prfofClaim
If there is $\condp\in\poP$ \st\ $\condp\forces{\poP}{\check{a}\varin\utilde{F}}$, 
then $\notforces{\poP}{\check{a}\notvarin\utilde{F}}$ and hence $a\not\in I_0$. 
Thus $a$ is $I_0$-stationary.

If there is no $\condp\in\poP$ \st\
$\condp\forces{\poP}{\check{a}\varin\utilde{F}}$, then
we have\\
$\forces{\poP}{\check{a}\notvarin\utilde{F}}$. Thus $a_0\in I_0$. 
\qedofClaim

\begin{Claim}
  \Label{Cl-lt-conti-2} $I_0$ is $\mu$-saturated. 
\end{Claim}
\prfofClaim
Suppose that $a_\xi$, $\xi<\mu$ are $I_0$-stationary. By \Claimof{Cl-lt-conti-1}, 
there are $\condp_\xi\in\poP$, $\xi<\mu$ \st\
$\condp_\xi\forces{\poP}{a_\xi\in\utilde{F}}$. By the $\mu$-cc of $\poP$, there 
are $\xi_0<\xi_1<\mu$ \st\ $\condp_{\xi_0}$ and $\condp_{\xi_1}$ are compatible in
$\poP$. Let $\condq\in\poP$ be \st\ $\condq\leq_\poP\condp_{\xi_0}$,
$\condp_{\xi_1}$. Then we have
$\condq\forces{\poP}{a_{\xi_0}\cap a_{\xi_1}\in\utilde{F}}$ and hence
$a_{\xi_0}\cap a_{\xi_1}$ is $I_0$-stationary by \Claimof{Cl-lt-conti-1}. 
\qedofClaim\qedskip

By Proposition 16.8 in Kanamori \cite{higher-inf:}, it 
follows that $\kappa$ is a $\kappa$-weakly Mahlo. 
}
\else
\prf \assertof{1}: Easy. \assertof{2}: Let $F$ be as in \assertof{1} and
let $\utilde{F}\in\uniV$ be a $\poP$-name of $F$. Then, by $\mu$-cc of $\poP$, 
$I_0=\setof{u\in\psof{\kappa}}{\forces{\poP}{\check{u}\notvarin\utilde{F}}}$ is 
a $\mu$-saturated normal ideals. By Proposition 16.8 in \cite{higher-inf:}, it 
follows that $\kappa$ is a $\kappa$-weakly Mahlo. 
\fi
\qedofLemma
\qedskip

\begin{Prop}
  \Label{T-laver-1-1}
  Suppose that $\kappa$ is generically supercompact for a class $\calP$ of \pos\ \st\ all\/ 
  $\poP\in\calP$ are $\mu$-cc for some fixed $\mu\in\Card$. Then\smallskip

  \assert{1} \SCH\ holds above 
  $\max\ssetof{2^{\LT\kappa},\mu}$.\smallskip

  \assert{2} For all regular $\lambda\geq\kappa$, there is a $\mu$-saturated 
  normal filter over $\Pkl{}{}$. 
\end{Prop}
\prf \assertof{1}: The following proof is a slight modification of the proof of Solovay's 
theorem on \SCH\ above a strongly compact cardinal (see Theorem 20.8 and its 
proof in \cite{millennium-book:}). 

Let $\lambda\geq\max\ssetof{2^{\LT\kappa},\mu}$ be a regular cardinal. It is 
enough to show that $\lambda^{\LT\kappa}=\lambda$.

Suppose that $\poP\in\calP$ and $(\uniV,\poP)$-generic filter $\genG$ be \st\ 
there are classes $j$, $M\subseteq\uniV[\genG]$ \st\ 
$\elembed{j}{\uniV}{M}\subseteq\uniV[\genG]$, $\crit(j)=\kappa$,
$j(\kappa)>\lambda$ and $j\imageof\lambda\in M$.

Let
\begin{xitemize}
\xitem[laver-1-6] 
  $\calU=\setof{A\in(\psof{\Pkl{}{}})^\uniV}{j\imageof\lambda\in j(A)}$. 
\end{xitemize}
\ifextended
\else

The following can be proved by standard arguments. 

\fi
\begin{Claim}
  \Label{Cl-laver-0} $\calU$ is a $\uniV$-$\kappa$-complete fine ultrafilter on
  $(\psof{\Pkl{}{}})^\uniV$. (Actually $\calU$ is even a $\uniV$-normal ultrafilter.) \ifextended\else\smallqed\fi
\end{Claim}
\ifextended{\tt 
\prfofClaim
For $A$, $B\in(\psof{\Pkl{}{}})^\uniV$ with $A\cup B=(\Pkl{}{})^\uniV$,
we have
$j(A)\cup j(B)=j(A\cup B)=j((\Pkl{}{})^\uniV)=(\Pkl{j(\kappa)}{j(\lambda)})^M\ni 
j\imageof\lambda$. Hence $M\models j\imageof\lambda\in j(A)$ or
$M\models j\imageof\lambda\in j(A)$. That is, $A\in\calU$ or $B\in\calU$.

For $A$, $B\in(\psof{\Pkl{}{}})^\uniV$, we can also show similarly that,
if $A\in\calU$ and $B\in\calU$ then $A\cap B\in\calU$ and that, if $A\in\calU$ 
and $A\subseteq B$, then $B\in\calU$.

To show the $\uniV$-$\kappa$-completeness, suppose that
$\vec{A}=\seqof{A_\xi}{\xi<\eta}\in\uniV$ be \st\ $\eta<\kappa$ and 
$A_\xi\in\calU$ for all $\xi<\eta$. The last condition means that
$j\imageof\lambda\in j(A_\xi)$ for all $\xi<\eta$ by definition of $\calU$. 
Note that $\seqof{j(A_\xi)}{\xi<\eta}=j(\vec{A})\in M$. Thus
$M\models j\imageof\lambda\in\bigcap_{\xi<\eta}j(A_\eta)=j(\bigcap\setof{A_\xi}{\xi<\mu})$,    
i.e.\ $\bigcap\setof{A_\xi}{\xi<\mu}\in\calU$. 

For any $c\in(\Pkl{}{})^\uniV$, letting
$A=\setof{a\in(\Pkl{}{})^\uniV}{c\subseteq a}$, we have 
$M\models j\imageof\lambda\in j(A)$ since
$M\models j(c)=j\imageof c\subseteq j\imageof\lambda$. Thus $A\in\calU$. This 
shows that $\calU$ is fine. 

To show the $\uniV$-normality of $\calU$, suppose that 
$\vec{U}=\seqof{U_\alpha}{\alpha\in\lambda}\in\uniV$ is a sequence of elements 
of $\calU$. Let $U=\bigtriangleup\vec{U}$. Then we have
\begin{xitemize}
\xitemA[] 
	$\uniV\modelof{(\forall x\in\Pkl{}{})(x\in U\leftrightarrow
  \forall\alpha\in x\,(x\in U_\alpha))}$.	
\end{xitemize}
By elementarity
\begin{xitemize}
\xitemA[] 
	$M\modelof{(\forall x\in\Pkl{j(\kappa)}{j(\lambda)})(x\in j(U)\leftrightarrow
  \forall\alpha\in x\,(x\in j(\vec{U})_\alpha))}$.	
\end{xitemize}

Now, for $x=j\imageof\lambda$, if $\alpha\in j\imageof\lambda$, then there is
$\beta\in\lambda$ \st\ $\alpha=j(\beta)$. Since $U_\beta\in\calU$, we have
$j\imageof\lambda\in j(U_\beta)=j(\vec{U})_\alpha$.

Thus we have $M\modelof{j\imageof\lambda\in j(U)}$ or equivalently $U\in\calU$. 
\qedofClaim\qedskip
}\fi

Let $f^*\in\uniV$ with $\mapping{f^*}{(\Pkl{}{})^\uniV}{\On}$ be \st\
\begin{xitemize}
\xitem[laver-1-6-0] 
  $[f^*]_\calU=\sup j_\calU\imageof\lambda$. 
\end{xitemize}
\begin{Claim}
  \Label{Cl-laver-1} $\setof{a\in(\Pkl{}{})^\uniV}{f^*(a)<\lambda}\in\calU$.
\end{Claim}
\prfofClaim
Let $g\in\uniV$ with $\mapping{g}{(\Pkl{}{})^\uniV}{\lambda}$; $a\mapsto\sup(a)$. 
Then, for $\gamma<\lambda$, 
$j_\calU(\gamma)=[c_\gamma]_\calU\leq[g]_\calU$ for all $\gamma<\lambda$ since 
$\calU$ is fine. Hence
$[f^*]_\calU\leq[g]_\calU$. Thus, 
by \L o\'s's Theorem,
\begin{xitemize}
\xitem[] 
  $\calU\ni \setof{a\in(\Pkl{}{})^\uniV}{f^*(a)\leq\hspace{-1.8em}
  \underbrace{\ g(a)\ }_{\mbox{\tt$\ \ =\sup a<\lambda$}}\hspace{-1.8em}}\subseteq
  \setof{a\in(\Pkl{}{})^\uniV}{f^*(a)<\lambda}$. 
\end{xitemize}\vspace*{-3ex}
\qedofClaim
\qedskip

By the Claim above, we may assume \wolog\ that\\ $\mapping{f^*}{(\Pkl{}{})^\uniV}{\lambda}$. 

In $\uniV[\genG]$, let $\calD\subseteq\psof{\lambda}^\uniV$ be defined by
\begin{xitemize}
\xitem[laver-1-7] 
  $\calD=\setof{X\in\psof{\lambda}^\uniV}{{f^*}^{-1}\imageof X\in\calU}$.\footnote{Here,
a more straightforward proof is possible e.g. by directly defining $\calD$. We 
defined $\calD$ via definition of $\calU$ instead since we need $\calU$ in the 
proof of \assertof{2} in any way. }
\end{xitemize}

\ifextended\else 
Clearly we have
\fi
\begin{Claim}
  \Label{Cl-laver-2}
  $\calD$ is a $\uniV$-$\kappa$-complete uniform ultrafilter on $\psof{\lambda}^\uniV$.\ifextended\else\smallqed\fi
\end{Claim}
\ifextended{\tt
\prfofClaim 
That $\calD$ is an ultrafilter on $\psof{\lambda}^\uniV$ is clear by definition.

Suppose that $\seqof{X_\alpha}{\alpha<\mu}\in\uniV$ is a sequence of elements of
$\calD$ for some $\mu<\kappa$. Then $\seqof{{f^*}^{-1}(X_\alpha)}{\alpha<\mu}\in\uniV$ is 
a sequence of elements of $\calU$.
Hence 
${f^*}^{-1}\imageof(\bigcap_{\alpha<\mu}X_\alpha)=\bigcap_{\alpha<\mu}{f^*}^{-1}\imageof X_\alpha\in\calU$. 
Thus $\bigcap_{\alpha<\mu}X_\alpha\in\calD$. 

If $X\in([\lambda]^{\LT\lambda})^\uniV$, then there is $\gamma<\lambda$ \st\ 
$X\subseteq\gamma$ by regularity of $\lambda$. Since 
${f^*}^{-1}\imageof\gamma=\setof{a\in(\Pkl{}{})^\uniV}{f^*(a)<\gamma}$ is 
disjoint with $\setof{a\in(\Pkl{}{})^\uniV}{f^*(a)\geq\gamma}$ and since the 
latter set is in $\calU$ by the definition \xitemof{laver-1-6-0}, it follows that
$X\not\in\calD$. 
\qedofClaim
}\fi
\begin{Claim}
  \Label{Cl-laver-3}
  $[id_\lambda]_\calD=\sup(j_\calD\imageof\lambda)$. 
\end{Claim}
\prfofClaim
Suppose $\gamma<\lambda$. Then $\setof{\alpha<\lambda}{c_\gamma(\alpha)<id_\lambda(\alpha)}
=\setof{\alpha<\lambda}{\gamma<\alpha}=\lambda\setminus(\gamma+1)\in\calD$ since 
$\calD$ is homogeneous by \Claimof{Cl-laver-2}. Thus
$j_\calD(\gamma)=[c_\gamma]_\calD<[id_\lambda]_\calD$.

Suppose now that $\mapping{g}{\lambda}{\On}$ is 
\st\ $[g]_\calD<[id_\lambda]_\calD$. This means that
\begin{xitemize}
\xitem[laver-1-8] 
  $\setof{\alpha<\lambda}{g(\alpha)<\underbrace{id_\lambda(\alpha)}_{\mbox{\small$=\alpha$}}}\in\calD$\\[\jot]
  $\Leftrightarrow$\ \
  $\calU\ni{f^*}^{-1}\imageof\setof{\alpha<\lambda}{g(\alpha)<\alpha}
  =\setof{a\in(\Pkl{}{})^\uniV}{g(f^*(a))<f^*(a)}$.
\end{xitemize}
By definition of $f^*$, there is $\gamma<\lambda$ \st\
$[gf^*]_\calU\leq[c_\gamma]_\calU$. That is,\\
$\setof{a\in(\Pkl{}{})^\uniV}{g(f^*(a))\leq\gamma}\in\calU$. 
This is 
equivalent to $\setof{\alpha<\lambda}{g(\alpha)\leq\gamma}\in\calD$ or
$[g]_\calD\leq j_\calD(\gamma)$. 
\qedofClaim
\begin{Claim}
  \Label{Cl-laver-4}\wassertof{1}\ \ $\setof{a\in(\Pkl{}{})^\uniV}{f^*(a)=\sup a\cap f^*(a)}\in\calU$.\smallskip
  
  \assert{2} $\setof{\alpha<\lambda}{\cf(\alpha)<\kappa}\in\calD$. 
\end{Claim}
\prfofClaim
\assertof{1}: In $\uniV$, let $\mapping{g}{\Pkl{}{}}{\lambda}$;
$a\mapsto \sup(a\cap f^*(a))$. Then we have $[g]_\calU\leq[f^*]_\calU$.

On the other hand, for each $\gamma<\lambda$, we have 
\begin{xitemize}
\xitem[]
  $\setof{a\in(\Pkl{}{})^\uniV}{\hspace{-2em}\underbrace{c_\gamma(a)\leq g(a)}_{
    \mbox{\small\mbox{}$\ \ \ \Leftrightarrow\ \ 
    \gamma\leq\sup(a\cap f^*(a))$}}\hspace{-2em}}$\\[3\jot]
  $\supseteq\underbrace{\setof{a\in(\Pkl{}{})^\uniV}{\gamma\leq f^*(a)}}_{\mbox{\small 
    \ \ \ $\in\,\calU$ (by the choice of $f^*$)}}
  \cap\underbrace{\setof{a\in(\Pkl{}{})^\uniV}{\gamma\in a}}_{\mbox{\small\ 
      \ \ $\in\,\calU$ (since $\calU$ is fine)}}\in\calU$.
\end{xitemize}
Thus $j_\calU(\gamma)\leq[g]_\calU$. By the choice of $f^*$ it follows that $[f^*]_\calU\leq[g]_\calU$.
\smallskip

\assertof{2}: by the definition \xitemof{laver-1-7} of $\calD$
\begin{xitemize}
\xitem[laver-1-9] 
  $\setof{\alpha<\lambda}{\cf(\alpha)<\kappa}\in\calD$\ \ $\Leftrightarrow$\ \ 
  $\setof{a\in(\Pkl{}{})^\uniV}{\cf(f^*(a))<\kappa}\in\calU$.
\end{xitemize}

By \assertof{1}, 
\begin{xitemize}
\xitem[] 
  $\setof{a\in(\Pkl{}{})^\uniV}{\cf(f^*(a)),\kappa}$\\ $\supseteq
  \setof{a\in(\Pkl{}{})^\uniV}{f^*(a)=\sup(a\cap f^*(a))}\in\calU$. \qedofClaim
\end{xitemize}

In $\uniV$, 
let
\begin{xitemize}
\xitem[laver-1-10] 
  $A_\alpha=\left\{\,
  \begin{array}{@{}ll}
    \mbox{a cofinal subset of }\alpha\mbox{ of order type }\cf(\alpha), &\mbox{if }\cf(\alpha)<\kappa;\\
    \emptyset, &\mbox{otherwise}
  \end{array}
  \right.$
\end{xitemize}
for $\alpha<\lambda$.  

Let $\vec{A}=\seqof{A_\alpha}{\alpha<\lambda}$. 

By \Claimof{Cl-laver-3} and \Claimof{Cl-laver-4},\,\assertof{2}, 
\begin{xitemize}
\xitem[laver-1-11] 
  $\fnsp{\lambda}{\uniV}/\calD\modelof{[\vec{A}]_\calD
  \mbox{ is a cofinal subset of }\sup(j_\calD\imageof\lambda)}$.
\end{xitemize}

By the $\mu$-cc of $\poP$, there is a strictly and continuously increasing sequence 
$\seqof{\eta_\xi}{\xi<\lambda}$ of ordinals $<\lambda$ in $\uniV$ \st, letting $I_\xi=[\eta_\xi,\eta_{\xi+1})$, 
\begin{xitemize}
\xitem[laver-1-12] 
  $\forces{\poP}{j_{\utilde{\calD}}(I_\xi)\cap[\vec{A}]_{\utilde{\calD}}\not\equiv\emptyset}$
\end{xitemize}
where $\utilde{\calD}$ is a $\poP$-name for $\calD$. Let
\begin{xitemize}
\xitem[laver-1-13] 
  $M_\alpha=\setof{\xi<\lambda}{I_\xi\cap A_\alpha\not=\emptyset}$
\end{xitemize}
for $\alpha<\lambda$. Note that $\seqof{M_\alpha}{\alpha<\lambda}\in\uniV$.

Since $\lambda\geq2^{\LT\kappa}$, 
\assertof{3} of the following \Claimof{Cl-laver-5} implies that $\lambda^{\LT\kappa}=\lambda$ 
as desired.

\begin{Claim}
  \Label{Cl-laver-5}
  \wassertof{1}\ \ $\cardof{M_\alpha}<\kappa$ for all $\alpha<\lambda$.\smallskip

  \assert{2} For any $\xi<\lambda$,
  $\setof{\alpha<\lambda}{\xi\in M_\alpha}\in\calD$.\smallskip

  \assert{3} For any $c\in(\Pkl{}{})^\uniV$,
  $\setof{\alpha<\lambda}{c\subseteq M_\alpha}\in\calD$. In particular, there is 
  some $\alpha<\lambda$ \st\ $c\subseteq M_\alpha$. 
\end{Claim}
\prfofClaim
\assertof{1}: Since $\cardof{A_\alpha}<\kappa$ (see \xitemof{laver-1-10}) and
$I_\xi$'s are pairwise disjoint, there can be only $<\kappa$ elements 
of $M_\alpha$.\smallskip

\assertof{2}: Suppose $\xi<\lambda$. 
By \xitemof{laver-1-12}, we have
\begin{xitemize}
\xitem[laver-1-14] 
  $j_\calD(I_\xi)\cap [\vec{A}]_\calD\not=\emptyset$.
\end{xitemize}

By \L o\'s's Theorem, it follows that
\begin{xitemize}
\xitem[laver-1-15] 
  $\setof{\alpha<\lambda}{
  \underbrace{I_\xi\cap A_\alpha\not=\emptyset}_{\mbox{\small\ \ $\Leftrightarrow\ \ \xi\in M_\alpha$}}}\in\calD$. 
\end{xitemize}

\assertof{3}: Suppose $c\in(\Pkl{}{})^{\uniV}$. For each $\xi\in c$, we have
$\setof{\alpha<\lambda}{\xi\in M_\alpha}\in\calD$ by \assertof{2}. By 
$\kappa$-completeness of $\calD$, it follows that
\begin{xitemize}
\xitem[] 
  $\setof{\alpha<\lambda}{c\subseteq M_\alpha}=\bigcap_{\,\xi\in c\,}\setof{\alpha<\lambda}{\xi\in M_\alpha}\in\calD$.
\qedofClaim
\end{xitemize}

\assertof{2} (of \Propof{T-laver-1-1}): Suppose that $\lambda\geq\kappa$ is a regular cardinal,
$\poP\in\calP$, $\genG$ is a $(\uniV,\poP)$-generic filter and $j$, 
$M\subseteq\uniV[\genG]$ are \st, $M$ is transitive in $V[\genG]$,
$\elembed{j}{\uniV}{M}$, $crit(j)=\kappa$, $j(\kappa)>\lambda$ and
$j\imageof\lambda\in M$. Let $\calU$ be defined by \xitemof{laver-1-6}. By 
\Claimof{Cl-laver-0}, $\calU$ is a $V$-normal ultrafilter.

In $\uniV$, let $\utilde{\calU}$ be a $\poP$-name of $\calU$ and let
\begin{xitemize}
\xitem[laver-1-15-0] 
  $\calF=\setof{U\in\psof{\Pkl{}{}}}{\forces{\poP}{\check{U}\varin\utilde{\calU}}}$.
\end{xitemize}
Then $\calF$ is a $\mu$-saturated normal filter over
$\Pkl{}{}$. 
\qedofProp\qedskip

\ifextended{\tt
\Thmabove,\,\assertof{2} follows from the next Lemma. Note that the normality of a filter 
over $\Pkl{}{}$ we use here is as defined in \cite{higher-inf:} p.301 
and $\uniV$-normality a natural modification of this definition. 

\begin{LemmaA}
  \Label{LA-laver-0} Suppose $\kappa<\lambda$ and $\mu$ are uncountable regular 
  cardinals and $\poP$ a $\mu$-cc \po\ \st\ $\poP$ preserves the cardinal $\kappa$. If
  \begin{xitemize}
  \xitemA[laverA-0]
    $\forces{\poP}{\mbox{ there is a }\uniV\mbox{-normal ultrafilter }
    \calU\mbox{ over }(\Pkl{}{})^\uniV}$,
  \end{xitemize}
  then there is a $\LT\mu$-saturated normal filter $\calF$ on $\Pkl{}{}$. 
\end{LemmaA}
\prf Let $\utcalU$ be a $\poP$ name \st\
\begin{xitemize}
\xitemA[laverA-1] 
  $\forces{\poP}{\utcalU\mbox{ is a }\uniV\mbox{-normal ultrafilter over }(\Pkl{}{})^\uniV}$.
\end{xitemize}

Let 
\begin{xitemize}
\xitemA[laverA-2] 
  $\calF=\setof{U\in\psof{\Pkl{}{}}}{\forces{\poP}{\check{U}\varin\utcalU}}$.
\end{xitemize}

The following claim shows that $\calF$ above is as desired. 
\begin{ClaimA}
  \Label{ClA-laver-0}
  \assert{0} $\calF$ is a filter over $\Pkl{}{}$.\smallskip

  \assert{1} $\calF$ is $\LT\kappa$-complete.\smallskip

  \assert{2} $\calF$ is fine.\smallskip

  \assert{3} $\calF$ is normal.\smallskip

  \assert{4} $\calF$ is $\LT\mu$-saturated.
\end{ClaimA}
\prfofClaimA
\assertof{0}: By \xitemAof{laverA-1} and \xitemAof{laverA-2}.\smallskip

\assertof{1}: Suppose that $U_\alpha\in\calF$, $\alpha<\delta$ for some
$\delta<\kappa$. Then
$\utS=\setof{\pairof{\check{U}_\alpha,\bbone_\poP}}{\alpha<\delta}$ is 
a $\poP$-name and 
\begin{xitemize}
\xitemA[laverA-3] 
  $\forces{\poP}{\utS\subseteq\utcalU,\ \cardof{\utS}<\kappa,\ \utS\varin\uniV}$. 
\end{xitemize}
Thus $\forces{\poP}{\bigcap\utS\varin\utcalU}$ by \xitemAof{laverA-1}.
Since
$\forces{\poP}{\bigcap\utS\equiv
  \CHECK{\poP}{\bigcap_{\alpha<\delta}U_\alpha}}$, we have $\bigcap_{\alpha<\delta}U_\alpha\in\calF$.
\smallskip

\assertof{2}: By \xitemAof{laverA-1} and \xitemAof{laverA-2}. \assertof{3}: 
Similarly to \assertof{1}.\smallskip

For the proof of \assertof{4}, we need the following:
\begin{ClaimA}
  \label{ClA-laver-1} For $X\subseteq\Pkl{}{}$, $X$ is $\calF$-stationary if and 
  only if $\condp\forces{\poP}{\check{X}\varin\utcalU}$ for some $\condp\in\poP$. 
\end{ClaimA}
\prfofClaimA
Suppose that $\condp\notforces{\poP}{\check{X}\varin\utcalU}$ for all
$\condp\in\poP$. Then we have $\forces{\poP}{\check{X}\notvarin\utcalU}$. 
Since $\utcalU$ is a $\poP$-name of an ultrafilter, if follows that
$\forces{\poP}{(\Pkl{}{}^\uniV\setminus\check{X})\varin\utcalU}$.
Since
$\forces{\poP}{(\Pkl{}{}^\uniV\setminus\check{X})\equiv\CHECK{\poP}{\Pkl{}{}\setminus X}}$,
it follows that $(\Pkl{}{}\setminus X)\in\calF$. Thus $X$ is 
not $\calF$-stationary.

Suppose now that $\condp\forces{\poP}{\check{X}\varin\utcalU}$. Then for any
$U\in\calF$, we have $\condp\forces{\poP}{\check{U}\cap\check{X}\varin\utcalU}$. 
Thus $\condp\forces{\poP}{\check{U}\cap\check{X}\not\equiv\emptyset}$ and hence 
$U\cap X\not=\emptyset$. This shows that $X$ is $\calF$-stationary.
\qedofClaimA
\qedskip

\noindent
{\bf Continuation of the proof of \ClaimAof{ClA-laver-0}:}\\
\assertof{4}: Suppose that $\seqof{X_\alpha}{\alpha<\mu}$ is a sequence 
of $\calF$-stationary sets. By \ClaimAof{ClA-laver-1}, there are
$\condp_\alpha\in\poP$, $\alpha<\mu$ \st\
$\condp_\alpha\forces{\poP}{\check{X}_\alpha\varin\utcalU}$ for $\alpha<\mu$.

Since $\poP$ is $\mu$-cc, there are $\alpha<\alpha'<\mu$ \st\ $\condp_\alpha$ and 
$\condp_{\alpha'}$ are compatible in $\poP$. Let $\condr\in\poP$ be \st\
$\condr\leq_\poP\condp_\alpha$, $\condp_{\alpha'}$. Then we have
$\condr\forces{\poP}{\check{X}_\alpha\cap\check{X}_{\alpha'}\varin\utcalU}$. 
Again by \ClaimAof{ClA-laver-1}, it follows that $X_\alpha\cap X_{\alpha'}$ is
$\calF$-stationary. 
\qedofClaimAof{ClA-laver-0}
\\
\qedofLemmaA
\qedskip
}\fi 

In case of $\mu\leq2^{\LT\kappa}$, 
\Propof{T-laver-1-1},\,\assertof{1} can be also proved as follows:
By (the 
proof of) Corollary 3.5 and Theorem 5.1 in \cite{MaSaUs}, we obtain the theorem asserting:
\begin{Thm}{\rm (Matsubara, Sakai and Usuba \cite{MaSaUs})}
  For a regular uncountable $\kappa$, if there are $\lambda$-saturated fine ideal over $\Pkl{}{}$ for 
  all $\lambda\geq\kappa$, then \SCH\ holds above $2^{\LT\kappa}$. \qed
\end{Thm}
\Propof{T-laver-1-1},\,\assertof{1} follows immediately from this theorem and 
\Propof{T-laver-1-1},\,\assertof{2}. 

In \cite{I}, 
it is proved that, for a cardinal $\kappa$, if $\kappa^+$ is 
generically supercompact by $\LT\kappa$-closed \pos\ then
$\SDLS_+(\calL^{\aleph_0}_{stat},\LT\kappa^+)$ holds (this follows from \Thmof{T-grp-0}, 
\Thmof{T-grp-a} and \Lemmaof{L-sdls-a} in \cite{I}). 

In the next section, we show that, for a regular $\kappa$, the 
assertion ``$(\ast)^{\intnl+}_{\LT\kappa,\lambda}$ holds for all 
regular $\lambda\geq\kappa$'' can be characterized in terms of SDLS for $\calL^{\aleph_0}_{stat}$ in internal 
interpretation. 
This property 
``$(\ast)^{\intnl+}_{\LT\kappa,\lambda}$ for all regular $\lambda\geq\kappa$'' holds  
under a wider class of generic supercompactness (c.f. 
\Thmof{T-laver-2},\,\assertof{1} and \assertof{2}).

\begin{Thm}
  \Label{P-lt-conti-3} Suppose that $\kappa$ is a generically supercompact 
  cardinal by
  proper \pos. Then $(\ast)^{\intnl+}_{\LT\kappa,\lambda}$ holds for all 
  regular $\lambda\geq\kappa$. 
\end{Thm}

In the proof of \Thmof{P-lt-conti-3}, we use the case $\kappa=\aleph_1$ of the following well-known fact:
\begin{Lemma}
  \Label{P-lt-conti-4}
  Suppose that $M$ is an inner model in $\uniV$. For ordinals $\kappa$, $\lambda$ 
  with $\kappa\leq\lambda$ and $\uniV\modelof{\kappa\xmbox{ is a regular uncountable cardinal}}$, 
  and for $A\in M$ with $M\modelof{A\subseteq\Pkl{}{}}$, if\/
  $\uniV\modelof{A\xmbox{ is stationary in }\Pkl{}{}}$, then
  $M\modelof{A\xmbox{ is stationary in }\Pkl{}{}}$. \ifextended\else\qed\fi
\end{Lemma}
\ifextended{\tt
  \prf
  Suppose that $\uniV\modelof{A\mbox{ is stationary in }\Pkl{}{}}$.
  Let $\calC\in M$ be \st\
  $M\modelof{\calC\mbox{ is a club in }\Pkl{}{}}$. By Kueker's 
  Theorem (see Exercise 38.10 in \cite{millennium-book:}) there is $F\in M$ \st\
  \begin{xitemize}
  \xitemA[P-lt-conti-A-0] 
    $M\modelof{
    \begin{array}[t]{@{}l}
      \mapping{F}{[\lambda]^{<\aleph_0}}{\lambda}\mbox{ and}\\
      \begin{array}{@{}l@{}l}
        \calC_F=\setof{a\in\Pkl{}{}}{&a\cap\kappa\in\kappa\mbox{ and }\\
          &a\mbox{ is closed under }F}\subseteq\calC}.
      \end{array}
    \end{array}
    $
  \end{xitemize}
  In $\uniV$, there is $a\in A\cap{\calC_F}^\uniV$. $a\in A\subseteq M$
  and hence $M\modelof{a\in A\cap\calC_F}$. \qedofLemma\qedskip
}\fi

\noindent
\prfof{\Thmof{P-lt-conti-3}} For a regular $\lambda\geq\kappa$, let
$\lambda^*=\cardof{\calH(\lambda)}$ and $\lambda^{**}=(2^{\lambda^*})^+$. 

Suppose that\quad
$\gmA=\pairof{\calH(\lambda),\in\ctenten}$\quad is a countable expansion of the structure
$\pairof{\calH(\lambda),\in}$, $\seqof{S_a}{a\in\calH(\lambda)}$ a sequence of 
stationary subsets of $[\calH(\lambda)]^{\aleph_0}$ and 
$\calC\subseteq[\calH(\lambda)]^{\LT\kappa}$ a club.

Let $\mapping{\iota}{\calH(\lambda)}{\lambda^*}$ be a bijection and let
$\gmA^*=\pairof{\lambda^*,E^*\ctenten}$ be a copy of $\gmA$ translated 
by $\iota$. 

Let 
$\seqof{S^*_\alpha}{\alpha\in\lambda^*}$ be \st\ 
\begin{xitemize}
\xitem[lt-conti-2-3] 
  $S^*_\alpha=\iota\imageof S_{\iota^{-1}(\alpha)}$ for $\alpha\in\lambda^*$ and 
\end{xitemize}
\begin{xitemize}
\xitem[lt-conti-2-4] 
  $\calC^*=\setof{\iota\imageof X}{X\in\calC}$. 
\end{xitemize}
Thus, we have 
\begin{xitemize}
\xitem[lt-conti-2-5] 
  $EXT_{E^*}(S^*_\alpha)$ is a stationary subset of $[\lambda^*]^{\aleph_0}$ for 
  all $\alpha\in\lambda^*$ and 
\end{xitemize}
\begin{xitemize}
\xitem[lt-conti-2-6] 
  $\calC^*$ is a club in $[\lambda^*]^{\LT\kappa}$
\end{xitemize}
where, for $S\subseteq\lambda^*$,
$EXT_{E^*}(S)=\setof{\setof{\beta\in\lambda^*}{\beta\mathrel{E^*}\alpha}}{\mbox{ for }\alpha\in S}$
is the set of all extents of elements of $S$ \wrt\ the element relation $E^*$.

It is enough to show that there is an $X\in\calC^*$ \st\ 
\begin{xitemize}
\item[\wassertof{2$'$}] $\gmA^*\restr X\prec\gmA^*$, and
\item[\wassertof{3$''$}] $EXT_{E^*}(S^*_\alpha\cap X)$ is stationary in 
$[X]^{\aleph_0}$ for all $\alpha\in X$.  
\end{xitemize}

Let $\poP$ be a proper \po\ and $\genG$ a $(\uniV,\poP)$-generic 
filter \st\ there are transitive $M\subseteq\uniV[\genG]$ and
$\elembed{j}{\uniV}{M}$ \st\ 
\begin{xitemize}
\xitem[lt-conti-3] 
  $\kappa=crit(j)$,
\xitem[lt-conti-4] 
  $j(\kappa)>\lambda^{**}$ and
\xitem[lt-conti-5] 
  $j\imageof\lambda^{**}\in M$. 
\end{xitemize}
Let $X=j\imageof\lambda^*$. $X\in M$ by \xitemof{lt-conti-5} and \Lemmaof{L-lt-conti-0},\,\assertof{1}. 
Thus, we also have $j\restr\lambda^*\in M$ and
\begin{xitemize}
\xitem[costat-2-3] 
  $M\models\cardof{X}=\cardof{\lambda^*}< j(\kappa)$. 
\end{xitemize}

Since $\cardof{\calC^*}<\lambda^{**}$, we have $j\imageof\calC^*\in M$ by 
\Lemmaof{L-lt-conti-0},\,\assertof{1}. 
$M\modelof{\cardof{j\imageof\calC^*}\leq\cardof{j(\calC^*)}<j(\kappa)}$.  
Also we have
$M\modelof{j\imageof\calC^*\mbox{ is directed and }\bigcup(j\imageof\calC^*)=X}$. 
Since $M\modelof{\calC^*\xmbox{ is a club in }[j(\lambda^*)]^{\LT j(\kappa)}}$ by 
elementarity and \xitemof{lt-conti-2-6}, 
it follows that 
\begin{xitemize}
\xitem[costat-2-4] 
  $X\in j(\calC^*)$.
\end{xitemize}

Let
$\seqof{\tilde{S}^*_\alpha}{\alpha<j(\lambda^*)}=j(\seqof{S^*_\alpha}{\alpha<\lambda^*})$. 

By elementarity and applying Vaught's test, we can show
\begin{xitemize}
\xitem[] 
  $M\models j(\gmA^*)\restr X\prec j(\gmA^*)$.
\end{xitemize}
For $\alpha\in\lambda^*$, since
$\tilde{S}^*_{j(\alpha)}=j(S^*_\alpha)\supseteq j\imageof S^*_\alpha$, we have 
\begin{xitemize}
\xitem[costat-2-5] 
  $M\modelof{EXT_{j(E^*)}(\tilde{S}^*_{j(\alpha)}\cap X)
  \supseteq EXT_{j(E^*)}(j\imageof S^*_\alpha)}$.
\end{xitemize}

By \xitemof{lt-conti-2-5} and since $\poP$ is proper, we have
$\uniV[\genG]\modelof{EXT_{E^*}(S^*_\alpha)\mbox{ is stationary in }[\lambda^*]^{\aleph_0}}$. 
Since $EXT_{j(E^*)}(j\imageof S^*_\alpha)$ is a translation of
$EXT_{E^*}(S^*_\alpha)$ induced by $j\restr\lambda^*$, it follows that 
\begin{xitemize}
\xitem[costat-2-6] 
  $\uniV[\genG]\modelof{EXT_{j(E^*)}(j\imageof S^*_\alpha)
  \mbox{ is stationary in }[X]^{\aleph_0}}$.
\end{xitemize}
By \Lemmaof{P-lt-conti-4}, it follows that 
\begin{xitemize}
\xitem[costat-2-7] 
  $M\modelof{EXT_{j(E^*)}(j\imageof S^*_\alpha)
  \mbox{ is stationary in }[X]^{\aleph_0}}$.
\end{xitemize}
Thus, we have
\begin{xitemize}
\xitem[] 
  $M\modelof{\exists X\in j(\calC^*)\,\big(\,
  \begin{array}[t]{@{}l}
    j(\gmA^*)\restr X\prec j(\gmA^*)  \land EXT_{j(E^*)}(\tilde{S}^*_\xi\cap 
    X)\mbox{ is}
    \\\mbox{a stationary subset of }[X]^{\aleph_0}\mbox{ for all }\xi\in X\big)}.
  \end{array}
  $
\end{xitemize}
By elementarity, it follows that 
\begin{xitemize}
\xitem[] 
  $\uniV\modelof{\exists X\in\calC^*\,\big(\,
  \begin{array}[t]{@{}l}
    \gmA^*\restr X\prec \gmA^*\  \land\ EXT_{E^*}(S^*_\xi\cap X)
    \mbox{ is}\\\mbox{a stationary subset of }[X]^{\aleph_0}\mbox{ for all }\xi\in X\big)}.
  \end{array}
  $
\end{xitemize}
\mbox{}\vspace{-3ex}

\mbox{}\qedof{\Thmof{P-lt-conti-3}}

\section{Internal interpretation of stationary logic}
\Label{internal}
For a structure $\gmA=\pairof{A\ctenten}$ of a countable signature, an
$\calL^{\aleph_0}_{stat}$-formula
$\varphi=\varphi(x_0\ctentenc X_0\ctenten)$\footnote{As before, when we write
  $\varphi=\varphi(x_0\ctentenc X_0\ctenten)$, we always assume that the list $x_0\ctenten$ 
  contains all the free first order variables of $\varphi$ and $X_0\ctenten$ 
  all the free weak second order variables of $\varphi$.} and
$a_0\ctenten\in A$, $U_0\ctenten\in[A]^{\aleph_0}\cap A$, we define the internal 
interpretation of $\varphi(a_0\ctentenc U_0\ctenten)$ in $\gmA$ (notation:
$\gmA\models^{\intnl}\varphi(a_0\ctentenc U_0\ctenten)$ for
``$\varphi(a_0\ctentenc U_0\ctenten)$ holds internally in $\gmA$'') by induction 
on the construction of $\varphi$ as follows:

If $\varphi$ is ``$x_i\varin X_j$'' then
\begin{xitemize}
\xitem[internal-0] 
  $\gmA\models^{\intnl}\varphi(a_0\ctentenc U_0\ctenten)$\ \ $\Leftrightarrow$\ \
  $a_i\in U_j$
\end{xitemize}
for a structure $\gmA=\pairof{A\ctenten}$, $a\in A$ and
$U\in [A]^{\aleph_0}\cap A$. 

For first-order connectives and quantifiers in
$\calL^{\aleph_0}_{stat}$, the semantics ``$\models^{\intnl}$'' is defined exactly as for 
the first order ``$\models$''.

For 
an $\calL^{\aleph_0}_{stat}$ formula $\varphi$ with
$\varphi=\varphi(x_0\ctentenc X_0\ctentenc X)$, assuming that the notion of 
$\gmA\models^{\intnl}\varphi(a_0\ctentenc U_0\ctentenc U)$ has 
been defined for all $a_0\ctenten\in A$, $U_0\ctentenc U\in[A]^{\aleph_0}\cap A$, 
we stipulate 
\begin{xitemize}
\xitem[internal-1] 
  $\gmA\models^{\intnl}stat\,X\,\varphi(a_0\ctentenc U_0\ctentenc X)$\ \
  $\Leftrightarrow$\\
  \qquad$\setof{U\in[A]^{\aleph_0}\cap A}{\gmA\models^{\intnl}\varphi(a_0\ctentenc U_0\ctentenc U)}$
  is stationary in $[A]^{\aleph_0}$
\end{xitemize}
for a structure $\gmA=\pairof{A\ctenten}$ of a relevant signature,  
$a_0\ctenten\in A$ and $U_0\ctenten\in[A]^{\aleph_0}\cap A$.

For structures $\gmA$, $\gmB$ of the same signature 
with $\gmB=\pairof{B\ctenten}$ and $\gmB\subseteq\gmA$, 
we define
\begin{xitemize}
\xitem[internal-2] 
  {\it 
  $\gmB\prec^{\intnl}_{\calL^{\aleph_0}_{stat}}\gmA$\ \ $\Leftrightarrow$\\[2\jot]
  \qquad $\gmB\models^{\intnl}\varphi(b_0\ctentenc U_0\ctenten)$ if and only if 
  $\gmA\models^{\intnl}\varphi(b_0\ctentenc U_0\ctenten)$\\
  \qquad for all $\calL^{\aleph_0}_{stat}$-formulas $\varphi$ in the signature 
  of the structures with\\
  \qquad $\varphi=\varphi(x_0\ctentenc X_0\ctenten)$, $b_0\ctenten\in B$ and $U_0\ctenten\in[B]^{\aleph_0}\cap B$.}
\end{xitemize}

Finally, for a regular $\kappa>\aleph_1$, the internal strong downward 
L\"owenheim-Skolem Theorem 
$\SDLS^{\intnl}_+(\calL^{\aleph_0}_{stat},\LT\kappa)$ 
is defined by 
\begin{xitemize}
\item[$\SDLS^{\intnl}_+(\calL^{\aleph_0}_{stat},\LT\kappa)$:]  {\it For any structure $\gmA=\pairof{A\ctenten}$ of 
  countable signature with $\cardof{A}\geq\kappa$, there are {\it stationarily many} 
  $M\in[A]^{\LT\kappa}$ \st\\ $\gmA\restr M\prec^{\intnl}_{\calL^{\aleph_0}_{stat}}\gmA$. }
\end{xitemize}

`$+$' in ``$\SDLS^{\intnl}_+(\calL^{\aleph_0}_{stat},\LT\kappa)$'' refers to 
the ``stationarily many'' existence of the reflection points $M$. Similarly to 
\Lemmaof{P-sdls-0} in \cite{I}, this additional condition can be drooped if
$\kappa=\aleph_2$. This is because the quantifier $Qx\,\varphi$ defined by
$stat\,X\exists x\,(x\notvarin X\land\varphi$, 
$\gmA\models^{\intnl}Qx\,\varphi(x\ctenten))$ still implies that ``there are 
uncountably many $a\in A$ with $\varphi(a\ctenten)$''. 
Note that, if $\gmA\models^{\intnl}\neg stat\,X\,(x\equiv x)$, for a structure
$\gmA=\pairof{A\ctenten}$,  we can easily find even club many  
$X\in[A]^{<\kappa}$ for any regular $\aleph_1\leq\kappa\leq\cardof{A}$ \st\ 
$\gmA\restr X\prec^{\intnl}_{\calL^{\aleph_0}_{stat}}\gmA$. 

\begin{Prop}
  \Label{P-internal-0}
  For a regular cardinal $\kappa>\aleph_1$, \tfae:\smallskip

  \assert{a} $(\ast)^{\intnl+}_{\LT\kappa,\lambda}$ holds for all 
  regular $\lambda\geq\kappa$.\smallskip

  \assert{b} $\SDLS^{\intnl}_+(\calL^{\aleph_0}_{stat},\LT\kappa)$ holds.
\end{Prop}
\prf This proof is a straightforward modification of the proof of \Lemmaof{P-DRP-3},\assertof{1} in 
\cite{I}. Nevertheless, we will give the complete proof since we have to modify 
it further to prove \Propof{P-PKL-0}.

Suppose first that $\SDLS^{\intnl}_+(\calL^{\aleph_0}_{stat}, \LT\kappa)$ 
holds. We show that 
$(\ast)^{\intnl+}_{\LT\kappa,\lambda}$ holds for all $\lambda\geq\kappa$. Let 
$\lambda\geq\kappa$. Let $\tilde{\gmA}$ be a countable expansion of
$\pairof{\calH(\lambda),\in}$, $\seqof{S_a}{a\in\calH(\lambda)}$ a 
sequence of stationary subsets of $[\calH(\lambda)]^{\aleph_0}$ and
$\calD\subseteq[\calH(\lambda)]^{\LT\kappa}$ a club.

Let
\begin{xitemize}
\xitem[internal-3] 
  $\tilde{\gmA}^*=\pairof{\underbrace{\calH(\lambda)\ctentenc\in}_{\tilde{\gmA}},\symb{\vec{S}}^{\tilde{\gmA}^*}}$
\end{xitemize}
where $\symb{\vec{S}}$ is a binary relation symbol and 
\begin{xitemize}
\xitem[internal-4] 
  $\symb{\vec{S}}^{\tilde{\gmA}^*}=\setof{\pairof{a,s}\in(\calH(\lambda))^2}{s\in S_a}$.
\end{xitemize}
Let $M\in[\calH(\lambda)]^{\LT\kappa}$ be \st\
\begin{xitemize}
\xitem[internal-5] 
  $M\in\calD$ and 
\xitem[internal-6] 
  $\tilde{\gmA}^*\restr M\prec^{\intnl}_{\calL^{\aleph_0}_{stat}}\tilde{\gmA}^*$.
\end{xitemize}

By the choice of $\tilde{\gmA}^*$ and \xitemof{internal-6},
$\tilde{\gmA}^*\restr M\models^{\intnl}\forall x\,stat\,X\exists y\,(\symb{\vec{S}}(x,y)\land \forall z\,
(z\varin X\leftrightarrow z\in y))$ holds and hence, 
for all $a\in M$,  $S_a\cap M$ is stationary in $[M]^{\aleph_0}$.
\smallskip

Suppose now that $(\ast)^{\intnl+}_{\LT\kappa,\lambda}$ holds for all
$\lambda\geq\kappa$. Let
$\gmA=\pairof{A\ctenten}$ be a structure in 
countable signature and of cardinality $\geq\kappa$, and 
$\calD\subseteq[A]^{\LT\kappa}$ a club. \Wolog, we may assume that $\gmA$ is a relational structure.
Let $\lambda$ be a regular cardinal \st\ $\gmA\in\calH(\lambda)$. In particular, 
we have $A\subseteq\calH(\lambda)$. 

Let
$\tilde{\gmA}=\pairof{\calH(\lambda),\underbrace{\symb{A}^{\tilde{\gmA}}\ctenten}_{=\gmA},\in}$ 
where $\symb{A}$ is a unary relation symbol and $\symb{A}^{\tilde{\gmA}}=A$. 

For each $a\in\calH(\lambda)$, let 
\begin{xitemize}
\xitem[internal-7] 
  $S_a=\left\{\,
  \begin{array}[c]{@{}l}
    \setof{U\in[\calH(\lambda)]^{\aleph_0}}{\cardof{U\cap A}=\aleph_0,\,U\cap A\in A,\\[\jot]
    \phantom{{}U\in[\calH(\lambda)]^{\aleph_0}:\ \ \ \ }
            \gmA\models^{\intnl}\psi(\variables{a}{m-1},\variables{U}{n-1},U\cap A)},\\[\jot]
    \qquad\quad
    \mbox{if }\psi=\psi(\variables{x}{m-1},\variables{Y}{n-1},X)\mbox{ is an }
    \calL^{\aleph_0}_{stat}\mbox{-formula}\\
    \qquad\quad\mbox{in the signature of }\gmA,\,a_0\ctenten\in 
    A,\,U_0\ctenten\in[A]^{\aleph_0}\cap A,\\
    \qquad\quad\gmA\models^{\intnl}stat\,X\,\psi(\variables{a}{m-1},\variables{U}{n-1},X)\mbox{ and}\\
    \qquad\quad a=\pairof{\psi,\variables{a}{m-1},\variables{U}{n-1}};\\[2\jot]

    [\calH(\lambda)]^{\aleph_0},\\\qquad\quad\mbox{otherwise.}
  \end{array}
  \right.$
\end{xitemize}

Let
\begin{xitemize}
\xitem[internal-8] 
  $\tilde{\calD}=\setof{U\in[\calH(\lambda)]^{\LT\kappa}}{U\cap A\in\calD}$. 
\end{xitemize}
$\tilde{\calD}$ contains a club in $[\calH(\lambda)]^{\LT\kappa}$.

By $(\ast)^{\intnl+}_{\LT\kappa,\lambda}$, there is an 
$M\in[\calH(\lambda)]^{\LT\kappa}$ \st\
\begin{xitemize}
\xitem[internal-9] 
  $M\in\tilde{\calD}$,
\xitem[internal-10] 
  $\tilde{\gmA}\restr M\prec\tilde{\gmA}$\ \ and 
\xitem[internal-11] 
  $S_a\cap\Pkl{\kappa\cap M}{M}\cap M$ is stationary in $\Pkl{\kappa\cap M}{M}$ 
  for all $a\in M$. 
\end{xitemize}

We expand $\tilde{\gmA}$ further by adding the relation
$\setof{\pairof{u,a}\in(\calH(\lambda))^2}{u\in S_a}$ to it as the interpretation of 
the binary relation symbol $\symb{\vec{S}}$. For simplicity, we shall also call 
this expanded structure $\tilde{\gmA}$\footnote{This expansion becomes necessary 
  below when we would like to have\\
  $\setof{U\cap B}{U\in S_a\cap[\tilde{B}]^{\aleph_0}\cap\tilde{B}}
   =\setof{V\in B}{V=U\cap B\mbox{ for some }U\in S_a\cap[\tilde{B}]^{\aleph_0}\cap\tilde{B}}$.  }. 

Let $\tilde{\gmB}=\tilde{\gmA}\restr M$ and 
let $B=\symb{A}^{\tilde\gmB}=A \cap M$ and $\gmB=\gmA\restr B$.  
Denoting the underlying set of $\tilde{\gmB}$ by $\tilde{B}$, we have 
$\tilde{B}=M$. 
$B\in\calD$ by \xitemof{internal-9} and the definition \xitemof{internal-8} of
$\tilde{\calD}$. 


By the elementarity $\tilde{\gmB}\prec\tilde{\gmA}$ \xitemof{internal-10}, the following Claim implies
$\gmB\prec^{\intnl}_{\calL^{\aleph_0}_{stat}}\gmA$.

\begin{Claim} For any $\calL^{\aleph_0}_{stat}$-formula
  $\varphi(\variables{x}{m-1},\variables{Y}{n-1})$ in the signature of the 
  structures $\gmA$, $\variables{a}{m-1}\in B$ and
  $\variables{U}{n-1}\in[B]^{\aleph_0}\cap B$, we have
  \begin{xitemize}
  \xitem[internal-12] 
    $\tilde{\gmB}\modelof{\gmA\models^{\intnl}\varphi(a_0\ctentenc U_0\ctenten)}$\ \ 
    $\Leftrightarrow$\ \
    $\gmB\models^{\intnl}\varphi(a_0\ctentenc U_0\ctenten)$.
  \end{xitemize}
\end{Claim}
\prfofClaim
By induction on $\varphi$. The crucial step in the induction is 
when $\varphi$ is of the form $stat\,X\psi$ and \xitemof{internal-12} holds for $\psi$:

Suppose first that
$\tilde{\gmB}\modelof{\gmA\models^{\intnl}\varphi(a_0\ctentenc U_0\ctenten)}$ holds.
Then, by elementarity and by the definition of $\tilde{\gmA}$, we have
$\gmA\models^{\intnl} stat\,X\psi(a_0\ctentenc\variables{U}{n-1},X)$. Thus, letting
$a=\pairof{\varphi,\variables{a}{m-1},\variables{U}{n-1}}$, we have $a\in\tilde{B}$ and 
\begin{xitemize}
\xitem[internal-14] 
  $S_a=\setof{U\in[\calH(\lambda)]^{\aleph_0}}{
  \begin{array}[t]{@{}l}
    \cardof{U\cap A}=\aleph_0,\,
    U\cap A\in A,\\\gmA\models^{\intnl}\varphi(\variables{a}{m-1},\variables{U}{n-1}, U\cap A)}
  \end{array}$ 
\end{xitemize}
by the definition \xitemof{internal-7} of $S_a$.  

By \xitemof{internal-11}, $S_a\cap[\tilde{B}]^{\aleph_0}\cap\tilde{B}$ is stationary in $[\tilde{B}]^{\aleph_0}$. 
It follows that
\begin{xitemize}
\xitem[internal-15] 
  $\setof{U\cap B}{\cardof{U\cap B}=\aleph_0, U\in S_a\cap[\tilde{B}]^{\aleph_0}\cap \tilde{B}}$\\[\jot]
  $=\setof{U\cap B}{{}
  \begin{array}[t]{@{}l}
    \cardof{U\cap B}=\aleph_0,\,U\cap B\in B,\,\\
    \tilde{\gmB}\modelof{\gmA\models^{\intnl}\psi(\variables{a}{m-1},\variables{U}{n-1},U\cap B)}}
  \end{array}
  $\\
  \mbox{}\hfill(by elementarity and \xitemof{internal-14})\\
  $\subseteq\setof{V\in[B]^{\aleph_0}\cap  B}{\gmB\models^{\intnl}
    \psi(\variables{a}{m-1},\variables{U}{n-1},V)}$\\
  \mbox{}\hfill(by induction hypothesis)
\end{xitemize}
is stationary.
Thus 
$\gmB\models^{\intnl} stat\,X\,\psi(\variables{a}{m-1},\variables{U}{n-1},X)$ holds, that is,\\
$\gmB\models^{\intnl} \varphi(\variables{a}{m-1},\variables{U}{n-1})$.

Suppose now that
$\tilde{\gmB}\not\modelof{\gmA\models^{\intnl}\varphi(\variables{a}{m-1},\variables{U}{n-1})}$. 
Then we have
\begin{xitemize}
\xitem[internal-16] 
  $\tilde{\gmB}\modelof{
  \begin{array}[t]{@{}l}
    \mbox{there is a club }\calC\subseteq[\symb{A}]^{\aleph_0}\mbox{ \st\ }
    \gmA\models^{\intnl}\neg\psi(a_0\ctentenc\variables{U}{n-1},x)\\
    \mbox{for all }x\in\calC}.
  \end{array}$
\end{xitemize}

By elementarity, there is a $\calC_0\in\tilde{B}$ \st\ $\calC_0$ is a club $\subseteq$ 
$[A]^{\aleph_0}$ and 
\begin{xitemize}
\xitem[internal-17] 
  $\tilde{\gmB}\modelof{\gmA\models^{\intnl}\neg\psi(\variables{a}{m-1},\variables{U}{n-1},V)}$
  for all $V\in\calC_0\cap B$.
\end{xitemize}

By induction hypothesis, it follows that
\begin{xitemize}
\xitem[internal-17-0] 
  $\setof{U\in[B]^{\aleph_0}\cap B}{
  \calB\models^{\intnl}\psi(a_0\ctentenc U_0\ctentenc U)}\cap\calC_0=\emptyset$. 
\end{xitemize}
Thus 
$\gmB\not\models^{\intnl}stat\,X\,\psi(a_0\ctentenc\variables{U}{n-1},X)$, i.e.\ 
$\gmB\not\models^{\intnl}\varphi(a_0\ctentenc U_0\ctenten)$. \\\qedofClaim\\
\qedofProp

\section{Stationarity quantifier in PKL logic}
\Label{PKL}
In this section, we consider a $\Pkl{}{}$ version of weak second-order 
logic with stationarity quantifier $\calL^{PKL}_{stat}$ and SDLS for this 
logic in internal interpretation.

One of the significant property of this SDLS is that 
it implies that 
the reflection cardinal is very large (see \Corof{Cl-a-stat-4}) and 
it is consistent (modulo supercompact cardinal) with the reflection cardinal
being ``$\LT2^{\aleph_0}$''.

For sets $s$ and $t$, we denote with $\Pkl{s}{t}$ the set
$[t]^{\LT\cardof{s}}=\setof{a\in\psof{t}}{\cardof{a}<\cardof{s}}$. We say
$S\subseteq\Pkl{s}{t}$ is stationary if it is stationary in the sense of Jech 
\cite{millennium-book:}. 

The logic $\calL^{PKL}_{stat}$ has a built-in unary predicate symbol $\symb{K}(\cdot)$. For 
a structure $\gmA=\pairof{A,\symb{K}^\gmA\ctenten}$, the weak second-order 
variables $X$, $Y$\ctenten\ run over elements of
$\Pkl{\symb{\scriptstyle K}^{\gmA}}{A}$.

We shall call a structure $\gmA$ with $\symb{K}$ in its signature 
as a unary predicate symbol \st\ $\cardof{\symb{K}^\gmA}$ is a regular 
uncountable cardinal, a {\it PKL-structure}.

$\calL^{PKL}_{stat}$ has 
the unique second-order quantifier ``$stat$'' and the internal interpretation $\models^{int}$ of 
formulas in this logic is defined similarly to $\calL^{\aleph_0}_{stat}$ 
with the crucial step in the inductive definition of $\models^{\intnl}$ being
\begin{xitemize}
\xitem[PKL-0]   $\gmA\models^{\intnl}stat\,X\,\varphi(a_0\ctentenc U_0\ctentenc X)$\ \
  $\Leftrightarrow$\\
  \qquad$\setof{U\in\Pkl{\symb{\scriptstyle K}^\gmA}{A}\cap A}{\gmA\models^{\intnl}\varphi(a_0\ctentenc U_0\ctentenc U)}$
  is stationary in\\
  \qquad$\Pkl{\symb{\scriptstyle K}^\gmA}{A}$
\end{xitemize}
for an $\calL^{PKL}_{stat}$-formula
$\varphi=\varphi(x_0\ctentenc X_0\ctentenc X)$ (for which the relation
$\models^{\intnl}$ has been defined), a PKL-structure $\gmA=\pairof{A,\symb{K}^\gmA\ctenten}$ of a relevant signature,  
$a_0\ctenten\in A$ and $U_0\ctenten\in\Pkl{\symb{\scriptstyle K}^\gmA}{A}\cap A$.

For PKL-structures $\gmA$, $\gmB$ of the same signature 
with $\gmB=\pairof{B,\symb{K}^\gmB\ctenten}$ and $\gmB\subseteq\gmA$, 
we define:
\begin{xitemize}
\xitem[PKL-1] 
  {\it $\gmB\prec^{\intnl}_{\calL^{PKL}_{stat}}\gmA$\ \ $\Leftrightarrow$\\[2\jot]
  \qquad $\gmB\models^{\intnl}\varphi(b_0\ctentenc U_0\ctenten)$ if and only if 
  $\gmA\models^{\intnl}\varphi(b_0\ctentenc U_0\ctenten)$\\
  \qquad for all $\calL^{PKL}_{stat}$-formulas $\varphi$ in the signature 
  of the structures with\\
  \qquad $\varphi=\varphi(x_0\ctentenc X_0\ctenten)$, $b_0\ctenten\in B$ and
  $U_0\ctenten\in\Pkl{\symb{\scriptstyle K}^\gmB}{B}\cap B$.}
\end{xitemize}

Finally, we define the internal SDLS for this logic as 
follows:

Suppose that $\kappa$ is a regular cardinal $>\aleph_1$. 
\begin{xitemize}
\item[$\SDLS^{\intnl}_+(\calL^{PKL}_{stat},\LT\kappa)$:]  {\it For any PKL-structure
  $\gmA=\pairof{A,\symb{K}^\gmA\ctenten}$ of  
  countable signature with $\cardof{A}\geq\kappa$ and
  $\cardof{\symb{K}^\gmA}=\kappa$, there are {\it stationarily many}  
  $M\in[A]^{\LT\kappa}$ \st\  
  $\gmA\restr M$ is a PKL-structure and
  $\gmA\restr M\prec^{\intnl}_{\calL^{PKL}_{stat}}\gmA$.  }
\end{xitemize}

The following diagonal reflection characterizes
$\SDLS^{\intnl}_+(\calL^{PKL}_{stat},\LT\kappa)$.
For regular cardinals $\kappa$, $\lambda$ with $\kappa\leq\lambda$, let
\begin{xitemize}
\item[$(\ast)^{\intnl+PKL}_{\LT\kappa,\lambda}$: ] {\it For any countable expansion $\gmA$ 
  of the structure $\pairof{\calH(\lambda),\kappa,\in}$ and any family
  $\seqof{S_a}{a\in\calH(\lambda)}$ \st\ $S_a$ is a stationary subset of
  $\Pkl{\kappa}{\calH(\lambda)}$ for all $a\in\calH(\lambda)$, there are stationarily many 
  $M\in\Pkl{\kappa}{\calH(\lambda)}$ \st\ $\cardof{\kappa\cap M}$ is regular, $\gmA\restr M\prec\gmA$ and
  $S_a\cap \Pkl{\kappa\cap M}{M}\cap M$ is stationary in $\Pkl{\kappa\cap M}{M}$ 
  for all $a\in M$.}
\end{xitemize}

The following \Propof{P-PKL-0} can be proved by a modification of the proof 
of \Propof{P-internal-0}.

\begin{Prop}
  \Label{P-PKL-0}
  For a regular cardinal $\kappa>\aleph_1$, \tfae:\smallskip

  \assert{a} $(\ast)^{\intnl+PKL}_{\LT\kappa,\lambda}$ holds for all 
  regular $\lambda\geq\kappa$.\smallskip

  \assert{b} $\SDLS^{\intnl}_+(\calL^{PKL}_{stat},\LT\kappa)$ holds.\ifextended\else\qed\fi
\end{Prop}

\ifextended
{\tt\noindent\prf Suppose first that $\SDLS^{\intnl}_+(\calL^{PKL}_{stat}, \LT\kappa)$ 
holds. We show that 
$(\ast)^{\intnl+PKL}_{\LT\kappa,\lambda}$ holds for all $\lambda\geq\kappa$. Let 
$\lambda\geq\kappa$ and let  $\tilde{\gmA}$ be a countable expansion of
$\pairof{\calH(\lambda),\kappa,\in}$ where $\kappa$ is the interpretation of
$\symb{K}$ in this PKL-structure, $\seqof{S_a}{a\in\calH(\lambda)}$ a 
sequence of stationary subsets of $\Pkl{\kappa}{\calH(\lambda)}$ and
$\calD\subseteq[\calH(\lambda)]^{\LT\kappa}$ a club.

Let
\begin{xitemize}
\xitemA[PKL-int-3] 
  $\tilde{\gmA}^*=\pairof{\underbrace{\calH(\lambda),\kappa\ctentenc\in}_{\tilde{\gmA}},\symb{\vec{S}}^{\tilde{\gmA}^*}}$
\end{xitemize}
where 
$\symb{\vec{S}}$ is a binary relation symbol and 
\begin{xitemize}
\xitemA[PKL-int-4] 
  $\symb{\vec{S}}^{\tilde{\gmA}^*}=\setof{\pairof{a,s}\in(\calH(\lambda))^2}{s\in S_a}$.
\end{xitemize}
Let $M\in[\calH(\lambda)]^{\LT\kappa}$ be \st\
\begin{xitemize}
\xitemA[PKL-int-5] 
  $M\in\calD$ and 
\xitemA[PKL-int-6] 
  $\tilde{\gmA}^*\restr M\prec^{\intnl}_{\calL^{PKL}_{stat}}\tilde{\gmA}^*$.
\end{xitemize}

By the choice of $\tilde{\gmA}^*$ and \xitemof{PKL-int-6},
$\tilde{\gmA}^*\restr M\models^{\intnl}\forall x\,stat\,X\exists y\,(\symb{\vec{S}}(x,y)\land \forall z\,
(z\varin X\leftrightarrow z\in y))$ holds and hence, 
for all $a\in M$,  $S_a\cap M$ is stationary in $\Pkl{\kappa\cap M}{M}$.
\smallskip

Suppose now that $(\ast)^{\intnl+PKL}_{\LT\kappa,\lambda}$ holds for all
$\lambda\geq\kappa$. Let
$\gmA=\pairof{A,\symb{K}^\gmA\ctenten}$ be a structure in 
countable signature and of cardinality $\geq\kappa$ with $\cardof{\symb{K}^\gmA}=\kappa$ and 
$\calD\subseteq[A]^{\LT\kappa}$ a club. \Wolog, we may assume that $\gmA$ is a relational structure.
We may also assume \wolog\ that $\kappa\subseteq A$ and $\symb{K}^\gmA=\kappa$. 
Let $\lambda$ be a regular cardinal \st\ $\gmA\in\calH(\lambda)$. Note that we 
have in particular $A\subseteq\calH(\lambda)$. 

Let
$\tilde{\gmA}=\pairof{\calH(\lambda),\underbrace{\symb{A}^{\tilde{\gmA}},\kappa\ctenten}_{=\gmA},\in}$ 
where $\symb{A}$ is a unary relation symbol and \vspace{-2ex}\\$\symb{A}^{\tilde{\gmA}}=A$. 

For each $a\in\calH(\lambda)$, let 
\begin{xitemize}
\xitemA[PKL-int-7] 
  $S_a=\left\{\,
  \begin{array}[c]{@{}l}
    \setof{U\in\Pkl{\kappa}{\calH(\lambda)}}{U\cap A\in A,\\[\jot]
    \phantom{{}U\in[\calH(\lambda)]^{\LT\kappa}:\ \ \ \ }
            \gmA\models^{\intnl}\psi(\variables{a}{m-1},\variables{U}{n-1},U\cap A)},\\[\jot]
    \qquad\quad
    \mbox{if }\psi=\psi(\variables{x}{m-1},\variables{Y}{n-1},X)\mbox{ is an }
    \calL^{PKL}_{stat}\mbox{-formula}\\
    \qquad\quad\mbox{in the signature of }\gmA,\,\variables{a}{m-1}\in 
    A,\,\variables{U}{n-1}\in\Pkl{\kappa}{A}\cap A,\\
    \qquad\quad\gmA\models^{\intnl}stat\,X\,\psi(\variables{a}{m-1},\variables{U}{n-1},X)\mbox{ and}\\
    \qquad\quad a=\pairof{\psi,\variables{a}{m-1},\variables{U}{n-1}};\\[2\jot]

    \Pkl{\kappa}{\calH(\lambda)},\\\qquad\quad\mbox{otherwise.}
  \end{array}
  \right.$
\end{xitemize}

We expand $\tilde{\gmA}$ further by adding the relation
$\setof{\pairof{u,a}\in(\calH(\lambda))^2}{u\in S_a}$ to it as the interpretation of 
the binary relation symbol $\symb{\vec{S}}$. For simplicity, we shall also call 
this expanded structure $\tilde{\gmA}$. 

Let
\begin{xitemize}
\xitemA[PKL-int-8] 
  $\tilde{\calD}=\setof{U\in[\calH(\lambda)]^{\LT\kappa}}{U\cap A\in\calD}$. 
\end{xitemize}
$\tilde{\calD}$ contains a club in $[\calH(\lambda)]^{\LT\kappa}$.

By $(\ast)^{\intnl+PKL}_{\LT\kappa,\lambda}$, there is an 
$M\in[\calH(\lambda)]^{\LT\kappa}$ \st\
\begin{xitemize}
\xitemA[PKL-int-9] 
  $M\in\tilde{\calD}$,
\xitemA[PKL-int-10] 
  $\tilde{\gmA}\restr M\prec\tilde{\gmA}$\ \ and 
\xitemA[PKL-int-11] 
  $S_a\cap[M]^{\LT\kappa}\cap M$ is stationary in $[M]^{\LT\kappa}$ 
  for all $a\in M$. 
\end{xitemize}

Let $\tilde{\gmB}=\tilde{\gmA}\restr M$ and 
let $B=\symb{A}^{\tilde\gmB}=A \cap M$ and $\gmB=\gmA\restr B$.  
By \xitemof{PKL-int-10}
\begin{xitemize}
\xitemA[PKL-int-11-0] 
  $\kappa\cap M=\kappa\cap\tilde{B}=\kappa\cap B$. 
\end{xitemize}
$B\in\calD$ by \xitemof{PKL-int-9} and the definition \xitemof{PKL-int-8} of
$\tilde{\calD}$.


By the elementarity $\tilde{\gmB}\prec\tilde{\gmA}$ \xitemof{PKL-int-10}, the following Claim implies
$\gmB\prec^{\intnl}_{\calL^{PKL}_{stat}}\gmA$.

\begin{Claim} For any $\calL^{PKL}_{stat}$-formula
  $\varphi(\variables{x}{m-1},\variables{Y}{n-1})$ in the signature of the 
  structures $\gmA$, $\variables{a}{m-1}\in B$ and
  $\variables{U}{n-1}\in\Pkl{\kappa\cap B}{B}$, we have
  \begin{xitemize}
  \xitemA[PKL-int-12] 
    $\tilde{\gmB}\modelof{\gmA\models^{\intnl}\varphi(a_0\ctentenc U_0\ctenten)}$\ \ 
    $\Leftrightarrow$\ \
    $\gmB\models^{\intnl}\varphi(a_0\ctentenc U_0\ctenten)$.
  \end{xitemize}
\end{Claim}
\prfofClaim
By induction on $\varphi$. The crucial step in the induction is 
when $\varphi$ is of the form $stat\,X\psi$ and \xitemof{PKL-int-12} holds for $\psi$:

Suppose first that
$\tilde{\gmB}\modelof{\gmA\models^{\intnl}\varphi(a_0\ctentenc U_0\ctenten)}$ holds.
Then, by elementarity and by the definition of $\tilde{\gmA}$, we have
$\gmA\models^{\intnl} stat\,X\psi(a_0\ctentenc\variables{U}{n-1},X)$. Thus, letting
$a=\pairof{\varphi,\variables{a}{m-1},\variables{U}{n-1}}$, we have $a\in\tilde{B}$ and 
\begin{xitemize}
\xitemA[PKL-int-14] 
  $S_a=\setof{U\in\Pkl{\kappa}{\calH(\lambda)}}{
  \begin{array}[t]{@{}l}
    U\cap A\in A,\\\gmA\models^{\intnl}\varphi(\variables{a}{m-1},\variables{U}{n-1}, U\cap A)}
  \end{array}$ 
\end{xitemize}
by the definition \xitemof{PKL-int-7} of $S_a$.  

By 
\xitemof{PKL-int-11}, $S_a\,\cap\,\Pkl{\kappa\cap\tilde{B}}{\tilde{B}}\,\cap\,\tilde{B}$ 
is stationary in $\Pkl{\kappa\cap\tilde{B}}{\tilde{B}}$.  
It follows that
\begin{xitemize}
\xitemA[PKL-int-15] 
  $\setof{U\cap B}{U\in S_a\,\cap\,\Pkl{\kappa\cap\tilde{B}}{\tilde{B}}\,\cap\,\tilde{B}}$\\[\jot]
  $=\setof{U\cap B}{{}
  \begin{array}[t]{@{}l}
    U\cap B\in\Pkl{\kappa\cap B}{B}\cap B,\,\\
    \tilde{\gmB}\modelof{\gmA\models^{\intnl}\psi(\variables{a}{m-1},\variables{U}{n-1},U\cap B)}}
  \end{array}
  $\\
  \mbox{}\hfill(by elementarity, \xitemof{PKL-int-11-0} and \xitemof{PKL-int-14})\\
  $\subseteq\setof{V\in\Pkl{\kappa\cap B}{B}\cap  B}{\gmB\models^{\intnl}
    \psi(\variables{a}{m-1},\variables{U}{n-1},V)}$\\
  \mbox{}\hfill(by induction hypothesis)
\end{xitemize}
is stationary.
Thus 
$\gmB\models^{\intnl} stat\,X\,\psi(\variables{a}{m-1},\variables{U}{n-1},X)$ holds, 
that is, 
$\gmB\models^{\intnl} \varphi(\variables{a}{m-1},\variables{U}{n-1})$.

Suppose now that
$\tilde{\gmB}\not\modelof{\gmA\models^{\intnl}\varphi(\variables{a}{m-1},\variables{U}{n-1})}$. 
Then we have
\begin{xitemize}
\xitemA[PKL-int-16] 
  $\tilde{\gmB}\modelof{
  \begin{array}[t]{@{}l}
    \mbox{there is a club }\calC\subseteq[\symb{A}]^{\LT\kappa}\mbox{ \st\ }\\
    \gmA\models^{\intnl}\neg\psi(a_0\ctentenc\variables{U}{n-1},x)
    \mbox{ for all }x\in\calC}.
  \end{array}$
\end{xitemize}

By elementarity, there is a $\calC_0\in\tilde{B}$ \st\ $\calC_0$ is a club $\subseteq$ 
$[A]^{\LT\kappa}$ and 
\begin{xitemize}
\xitemA[PKL-int-17] 
  $\tilde{\gmB}\modelof{\gmA\models^{\intnl}\neg\psi(\variables{a}{m-1},\variables{U}{n-1},V)}$
  for all $V\in\calC_0\cap B$.
\end{xitemize}

By induction hypothesis, it follows that
\begin{xitemize}
\xitemA[PKL-int-17-0] 
  $\setof{U\in\Pkl{\kappa\cap B}{B}\cap B}{
  \calB\models^{\intnl}\psi(a_0\ctentenc U_0\ctentenc U)}\cap\calC_0=\emptyset$. 
\end{xitemize}
Thus 
$\gmB\not\models^{\intnl}stat\,X\,\psi(a_0\ctentenc\variables{U}{n-1},X)$, i.e.\ 
$\gmB\not\models^{\intnl}\varphi(a_0\ctentenc U_0\ctenten)$. \\\qedofClaim\\
}
\qedofProp\qedskip
\fi

For a regular cardinal $\kappa$ and a cardinal $\lambda\geq\kappa$, $\calS\subseteq\Pkl{}{}$ is 
said to be {\it $2$-stationary} if, for any stationary $\calT\subseteq\Pkl{}{}$, 
there is an $a\in\calS$ \st\ $\cardof{\kappa\cap a}$ is a regular uncountable 
cardinal and $\calT\cap\Pkl{\kappa\cap a}{a}$ is stationary in
$\Pkl{\kappa\cap a}{a}$. A regular cardinal $\kappa$ has the {\it$2$-stationarity 
  property} if $\Pkl{}{}$ is $2$-stationary (as a subset of itself) for all $\lambda\geq\kappa$. 
More generally, we can define $\alpha$-stationarity for $\alpha\leq\kappa$ and 
show that these generalized stationarities are compatible with $\kappa$ being 
continuum (\cite{brickhill-fuchino-sakai}). 

\begin{Lemma}
  \Label{P-PKL-1} For a regular uncountable $\kappa$,
  $\SDLS^{\intnl}_+(\calL^{PKL}_{stat},\LT\kappa)$ implies that $\kappa$ has 
  the $2$-stationarity property.
\end{Lemma}
\prf The property \assertof{a} in \Propof{P-PKL-0} is a strengthening of the 
$2$-stationarity of $\kappa$. \qedofLemma

\begin{Lemma}
  \Label{P-PKL-2} Suppose that $\kappa$ is a regular uncountable cardinal.

  \assert{1} If $\kappa$ has the $2$-stationarity property, then $\kappa$ is a limit cardinal.\smallskip

  \assert{2} For any $\lambda\geq\kappa$, $2$-stationary $\calS\subseteq\Pkl{}{}$, and any 
  stationary $\calT\subseteq\Pkl{}{}$, there are stationarily many $r\in\calS$ 
  \st\ $\calT\cap\Pkl{\kappa\cap r}{r}$ is stationary. \smallskip

  \assert{3} If $\kappa$ is $2$-stationary then $\kappa$ is a weakly Mahlo cardinal.

\end{Lemma}
\prf \assertof{1}: Suppose that $\kappa=\mu^+$. Then
$\calC=\setof{a\in\Pkl{}{}}{\cardof{a}=\mu}$ is a club and hence 
stationary. But, for 
any $r\in\Pkl{}{}$, $\cardof{\kappa\cap r}\leq\mu$ and hence 
$\calC\cap\Pkl{\kappa\cap r}{r}=\emptyset$. 
Thus $\kappa$ is not $2$-stationary. 
\smallskip

\assertof{2}: Suppose that $\calS\subseteq\Pkl{}{}$ is $2$-stationary 
and $\calT\subseteq\Pkl{}{}$ is stationary. Let
$\calC\subseteq\Pkl{}{}$ be a club. We have to show that there is
$r\in\calS\cap\calC$ \st\ $\calT\cap\Pkl{\kappa\cap r}{r}$ is stationary. 

Let $\mapping{f}{\fnsp{\omega>}{\lambda}}{\lambda}$ be \st\
\begin{xitemize}
\xitem[] 
  $\calC_f=\setof{a\in\Pkl{}{}}{\kappa\cap a\in\kappa\mbox{ and }
    a\mbox{ is closed under }f}\subseteq\calC$. 
\end{xitemize}
Since $\calC_f$ is a club, $\calT\cap\calC_f$ is stationary. Let $r\in\calS$ be 
\st\ $\cardof{\kappa\cap r}$ is regular and 
$(\calT\cap\calC_f)\cap\Pkl{\kappa\cap r}{r}$ is stationary in $\Pkl{\kappa\cap r}{r}$.  We have
\begin{Claim}
  \Label{Cl-a-stat-0}
  $\kappa\cap r\in \kappa$.  
\end{Claim}
\prfofClaim
Otherwise, there is an $\alpha\in\sup(\kappa\cap r)\setminus r$. Let
$\calY=\setof{b\in\Pkl{\kappa\cap r}{r}}{\sup(\kappa\cap b)>\alpha}$. $\calY$ 
is club in $\Pkl{\kappa\cap r}{r}$ but $\calY\cap\calC_f=\emptyset$. 
Thus $(\calT\cap\calC_f)\cap\Pkl{\kappa\cap r}{r}\cap\calY=\emptyset$. This is a 
contradiction to the stationarity of
$(\calT\cap\calC_f)\cap\Pkl{\kappa\cap r}{r}$ in $\Pkl{\kappa\cap r}{r}$. 
\qedofClaim
\begin{Claim}
  \Label{Cl-a-stat-1} $r$ is closed under $f$.
\end{Claim}
\prfofClaim
This is clear since $\calT\cap\calC_f$ is cofinal in $\Pkl{\kappa\cap r}{r}$ 
\wrt\ $\subseteq$) and elements of $\calT\cap\calC_f$ are closed under $f$. 
\qedofClaim\qedskip

From the Claims above, it follows that $r\in\calC_f\subseteq\calC$ is as desired. \smallskip

\assertof{3}: Let
$\calT=\setof{a\in\Pkl{}{}}{\kappa\cap a\in\kappa}$. $\calT$ is a club and 
hence stationary. Let $r\in\Pkl{}{}$ be \st\ $\cardof{\kappa\cap r}$ is regular 
and 
\begin{xitemize}
\xitem[a-stat-2] 
  $\calT\cap\Pkl{\kappa\cap r}{r}$ is stationary in
  $\Pkl{\kappa\cap r}{r}$.
\end{xitemize}
Similarly to \Claimof{Cl-a-stat-0},  
we have $\kappa\cap r\in\kappa$.  
\begin{Claim}
  \Label{Cl-a-stat-2}
  $\kappa\cap r$ is a cardinal.
\end{Claim}
\prfofClaim Otherwise there is $\mu<\kappa\cap r$ \st\
$\cardof{\kappa\cap r}=\mu$. But then the set 
$\setof{a\in\Pkl{\kappa\cap r}{r}}{\sup(a\cap\kappa)\geq\mu}$ is a club 
in $\Pkl{\kappa\cap r}{r}$ disjoint from $\calT$. \qedofClaim
\begin{Claim}
  \Label{Cl-a-stat-3}
  $\kappa\cap r$ is a regular cardinal.
\end{Claim}
\prfofClaim Otherwise there is an $s\subseteq\kappa\cap r$ cofinal in $\kappa\cap r$
with $\cardof{s}<\kappa\cap r$. 
But then the set 
$\setof{a\in\Pkl{\kappa\cap r}{r}}{\sup(a\cap\kappa)\supseteq s}$ is a club 
in $\Pkl{\kappa\cap r}{r}$ disjoint from $\calT$.\\\qedofClaim
\qedskip

Since there are stationarily many $r$ with \xitemof{a-stat-2} by \assertof{2}, it 
follows from \Claimabove\ that 
$\kappa$ is weakly Mahlo.  \qedofLemma
\qedskip

Note that we can continue in the proof of \assertof{3} above to show that 
$\kappa$ is weakly hyper Mahlo, weakly hyper hyper Mahlo. etc.

\begin{Cor}
  \Label{Cl-a-stat-4} $\SDLS^{\intnl}_+(\calL^{PKL}_{stat},\LT\kappa)$ implies that 
$\kappa$ is weakly Mahlo, weakly hyper Mahlo, etc.
\end{Cor}
\prf By \Lemmaof{P-PKL-1}, \Lemmaof{P-PKL-2},\,\assertof{3} and the remark below 
it. \qedofCor

{
\begin{Thm}
  \Label{P-lt-conti-3:} Suppose that $\kappa$ is a generically supercompact 
  cardinal by
  $\mu$-cc \pos\ for some $\mu<\kappa$. Then $(\ast)^{\intnl+PKL}_{\LT\kappa,\lambda}$ holds for all 
  regular $\lambda\geq\kappa$. This means that $\SDLS^{int}_+(\calL^{PKL}_{stat},\LT\kappa)$ holds 
  by \Propof{P-PKL-0}.
\end{Thm}

The proof of \Thmof{P-lt-conti-3:} can be done analogously to the proof of 
\Thmof{P-lt-conti-3} noting \Lemmaof{P-lt-conti-4} and 
\begin{Lemma}
  \Label{P-lt-conti-4:}
  Suppose that $M$ is an inner model in $\uniV$, $\lambda$ an ordinal and
  $\mu\leq\lambda$ a regular cardinal in $\uniV$. 

  If $\poP$ is a $\mu$-cc \po\ and $\calS\subseteq\Pkl{\mu}{\lambda}$ is 
  stationary in $\Pkl{\mu}{\lambda}$, then
  $\forces{\poP}{\check{\calS}\xmbox{ is stationary in }\Pkl{\mu}{\lambda}}$.
  \ifextended\else\qed\fi
\end{Lemma}
\ifextended{\tt
  \prf
  Suppose that $\utilde{\calC}$ is a $\poP$-name with
  $\forces{\poP}{\utilde{\calC}\mbox{ is a club in }\Pkl{\mu}{\lambda}}$.

  In $\uniV$, 
  let
  $\calC=\setof{C\in\Pkl{\mu}{\lambda}}{\forces{\poP}{\check{C}\varin\utilde{\calC}}}$. Then
  $\calC$ is club by the $\mu$-cc of $\poP$. 
  Hence $\calS\cap\calC\not=\emptyset$. Since
  $\forces{\poP}{\check{\calC}\subseteq\utilde{\calC}}$, it follows that
  $\forces{\poP}{\check{\calS}\cap\utilde{\calC}\not\equiv\emptyset}$. 
  \qedofLemma\qedskip
}\fi

\ifextended{\tt
\noindent
\prfof{\Thmof{P-lt-conti-3:}} For a regular $\lambda\geq\kappa$, let
$\lambda^*=\cardof{\calH(\lambda)}$ and $\lambda^{**}=(2^{\lambda^*})^+$. 

Suppose that 
$\gmA=\pairof{\calH(\lambda),\in, \kappa\ctenten}$ is a countable expansion of the structure
$\pairof{\calH(\lambda),\in,\kappa}$, $\seqof{S_a}{a\in\calH(\lambda)}$ a sequence of 
stationary subsets of $\Pkl{}{\calH(\lambda)}$ and 
$\calC\subseteq\Pkl{}{\calH(\lambda)}$ a club.

\Wolog, we may assume that $\gmA$ contains the relation 
$\setof{\pairof{b,a}\in(\calH(\lambda))^2}{b\in S_a}$ coding the sequence
$\seqof{S_a}{a\in\calH(\lambda)}$ in its structure. 

Let $\mapping{\iota}{\calH(\lambda)}{\lambda^*}$ be a bijection \st\
$\iota\restr\kappa=id_\kappa$ 
and let
$\gmA^*=\pairof{\lambda^*,E^*,\kappa\ctenten}$ be a copy of $\gmA$ translated 
by $\iota$. 

Let 
$\seqof{S^*_\alpha}{\alpha\in\lambda^*}$ be \st\ 
\begin{xitemize}
\xitemA[lt-conti-2-3:] 
  $S^*_\alpha=\iota\imageof S_{\iota^{-1}(\alpha)}$ for $\alpha\in\lambda^*$ and 
\end{xitemize}
\begin{xitemize}
\xitemA[lt-conti-2-4:] 
  $\calC^*=\setof{\iota\imageof X}{X\in\calC}$. 
\end{xitemize}
Thus, we have 
\begin{xitemize}
\xitemA[lt-conti-2-5:] 
  $EXT_{E^*}(S^*_\alpha)$ is a stationary subset of $\Pkl{}{\lambda^*}$ for 
  all $\alpha\in\lambda^*$ and 
\end{xitemize}
\begin{xitemize}
\xitemA[lt-conti-2-6:] 
  $\calC^*$ is a club in $\Pkl{}{\lambda^*}$
\end{xitemize}
where, for $S\subseteq\lambda^*$,
$EXT_{E^*}(S)=\setof{\setof{\beta\in\lambda^*}{\beta\mathrel{E^*}\alpha}}{\mbox{ for }\alpha\in S}$
is the set of all extents of elements of $S$ \wrt\ the element relation $E^*$.

It is enough to show that there is an $X\in\calC^*$ \st\ 
\begin{xitemize}
\item[\wassertof{2$'$}] $\gmA^*\restr X\prec\gmA^*$, $\kappa\cap X<\kappa$ and
\item[\wassertof{3$''$}] $EXT_{E^*}(S^*_\alpha\cap X)$ is stationary in 
$\Pkl{\kappa\cap X}{X}$ for all $\alpha\in X$.  
\end{xitemize}

Let $\poP$ be a $\mu$-cc \po\ and $\genG$ a $(\uniV,\poP)$-generic 
filter \st\ there are transitive $M\subseteq\uniV[\genG]$ and
$\elembed{j}{\uniV}{M}$ with
\begin{xitemize}
\xitemA[lt-conti-3:] 
  $\kappa=crit(j)$,
\xitemA[lt-conti-4:] 
  $j(\kappa)>\lambda^{**}$ and
\xitemA[lt-conti-5:] 
  $j\imageof\lambda^{**}\in M$. 
\end{xitemize}
Let $X=j\imageof\lambda^*$. $X\in M$ by \xitemof{lt-conti-5:} and \Lemmaof{L-lt-conti-0},\,\assertof{1}. 
Thus, we also have $j\restr\lambda^*\in M$ and hence   
\begin{xitemize}
\xitemA[costat-2-3:] 
  $M\models\cardof{X}=\cardof{\lambda^*}< j(\kappa)$. 
\end{xitemize}

Since $\cardof{\calC^*}<\lambda^{**}$, we have $j\imageof\calC^*\in M$ by 
\Lemmaof{L-lt-conti-0},\,\assertof{1}. 
$M\modelof{\cardof{j\imageof\calC^*}\leq\cardof{j(\calC^*)}<j(\kappa)}$.  
Also we have
$M\modelof{j\imageof\calC^*\mbox{ is directed and }\bigcup(j\imageof\calC^*)=X}$. 
Since $M\modelof{\calC^*\xmbox{ is a club in }\Pkl{j(\kappa)}{j(\lambda^*)}}$ by 
elementarity and \xitemof{lt-conti-2-6:}, 
it follows that 
\begin{xitemize}
\xitemA[costat-2-4:] 
  $X\in j(\calC^*)$.
\end{xitemize}

Let
$\seqof{\tilde{S}^*_\alpha}{\alpha<j(\lambda^*)}=j(\seqof{S^*_\alpha}{\alpha<\lambda^*})$. 

By elementarity and applying Vaught's test, we can show
\begin{xitemize}
\xitemA[] 
  $M\models j(\gmA^*)\restr X\prec j(\gmA^*)$.
\end{xitemize}
For $\alpha\in\lambda^*$, since
$\tilde{S}^*_{j(\alpha)}=j(S^*_\alpha)\supseteq j\imageof S^*_\alpha$, we have 
\begin{xitemize}
\xitemA[costat-2-5:] 
  $M\modelof{EXT_{j(E^*)}(\tilde{S}^*_{j(\alpha)}\cap X)
  \supseteq EXT_{j(E^*)}(j\imageof S^*_\alpha)}$.
\end{xitemize}

By \xitemof{lt-conti-2-5:} and since $\poP$ is $\mu$-cc, we have\\
$\uniV[\genG]\modelof{EXT_{E^*}(S^*_\alpha)\mbox{ is stationary in }\Pkl{}{\lambda^*}}$
by \Lemmaof{P-lt-conti-4:}
. \\
Since $EXT_{j(E^*)}(j\imageof S^*_\alpha)$ is a translation of
$EXT_{E^*}(S^*_\alpha)$ induced by $j\restr\lambda^*$, it follows that 
\begin{xitemize}
\xitemA[costat-2-6:] 
  $\uniV[\genG]\modelof{EXT_{j(E^*)}(j\imageof S^*_\alpha)
  \mbox{ is stationary in }\Pkl{}{X}}$.
\end{xitemize}
By \Lemmaof{P-lt-conti-4}, it follows that 
\begin{xitemize}
\xitemA[costat-2-7:] 
  $M\modelof{EXT_{j(E^*)}(j\imageof S^*_\alpha)
  \mbox{ is stationary in }\Pkl{}{X}}$.
\end{xitemize}
Noting that $j(\kappa)\cap X=\kappa$, 
\begin{xitemize}
\xitemA[] 
  $M\modelof{\exists X\in j(\calC^*)\,\big(\,
  \begin{array}[t]{@{}l}
    j(\gmA^*)\restr X\prec j(\gmA^*)  \land EXT_{j(E^*)}(\tilde{S}^*_\xi\cap 
    X)\mbox{ is}
    \\\mbox{a stationary subset of }\Pkl{j(\kappa)\cap X}{X}\mbox{ for all }\xi\in X\big)}.
  \end{array}
  $
\end{xitemize}
By elementarity, it follows that 
\begin{xitemize}
\xitemA[] 
  $\uniV\modelof{\exists X\in\calC^*\,\big(\,
  \begin{array}[t]{@{}l}
    \gmA^*\restr X\prec \gmA^*\  \land\ EXT_{E^*}(S^*_\xi\cap X)
    \mbox{ is}\\\mbox{a stationary subset of }\Pkl{\kappa\cap X}{X}\mbox{ for all }\xi\in X\big)}.
  \end{array}
  $
\end{xitemize}
\mbox{}\vspace{-3ex}

\mbox{}\qedof{\Thmof{P-lt-conti-3:}}
}\fi
}

\section{Laver-generic large cardinals}
\Label{laver}


For a cardinal $\kappa$ and a class $\calP$ of \pos, we call $\kappa$ a 
{\it Laver-generically supercompact for $\calP$} if, for any $\lambda\geq\kappa$ and any
$\poP\in\calP$, there are a \po\ $\poQ\in\calP$ 
with $\poP\circleq\poQ$ and $(\uniV,\poQ)$-generic filter $\genH$ \st\ there 
are an inner model
$M\subseteq\uniV[\genH]$ and a class $j\subseteq\uniV[\genH]$ with 
\begin{xitemize}
\xitem[laver-a] $\elembed{j}{\uniV}{M}$, 
\xitem[laver-0] $\crit(j)=\kappa$, $j(\kappa)
  >\lambda$,
\xitem[laver-0-0] $\poP$, $\genH\in M$ and 
\xitem[laver-1] $j\imageof\lambda\in M$.
\end{xitemize}

$\kappa$ is {\it Laver-generically superhuge} ({\it Laver-generically super almost-huge} 
resp.) {\it for $\calP$} if $\kappa$ satisfies 
the definition of Laver-generic supercompactness for $\calP$ with  
\xitemof{laver-1} replaced by 
\begin{xitemize}
\xitemd[laver-1]{'}
  $j\imageof j(\kappa)\in M$  ($j\imageof\mu\in M$ for all $\mu<j(\kappa)$ resp.).
\end{xitemize}

$\kappa$ is {\it tightly Laver-generically supercompact} ({\it tightly 
  Laver-generically superhuge, tightly Laver-generically super almost-huge}, 
resp.) if the definition of {Laver-generically supercompact} ({Laver-generically 
  superhuge, Laver-generically super almost-huge}, resp.) holds with 
\xitemof{laver-0} replaced by 
\begin{xitemize}
\xitemd[laver-0]{'} 
  $\crit(j)=\kappa$, $j(\kappa)=\cardof{\poQ}>\lambda$.
\end{xitemize}

All consistency proofs of the existence of Laver-generic very large cardinals we know actually 
show the existence of tightly Laver-generic very large cardinals (see 
the proof of \Thmof{T-laver-a}).

The following is clear by definition.
\begin{Lemma}
  \Label{L-laver-0} Suppose that $\calP$ is a class of \pos. \wassertof{1} 
  If $\kappa$ is Laver-generically superhuge for $\calP$ then $\kappa$ is 
  Laver-generically super almost-huge for $\calP$. If $\kappa$ is 
  Laver-generically super almost-huge for $\poP$ then $\kappa$ is 
  Laver-generically supercompact for $\calP$. 
  If   $\kappa$ is Laver-generically supercompact for $\calP$ then $\kappa$ is 
  generically supercompact by $\calP$. \smallskip

  \assert{2} If $\kappa$ is generically supercompact by $\calP$ then
  $\kappa$ is generically measurable by some $\poP\in\calP$.\smallskip

  \assert{3} If $\kappa$ is tightly Laver-generically supercompact (super 
  almost-huge, superhuge, resp.) for $\calP$ then $\kappa$ is 
  Laver-generically supercompact (super almost-huge, huge, resp.) for $\calP$.\smallskip

  \assert{4} If $\kappa$ is Laver-generically supercompact for $\calP$ then, for 
  any $\poP\in\calP$, there is $\poQ\in\calP$ with $\poP\circleq\poQ$ 
  \st\ $\kappa$ is generically measurable by $\poQ$.\qed
\end{Lemma}

The following \Thmof{T-laver-a} can be still improved and extended: among other things the 
assumption of the consistency of the existence of a  superhuge cardinal in 
\assertof{2} and \assertof{3} can be reduced to that of the existence of a super 
almost-huge cardinal. 
We shall discuss more about these improvements in \cite{IV}.

\begin{Thm}
  \Label{T-laver-a}
  \assert{1} Suppose that \ZFC\ $+$ ``there exists a supercompact cardinal (super 
  almost-huge cardinal, superhuge cardinal, resp.)'' is 
  consistent. Then \ZFC\ $+$ ``there exists a tightly Laver-generically supercompact cardinal 
  (super almost-huge cardinal, superhuge cardinal, resp.)
  for 
  $\sigma$-closed \pos'' is consistent as well.\smallskip
  
  \assert{2}  Suppose that \ZFC\ $+$ ``there exists a superhuge cardinal'' is  
  consistent. Then \ZFC\ $+$ ``there exists a tightly Laver-generically super almost-huge
  cardinal for proper \pos'' is consistent as well.\smallskip

  \assert{3} Suppose that \ZFC\ $+$ ``there exists a supercompact cardinal (superhuge cardinal, resp.)'' is 
  consistent. Then \ZFC\ $+$ ``there exists a tightly Laver-generically supercompact 
  cardinal (super almost-huge cardinal, resp.) for ccc \pos'' is consistent as well.\smallskip
\end{Thm}
\prf \assertof{1}: Suppose that $\kappa$ is a supercompact cardinal (the case 
with a super almost-huge cardinal or a superhuge cardinal can be treated 
similarly).

We show that $\poC=\Col(\omega_1,\kappa)$ forces that $\kappa$ is 
Laver-generically supercompact cardinal for $\sigma$-closed \pos.

Let $\genG_0$ be a $(\uniV,\poC)$-generic filter and $\poP$ be a $\sigma$-closed 
\po\ in $\uniV[\genG_0]$ and $\lambda\geq\kappa$. We may assume that
$\uniV[\genG_0]\modelof{\lambda\geq\cardof{\poP}}$. 

Let $\elembed{j}{\uniV}{M\subseteq\uniV}$ be \st\ $\crit(j)=\kappa$, 
$j(\kappa)>\lambda$ and $[M]^\lambda\subseteq M$. 
Note that, by the $\sigma$-closedness of $\poC$ and closedness property of $M$, we 
have $\Col(\omega_1,j(\kappa))^\uniV=\Col(\omega_1,j(\kappa))^{\uniV[\genG_0]}=\Col(\omega_1,j(\kappa))^M$.
Thus we denote this \po\ simply by $\Col(\omega_1,j(\kappa))$. 

We have 
\begin{xitemize}
\xitem[laver-1-16] 
  $M\modelof{j(\poC)=\Col(\omega_1,j(\kappa))\sim\poC\times\Col(\omega_1,j(\kappa)\setminus\kappa)
  \sim\poC\times\Col(\omega_1,j(\kappa))}$. 
\end{xitemize}
By \Corof{Cor-col-0} we also have 
\begin{xitemize}
\xitem[laver-1-16-0] 
  $\uniV[\genG_0]\modelof{\poP\circleq\Col(\omega_1, j(\kappa))}$. 
\end{xitemize}
Note that $\uniV[\genG_0]\modelof{\cardof{\Col(\omega_1,j(\kappa))}=j(\kappa)}$.

Let $\genH$ be a $(\uniV[\genG_0], \Col(\omega_1,j(\kappa)))$-generic filter and 
let $\genH^*$ be the $(\uniV,j(\poC))$-generic filter extending $\genG_0$ which 
corresponds to $\genG_0*\genH$ via \xitemof{laver-1-16}. We have
$\uniV[\genG_0][\genH]=\uniV[\genH^*]$.
Let 
\begin{xitemize}
\xitem[laver-1-17] 
  $\elembed{j^*}{\uniV[\genG_0]}{M[\genG^*]\subseteq\uniV[\genG^*]}$;
  $\utildea[\genG_0]\mapsto j(\utildea)[\genH^*]$.
\end{xitemize}
Then we have $\crit(j^*)=\crit(j)=\kappa$, $j^*(\kappa)>\lambda$,
$j^*\imageof\lambda=j\imageof\lambda\in M\subseteq M[\genH^*]$ and $M[\genH^*]$ 
is an inner model of $\uniV[\genG_0][\genH]$. We also have $\genH\in M[\genH^*]$ 
since $\utilde{\genH}^*\in M$. 
Since $\lambda$ can be taken 
arbitrarily large, this shows that $\kappa$ in $\uniV[\genG_0]$ is 
tightly Laver-generically supercompact for $\sigma$-closed forcing. \smallskip

\assertof{2}: Suppose that $\kappa$ is a superhuge cardinal. Then, by 
\cite{corazza}, there is a super almost-huge Laver-function
$\mapping{f}{\kappa}{\uniV_\kappa}$. Let
$\seqof{\poP_\alpha,\utildepoQ_\beta}{\alpha\leq\kappa,\beta<\kappa}$ be a 
CS-iteration of proper \pos\ \st\
\begin{xitemize}
\xitem[laver-1-18] 
  $\utildepoQ_\beta=\left\{
  \begin{array}{@{}l}
    f(\beta),\quad\ \ \mbox{if }f(\alpha)\mbox{ is a }\poP_\beta\mbox{-name and }
    \forces{\poP_\beta}{f(\alpha)\mbox{ is a proper \po}};\\
    \poP_\alpha\mbox{-name of the trivial \po},\quad \mbox{otherwise.}
  \end{array}
  \right.$
\end{xitemize}
Let $\genG_0$ be a $(\uniV,\poP_\kappa)$-generic filter.
We show that, in $\uniV[\genG_0]$, $\kappa$ is a tightly Laver-generically super 
almost-huge cardinal for proper \pos.

Working in $\uniV[\genG_0]$, suppose that $\poP$ is a proper \po\ and
$\lambda\geq\kappa$. Let $\utildepoP$ be be $\poP_\kappa$-name of $\poP$.

Back in $\uniV$, let $\elembed{j}{\uniV}{M\subseteq\uniV}$ be \st\
\begin{xitemize}
\xitem[laver-1-19] $\crit(j)=\kappa$,
\xitem[laver-1-20] $j(\kappa)>\lambda$,
\xitem[laver-1-21] $[M]^{\LT j(\kappa)}\subseteq M$, and
\xitem[laver-1-22] $f(\kappa)=\utildepoP$.
\end{xitemize}
By elementarity, we have
\begin{xitemize}
\xitem[laver-1-23] 
  $M\modelof{
      \begin{array}[t]{@{}l}
        j(\seqof{\poP_\alpha}{\alpha\leq\kappa})\mbox{ and }j(\seqof{\utildepoQ_\beta}{\beta<\kappa})
        \mbox{ make up a CS-iteration of }\\
        \mbox{proper \pos\ of length }j(\kappa)\mbox{, and each name in the 
          sequence }\\
        j(\seqof{\utildepoQ_\beta}{\beta<\kappa})\mbox{ is of size }<j(\kappa)}.
      \end{array}$
\end{xitemize}

By \xitemof{laver-1-21}, it follows that the statement in \xitemof{laver-1-23} 
also holds in $\uniV$. That is,
\begin{xitemize}
\xitem[laver-1-23-0] 
  $V\modelof{
      \begin{array}[t]{@{}l}
        j(\seqof{\poP_\alpha}{\alpha\leq\kappa})\mbox{ and }j(\seqof{\utildepoQ_\beta}{\beta<\kappa})
        \mbox{ make up a CS-iteration of }\\
        \mbox{proper \pos\ of length }j(\kappa)\mbox{, and each name in the 
          sequence }\\
        j(\seqof{\utildepoQ_\beta}{\beta<\kappa})\mbox{ is of size }<j(\kappa)}.
      \end{array}$
\end{xitemize}

In $\uniV$, let
\begin{xitemize}
\xitem[laver-1-24] 
  $\seqof{\poP^*_\alpha}{\alpha\leq j(\kappa)}=j(\seqof{\poP_\alpha}{\alpha\leq\kappa})$
  and
\xitem[laver-1-25] 
  $\seqof{\utildepoQ^*_\beta}{\beta<j(\kappa)}=j(\seqof{\utildepoQ_\beta}{\beta<\kappa})$.
\end{xitemize}
By elementarity of $j$ and since $\crit(j)=\kappa$, we have 
$\poP_\alpha=\poP^*_\alpha$ for all $\alpha\leq\kappa$ and
$\utildepoQ^*_\beta=\utildepoQ_\beta$ for all $\beta<\kappa$. 
In $\uniV[\genG_0]$, let $\poQ=j(\poP_\kappa)/\genG_0$.

By 
\xitemof{laver-1-23-0}, 
\xitemof{laver-1-22}, 
\xitemof{laver-1-24}, \xitemof{laver-1-25}  and Factor Lemma, $\poQ$ is a proper 
\po\ of cardinality $j(\kappa)$ and $\poP\circleq\poQ$.

Let $\genH$ be a $(\uniV[\genG_0],\poQ)$-generic filter and 
\begin{xitemize}
\xitem[laver-1-26] $\elembed{\tilde{j}}{\uniV[\genG_0]}{M[\genG_0\ast\genH]\subseteq\uniV[\genG_0][\genH]}$;
  $\utildea[\genG_0]\mapsto j(\utildea)[\genG_0\ast\genH]$. 
\end{xitemize}
Then $\poQ\in M[\genG_0\ast\genH]$ and 
$\tilde{j}\imageof\mu=j\imageof\mu\in M\subseteq M[\genG_0\ast\genH]$ for all
$\mu<j(\kappa)$. 
\smallskip

\assertof{3}: The proof is quite similar to that of \assertof{2}. Starting from a supercompact 
(super almost-huge, resp.) Laver function, we define a FS-iteration of ccc \pos\ 
with the definition of $\utilde{\poQ}_\beta$ similarly to \xitemof{laver-1-18}. 
The step corresponding to the one from \xitemof{laver-1-23} to 
\xitemof{laver-1-24} is now easily done by $[M]^{\aleph_1}\subseteq M$ instead of 
\xitemof{laver-1-21} since we only need to check subsets of a \po\ of size 
$\aleph_1$ to conclude that the \po\ has the ccc. 
\qedofThm
\qedskip

We show in the following that generic large cardinal property of $\kappa$ have 
very strong influence on the cardinal arithmetic around $\kappa$. 

A \po\ $\poP$ is {\it $\omega_1$-preserving} if it satisfies 
$\forces{\poP}{(\omega_1)^\uniV\equiv\omega_1}$. 

\begin{Lemma}
  \Label{L-laver-1}
  Suppose that $\kappa$ is generically measurable by a $\omega_1$-preserving $\poP$. Then $\kappa>\omega_1$. 
\end{Lemma}
\prf Suppose otherwise. Since $\kappa$ cannot be $\omega$, we have then
$\kappa=\omega_1$. Thus there is an $\omega_1$-preserving \po\ $\poP$ and a
$(\uniV,\poP)$-generic filter $\genG$ \st\ there are transitive
$M\subseteq\uniV[\genG]$ and $\elembed{j}{\uniV}{M}$ with
$\crit(j)=\omega_1$. Since $M\modelof{(\omega_1)^\uniV<j(\kappa)=\omega_1}$, 
we have $M\modelof{(\omega_1)^\uniV\mbox{ is countable}}$ and hence
$\uniV[\genG]\modelof{(\omega_1)^\uniV\mbox{ is countable}}$. This is a 
contradiction to the $\omega_1$-preserving of $\poP$. \qedofLemma

\begin{Lemma}
  \Label{L-laver-2}Suppose that $\kappa$ is Laver-generically supercompact for
  $\omega_1$-preserving $\calP$ with $\Col(\omega_1,\ssetof{\omega_2})\in\calP$. Then we 
  have $\kappa=\omega_2$.
\end{Lemma}
\prf $\kappa>\omega_1$ by \Lemmaof{L-laver-1}. Suppose, toward a contradiction, 
that $\kappa>\omega_2$. Let $\poP=\Col(\omega_1,\ssetof{\omega_2})$ and let $\poQ\in\calP$ 
be \st\ $\poP\circleq\poQ$ and, for a $(\uniV,\poQ)$-generic filter $\genH$, 
there are transitive $M\subseteq\uniV[\genH]$  
and $\elembed{j}{\uniV}{M}$ with $\crit(j)=\kappa$. Since $(\omega_2)^\uniV<\kappa$,
$j((\omega_2)^\uniV)=(\omega_2)^\uniV$. Thus, by elementarity, 
$M\modelof{(\omega_2)^\uniV\mbox{ is the second uncountable cardinal}}$. But since the $(\uniV,\poP)$-generic filter 
$\genH\cap\poP$ is in $M$, $M\modelof{\cardof{(\omega_2)^\uniV}=\aleph_1}$. This is 
a contradiction. \qedofLemma\qedskip

\begin{Lemma}
  \Label{T-laver-1} Suppose that $\calP$ is a class of \pos\ containing 
  a \po\/ $\poP$ \st\ any $(\uniV,\poP)$-generic filter $\genG$ codes a new real. If $\kappa$ is a 
  Laver-generically supercompact for $\calP$, then $\kappa\leq2^{\aleph_0}$.  
\end{Lemma}
\prf Let $\poP\in\calP$ be \st\ any generic filter over $\poP$ codes a new real. 

Suppose that $\mu<\kappa$. We show that $2^{\aleph_0}>\mu$. Let
$\veca=\seqof{a_\xi}{\xi<\mu}$ be a sequence of subsets of $\omega$. It 
is enough to show that $\veca$ does not enumerate $\psof{\omega}$.
By Laver-generic supercompactness of $\kappa$ for $\calP$, there are 
$\poQ\in\calP$ with $\poP\circleq\poQ$, 
$(\uniV,\poQ)$-generic $\genH$, transitive 
$M\subseteq\uniV[\genG]$ and $\elembed{j}{\uniV}{M}$ with $\crit(j)=\kappa$
and $\poP,\genH\in M$. Since
$\mu<\kappa$, we have $j(\veca)=\veca$. Since $\genG\in M$ where
$\genG=\genH\cap\poP$ and $\genG$ codes a new real not in $\uniV$, we have
\begin{xitemize}
\xitem[] 
  $M\modelof{j(\veca)\mbox{ does not enumerate }2^{\aleph_0}}$.
\end{xitemize}
By elementarity, it follows that 
\begin{xitemize}
\xitem[] 
  $\uniV\modelof{\veca\mbox{ does not enumerate }2^{\aleph_0}}$.\qedofLemma 
\end{xitemize}

\begin{Lemma}
  \Label{T-laver-1-a} Suppose that $\calP$ is a class of \pos\ \st\ elements of 
$\calP$ do not add any reals. If $\kappa$ is generically supercompact by $\calP$, then 
  we have $\continuum<\kappa$. 
\end{Lemma}
Let $\lambda>2^{\aleph_0}$ and let $\poP\in\calP$ be \st\ there are $j$, 
$M\subseteq\uniV[\genG]$ for a $(\uniV,\poP)$-generic filter $\genG$ with 
$\elembed{j}{\uniV}{M}$, $\crit(j)=\kappa$, $j(\kappa)>\lambda$ and
$j\imageof\lambda\in M$.

Since $\poP$ does not add any new reals
$\uniV[\genG]\modelof{\continuum<j(\kappa)}$. By 
\Lemmaof{L-lt-conti-0},\,\assertof{3}, it follows that
$M\modelof{\continuum<j(\kappa)}$. Thus, by elementarity, $\uniV\modelof{\continuum<\kappa}$.\qedofLemma
\qedskip

For a class $\calP$ of \pos\ and cardinals $\mu$, $\kappa$, we consider the 
following strengthening of the forcing axiom for $\calP$:

\begin{xitemize}
\item[$\MA^{+\mu}(\calP,\LT\kappa)$: ]\it For any $\poP\in\calP$, any family 
$\calD$ of dense subsets of\/ $\poP$ with $\cardof{\calD}<\kappa$ and any family 
$\calS$ of\/ $\poP$-names \st\ $\cardof{\calS}\leq\mu$ and
  $\forces{\poP}{\utilde{S}\xmbox{ is a stationary subset of }\omega_1}$ for all
  $\utilde{S}\in\calS$, there is a $\calD$-generic filter $\genG$ over $\poP$ \st\
  $\utilde{S}[\genG]$ is a stationary subset of $\omega_1$ for all
  $\utilde{S}\in\calS$. 
\end{xitemize}

\memo{added after the 1.submitted version:}
If $\kappa=\aleph_2$, we often drop ``$\LT\kappa$'' and 
write $\MA^{+\mu}(\calP)$. Also the family (class) of \pos\ $\calP$ in this 
notation is oftne identified with the property defining the class $\calP$. 
The notation of the principle ``$\MA^{+\omega_1}(\sigma\mbox{-closed})$'' we used already in 
\cite{I} can be understood an instance of this convention.

For a \po\ $\poP$, $\poP$-name $\utilde{S}$ of a set of subsets of $\On$ and a 
filter $\genG$ on $\poP$, let
\begin{xitemize}
\xitem[laver-1-27] 
  $\utilde{S}(\genG)=\setof{b}{{}
  \begin{array}[t]{@{}l}
    b=\setof{\alpha\in\On}{\condp\forces{\poP}{\check{\alpha}\varin\utilde{s}}
    \mbox{ for a }\condp\in\genG}\mbox{ for a }\poP\mbox{-name }\utilde{s}\\
    \mbox{\st\ }\forces{\poP}{\utilde{s}\varin\utilde{S}\mbox{ and }\sup(\utilde{s})\equiv\sup(b)}}.
  \end{array}$
\end{xitemize}

Note that if $\genG$ is a $(\uniV,\poP)$-generic filter, then
$\utilde{S}(\genG)=\utilde{S}[\genG]$.
\ifextended{\tt
  [ ...]
  }\fi

For uncountable cardinals $\mu$ and $\kappa>\aleph_1$, let
$\MA^{++\mu}(\calP,\LT\kappa)$ be the strengthening 
of $\MA^{+\mu}(\calP,\LT\kappa)$ defined by:
\begin{xitemize}
\item[$\MA^{++\mu}(\calP,\LT\kappa)$: ]\it For any $\poP\in\calP$, any family 
  $\calD$ of dense subsets of\/ $\poP$ with $\cardof{\calD}<\kappa$ and any family 
  $\calS$ of\/ $\poP$-names \st\ $\cardof{\calS}\leq\mu$ and
  $\forces{\poP}{\utilde{S}\xmbox{ is a stationary subset of }\Pkl{\eta_{\scriptstyle\utilde{S}}}{\theta_{\utilde{S}}}}$ 
  for some $\omega<\eta_{\utilde{S}}\leq\theta_{\utilde{S}}\leq\mu$ with $\eta_{\utilde{S}}$ regular, for all
  $\utilde{S}\in\calS$, there is a $\calD$-generic filter $\genG$ over $\poP$ \st\
  $\utilde{S}(\genG)$ is stationary in $\Pkl{\eta_{\scriptstyle\utilde{S}}}{\theta_{\utilde{S}}}$ for all
  $\utilde{S}\in\calS$. 
\end{xitemize}

Clearly $\MA^{++\omega_1}(\calP,\LT\kappa)$ is 
equivalent to $\MA^{+\omega_1}(\calP,\LT\kappa)$. 

For \Thmof{T-laver-0} below, we need a slight strengthening of Laver-generic 
supercompactness which also hold in the canonical models of the original 
Laver-generic spupercompactness:
For a notion of \pos\ $\calP$ which is closed \wrt\ two-step iteration (i.e, if
$\poP\models\calP$ and $\forces{\poP}{\utpoQ\models\calP}$ then
$\poP\ast\utpoQ\models\calP$), a cardinal $\kappa$ is said to be {\it strongly 
  Laver-generically supercompact for $\calP$}, if 
for any $\lambda\geq\kappa$ and any
$\poP\models\calP$, there are a $\poP$-name of a \po\ $\utpoQ$  with $\forces{\poP}{\utpoQ\models\calP}$
and $(\uniV,\poP\ast\utpoQ)$-generic filter $\genH$ \st\ there 
are an inner model
$M\subseteq\uniV[\genH]$ and a class $j\subseteq\uniV[\genH]$ with 
\begin{xitemize}
\xitemd[laver-a]{{}} $\elembed{j}{\uniV}{M}$, 
\xitemd[laver-0]{{}} $\crit(j)=\kappa$, $j(\kappa)
  >\lambda$,
\xitemd[laver-0-0]{{}} $\poP$, $\genH\in M$ and 
\xitemd[laver-1]{{}} $j\imageof\lambda\in M$.
\end{xitemize}

\memo{added after the 1.submitted version:}
\ifextended{\tt
  Recall that, for \pos\ $\poP$, $\poQ$ with $\poP\circleq\poQ$ 
  and $\condq\in\poP$,  
  there is always a condition $\condp\in\poP$ \st, for all $\condr\leq_\poP\condp$,
  $\condr$ and $\condq$ are compatible in $\poQ$. Such $\condp$ is called a
  {\it reduction of\/ $\condq$} (in $\poP$). 

  Actually, the condition $\poP\circleq\poQ$ (i.e.\ $\poP$ is a regular sub-\po\ of
  $\poQ$\ $\Leftrightarrow$\ $\poP$ is a sub-\po\ of $\poP$ \st\ the identity 
  mapping on it is a complete embedding of $\poP$ to $\poQ$) is equivalent to the 
  existence of a reduction $\condp\in\poP$ of $p$ for all $\condq\in\poQ$.

  The following \LemmaAof{laver-A-0} and \LemmaAof{laver-A-1} should be well-known. 

  \begin{LemmaA}
    \Label{laver-A-0}Suppose that\/ $\poP$, $\poQ$ are \pos\ 
    with $\poP\circleq\poQ$, $\varphi=\varphi(\variables{x}{n-1})$ is a
    $\Delta^\ZF_0$-formula, and\/ $\variables{\uta}{n-1}$ are  $\poP$-names.\smallskip

    \assert{0} $\variables{\uta}{n-1}$ are also $\poQ$-names.\smallskip

    \assert{1} For any $\condp\in\poP$, $\condp\forces{\poP}{\varphi(\variables{\uta}{n-1})}$ if and only 
    if\/ $\condp\forces{\poQ}{\varphi(\variables{\uta}{n-1})}$.\smallskip

    \assert{2} For any $\condq\in\poQ$, if\/ $\condp\in\poP$ is a reduction of
    $\condq\in\poQ$, then $\condq\forces{\poQ}{\varphi(\variables{\uta}{n-1})}$ implies 
    $\condp\forces{\poP}{\varphi(\variables{\uta}{n-1})}$.\smallskip
  \end{LemmaA}
  \prf \assertof{0}: We show by induction on the rank of $\poP$-names that all
  $\poP$-namese are also $\poQ$-names. \smallskip

  \assertof{1}: \Wolog, we may assume that $\poP$ and $\poQ$ are cBa \pos\ (i.e. 
  \pos\ of the form $\baB^{+}$ for a complete \Ba\ $\baB$). Then $\poQ$ embeds 
  in $\poP\ast\utpoR$ densely over $\poP$ for some $\poP$-name $\utpoR$ of a \po.
  Thus we see that, for all $(\uniV,\poP)$-generic filter $\genG$, there is a
  $(\uniV,\poQ)$-generic filter $\genH$ with $\genG\subseteq\genH$. 

  For $\condp\in\poP$, if $\condp\forces{\poP}{\varphi(\variables{\uta}{n-1})}$, 
  then for any $(\uniV,\poQ)$-generic $\genH$ with $\condp\in\genH$,
  $\genG=\genH\cap\poP$ is a $(\uniV,\poP)$-generic filter 
  with $\condp\in\genG$. Thus, by Forcing Theorem, we have
  $\varphi(\uta_0[\genG]\ctentenc\uta_{n-1}[\genG])$. Since
  $\uta_0[\genG]=\uta_0[\genH]$\ctentenc $\uta_{n-1}[\genG]=\uta_0[\genH]$, it 
  follows that $\varphi(\uta_0[\genH]\ctentenc\uta_{n-1}[\genH])$. Thus, again by 
  Forcing Theorem, we have $\condp\forces{\poQ}{\varphi(\variables{\uta}{n-1})}$. 

  Suppose now that $\condp\forces{\poQ}{\varphi(\variables{\uta}{n-1})}$. For an 
  arbitrary $(\uniV,\poP)$-generic filter $\genG$ with $\condp\in\genG$, let 
  $\genH$ be a $(\uniV,\poQ)$-generic filter \st\ $\genG\subseteq\genH$. By the 
  assumption and the Forcing Theorem we have
  $\varphi(\uta_0[\genH]\ctentenc\uta_0[\genH])$. Since
  $\uta_0[\genG]=\uta_0[\genH]$\ctentenc $\uta_{n-1}[\genG]=\uta_0[\genH]$, it  
  follows that $\varphi(\uta_0[\genG]\ctentenc\uta_0[\genG])$. Since $\genG$ is 
  arbitrary with $\condp\in\genG$, if follows by Forcing Theorem that
  $\condp\forces{\poP}{\varphi(\variables{\uta}{n-1})}$. \smallskip

  \assertof{2}: Suppose that $\condq\in\poQ$ and 
  \begin{xitemize}
  \xitemA[laverA-4] $\condq\forces{\poQ}{\varphi(\variables{\uta}{n-1})}$. 
  \end{xitemize}
  Assume, toward a contradiction, that $\condp$ is a reduction of $\condq$ in
  $\poP$ and $\condp\notforces{\poP}{\varphi(\variables{\uta}{n-1})}$. Then, 
  there is $\condr\leq_\poP\condp$ \st\
  $\condr\forces{\poP}{\neg\varphi(\variables{\uta}{n-1})}$.
  by \assertof{1}, it follows that 
  \begin{xitemize}
  \xitemA[laverA-5] $\condr\forces{\poQ}{\neg\varphi(\variables{\uta}{n-1})}$. 
  \end{xitemize}

  Since $\condp$ is a reduction of $\condq$ in $\poP$, $\condr$ and $\condq$ are 
  compatible in $\poQ$. Let $\condr'\leq_\poQ\condr$, $\condq$. Then,
  $\condr'\forces{\poQ}{\varphi(\variables{\uta}{n-1})}$ by \xitemof{laverA-4}. 
  On the other hand, we have $\condr'\forces{\poQ}{\neg\varphi(\variables{\uta}{n-1})}$
  by \xitemof{laverA-5}. This is a contradiction. \qedofLemmaA
  \begin{LemmaA}
    \Label{laver-A-1}
    Suppose that $\kappa$ is a regular uncountable cardinal, and
    $\poP$, $\poQ$ are $\kappa$-cc \pos\ with $\poP\circleq\poQ$. 
    If $\utS$ is a $\poP$-name with
    $\forces{\poP}{\utS\xmbox{ is a stationary subset of }\Pkl{}{}}$, then we 
    also have $\forces{\poQ}{\utS\xmbox{ is a stationary subset of }\Pkl{}{}}$ 
    for some $\lambda\in\On$. 
  \end{LemmaA}
  \prf Suppose that $\utS$ is a $\poP$-name of a stationary subset of $\Pkl{}{}$. 
  Toward a contradiction, assume that there is a $\poQ$-name $\utC$ of club 
  subset of $\Pkl{}{}$ \st\ 
  \begin{xitemize}
    \xitemA[laverA-6] $\forces{\poQ}{\utS\cap\utC\equiv\emptyset}$. 
  \end{xitemize}

  Let $C=\setof{x\in\Pkl{}{}}{\forces{\poQ}{\check{x}\varin\utC}}$. by 
  the $\kappa$-cc of $\poQ$, $C$ is a club $\subseteq\Pkl{}{}$ and
  $\forces{\poQ}{\check{C}\subseteq\utC}$. Let $\utC'$ be a $\poP$-name of the 
  club (in $\uniV^\poP$) generated from $\check{C}_\poP$. Then we have
  $\forces{\poQ}{\utC'\subseteq\utC}$. Thus
  $\forces{\poQ}{\utS\cap \utC'\equiv\emptyset}$. By 
  \Lemmaof{laver-A-0},\,\assertof{1}, it follows that
  $\forces{\poP}{\utS\cap \utC'\equiv\emptyset}$. This is a contradiction to the 
  choice of $\utS$.
  \qedofLemmaA
}\fi

\begin{Thm}
  \Label{T-laver-0} \wassertof{1} For a class $\calP$ of 
  ccc \pos, if
  $\kappa>\aleph_1$ is Laver-generically supercompact for $\calP$, 
  then $\MA^{++\mu}(\calP,\LT\kappa)$ holds for all $\mu<\kappa$.\smallskip

  \assert{2} If\/ $\aleph_2$ is strongly Laver-generically supercompact for a 
  notion $\calP$ of \pos\ which preserves stationarity of subsets of $\omega_1$ 
  and which is closed \wrt\ two-step iteration, then $\MA^{+\omega_1}(\calP)$ holds.
\end{Thm}
\prf \assertof{1} and \assertof{2} can be proved by practically the 
same proof. Let $\calP$ as in one of \assertof{1} or \assertof{2}, 
$\poP\in\calP$ and $\mu<\kappa$. Let $\calD$ and $\calS$ be as in 
the definition of $\MA^{++\mu}(\calP,\LT\kappa)$. \Wolog, we may assume that the underlying set of 
$\poP$ is some cardinal $\lambda_0$ and elements of $\calS$ are 
nice $\poP$-names.

Let $\lambda>\lambda_0$ be sufficiently large and let $\poQ\in\calP$ be \st\ 
$\poP\circleq\poQ$ (in case of \assertof{2}, \ifextended
{\begin{xitemize}
\tt\xitemA[laverA-20] 
  $\poQ=\poP\ast\utpoR$ for 
  a $\poP$-name $\utpoR$ \st\ $\forces{\poP}{\utpoR\models\calP}$) 
\end{xitemize}}\noindent\else
  $\poQ=\poP\ast\utpoR$ for 
  a $\poP$-name $\utpoR$ \st\ $\forces{\poP}{\utpoR\models\calP}$) 
\fi
and, for a $(\uniV,\poQ)$-generic filter $\genH$, there are  
transitive $M\subseteq\uniV[\genH]$ and $\elembed{j}{\uniV}{M}$ with
\begin{xitemize}
\xitemd[laver-0]{{}} $\crit(j)=\kappa$, $j(\kappa)
  >\lambda$,
\xitemd[laver-0-0]{{}} $\poP$, $\genH\in M$ and 
\xitemd[laver-1]{{}} $j\imageof\lambda\in M$.
\end{xitemize}



By the choice of $\lambda$, \xitemof{laver-1} and \Lemmaof{L-lt-conti-0},\,\assertof{5}, we 
have $\poP$, $\calD$, $\calS\in M$. Let $\genG=\genH\cap\poP$. Then
$\genG\in M$ by \xitemof{laver-0-0}. 
Thus $\genG$ witnesses 

\begin{xitemize}
\xitem[laver-1-3] 
  $M\modelof{
  \begin{array}[t]{@{}l}
    \mbox{there is a }\calD\mbox{-generic filter }G\mbox{ over }\poP
    \\
    \mbox{\st\ }\utilde{S}(G)\mbox{ is a stationary subset of }\Pkl{\eta_{\scriptstyle\utilde{S}}}{\theta_{\utilde{S}}}
    \mbox{ for all }\utilde{S}\in\calS}.
  \end{array}$
\end{xitemize}
\ifextended{\tt 
This follows from \LemmaAof{laver-A-1} in case of \assertof{1} or 
from \xitemAof{laverA-20} in case of \assertof{2}.
}
\else 
Note that it is here that we need separate settings in the cases \assertof{1}, 
and \assertof{2}. 
\fi 

Since $j(\calD)=\setof{j(D)}{D\in\calD}$ and 
$j(\calS)=\setof{j(S)}{S\in\calS}$ by $\cardof{\calD}$, $\cardof{\calS}<\kappa$, 
$j(D)\supseteq j\imageof D$ for all $D\in\calD$, $j(S)\supseteq j\imageof S$ 
for all $S\in\calS$ and $j\imageof\genG\in M$ by 
\Lemmaof{L-lt-conti-0},\,\assertof{6}, 
it follows that 
\begin{xitemize}
\xitem[laver-1-4] 
  $M\modelof{
  \begin{array}[t]{@{}l}
    \mbox{there is a }j(\calD)\mbox{-generic filter }G\mbox{ over }j(\poP)
    \\
    \mbox{\st\ }\utilde{S}(G)\mbox{ is a stationary subset 
      of }\Pkl{\eta_{\scriptstyle \utilde{S}}}{\theta_{\utilde{S}}}
    \mbox{ for all }\utilde{S}\in j(\calS)}.
  \end{array}$
\end{xitemize}

By elementarity, it follows that 
\begin{xitemize}
\xitem[laver-1-5] 
  $\uniV\modelof{
  \begin{array}[t]{@{}l}
    \mbox{there is a }\calD\mbox{-generic filter }G\mbox{ over }\poP
    \\
    \mbox{\st\ }\utilde{S}(G)\mbox{ is a stationary subset of }\Pkl{\eta_{\scriptstyle\utilde{S}}}{\theta_{\utilde{S}}}
    \mbox{ for all }\utilde{S}\in \calS}.
  \end{array}$\vspace{-0.6ex}
\end{xitemize}\mbox{}\qedofThm
\qedskip






At the moment we do not know if Laver-generic supercompactness of $\kappa$ for 
ccc \pos\ implies $\kappa=2^{\aleph_0}$ (note that we have $\kappa\leq 2^{\aleph_0}$ by 
\Lemmaof{T-laver-1}). However

\begin{Thm}
  \Label{T-laver-1-2} If $\kappa$ is tightly Laver-generically superhuge for ccc 
  \pos, then $\kappa=2^{\aleph_0}$. 
\end{Thm}
\prf Suppose that $\kappa$ is tightly Laver-generically superhuge for ccc \pos. 
By \Lemmaof{T-laver-1}, we have $\continuum\geq\kappa$.

To prove $\continuum\leq\kappa$, let $\lambda\geq\kappa$, $2^{\aleph_0}$ be 
large enough and let $\poQ$ be a ccc \po\ \st\ there are $(\uniV,\poQ)$-generic
$\genH$ and $\elembed{j}{\uniV}{M}\subseteq\uniV[\genH]$ with $\crit(j)=\kappa$,
$\cardof{\poQ}\leq j(\kappa)>\lambda$, $\genH\in M$ and
$j\imageof j(\kappa)\in M$. 

Since $M\modelof{j(\kappa)\mbox{ is regular}}$ by elementarity, $j(\kappa)$ is 
regular in $\uniV$ (e.g.\,by \Lemmaof{L-lt-conti-0},\,\assertof{3}). 
Thus, we have
$\uniV\modelof{j(\kappa)^{\aleph_0}=j(\kappa)}$  by \Propof{T-laver-1-1},\,\assertof{1}. Since $\poQ$ has the ccc and
$\cardof{\poQ}\leq j(\kappa)$, it follows that
$\uniV[\genG]\modelof{\continuum\leq j(\kappa)}$. Now by 
\Lemmaof{L-lt-conti-0},\,\assertof{4}, $(j(\kappa)^+)^M=(j(\kappa)^+)^{\uniV[\genG]}$. 
Thus $M\modelof{\continuum\leq j(\kappa)}$.

By elementarity, it follows that $\uniV\modelof{\continuum\leq\kappa}$.
\qedofThm

\begin{Thm}
  \Label{T-laver-2} \wassertof{1} Suppose that $\kappa$ is strongly Laver-generically 
  supercompact for $\sigma$-closed \pos. Then $2^{\aleph_0}=\aleph_1$,
  $\kappa=\aleph_2$, $\MA^{+\omega_1}(\sigma\mbox{-closed})$ and hence also
  $\SDLS(\calL^{\aleph_0}_{stat},\LT\aleph_2)$ holds.\smallskip

  \assert{2} Suppose that $\kappa$ is strongly Laver-generically supercompact for proper \pos. Then
  $2^{\aleph_0}=\kappa=\aleph_2$, $\PFA^{+\omega_1}$ and hence also 
  $\SDLS^-(\calL^{\aleph_0}_{stat},\LT\continuum)$ holds.\smallskip

  \assert{3} Suppose that $\kappa$ is Laver-generically supercompact for ccc 
  \pos. Then $2^{\aleph_0}\geq\kappa$ and $\Pkl{}{}$ for any regular $\lambda\geq\kappa$ carries an
  $\aleph_1$-saturated normal ideal. In particular, $\kappa$ is $\kappa$-weakly Mahlo.
  $\MA^{++\mu}(ccc,\LT\kappa)$ for all $\mu<\kappa$,  
  $\SDLS^{\intnl}(\calL^{\aleph_0}_{stat},\LT\kappa)$ and 
  $\SDLS^{\intnl}(\calL^{PKL}_{stat},\LT\kappa)$ also hold.
\end{Thm}
\prf \assertof{1}: Assume that $\kappa$ is Laver-generically supercompact for
$\sigma$-closed \pos. Then $\kappa=\aleph_2$ by \Lemmaof{L-laver-2}. Hence 
$\continuum=\aleph_1$ by \Lemmaof{T-laver-1-a}.
$\MA^{+\omega_1}(\sigma\mbox{-closed})$ holds by \Thmof{T-laver-0}. 
$\DRP(\IC_{\aleph_0})$ follows from $\MA^{+\omega_1}(\sigma\mbox{-closed})$ (see  
\cite{cox:}). 
Hence, by \Corof{P-Einf-2},\,\assertof{4}, 
$\SDLS(\calL^{\aleph_0}_{stat},\LT\aleph_2)$ holds. \smallskip

\assertof{2}: Assume that $\kappa$ is Laver-generically supercompact for
proper \pos. Then 
$\kappa=\aleph_2$ by \Lemmaof{L-laver-2}. $\kappa\leq\continuum$ by 
\Lemmaof{T-laver-1} and $\PFA^{+\omega_1}$ by \Thmof{T-laver-0}. 
Since \PFA\ implies
$\continuum=\aleph_2$, we obtain $\kappa=\continuum$. 

$\PFA^{+\omega_1}$ implies $\MA^{+\omega_1}(\sigma\mbox{-closed})$ and 
$\DRP(\IC_{\aleph_0})$ follows from $\MA^{+\omega_1}(\sigma\mbox{-closed})$ (\cite{cox:}). Thus, 
by  \Corof{P-Einf-2},\,\assertof{3}, 
$\SDLS^-(\calL^{\aleph_0}_{stat},\LT\aleph_2)$ holds. \smallskip

\assertof{3}: Assume that $\kappa$ is Laver-generically supercompact for
ccc \pos. Then $\kappa\leq\continuum$ by \Lemmaof{T-laver-1}.

$\SDLS^{\intnl}(\calL^{\aleph_0}_{stat},\LT\kappa)$ holds by \Thmof{P-lt-conti-3} 
and \Propof{P-internal-0}. $\SDLS^{\intnl}(\calL^{PKL}_{stat},\LT\kappa)$ holds 
by \Thmof{P-lt-conti-3:} and \Propof{P-PKL-0}.
For any regular $\lambda\geq\kappa$
$\Pkl{}{}$ carries an $\aleph_1$-saturated normal ideal by \Lemmaof{L-lt-conti-1-1},\,\assertof{2}.
$\MA^{++\mu}(ccc,\LT\kappa)$ for all $\mu<\kappa$ by \Thmof{T-laver-0}.
\qedofThm


\begin{thebibliography}{99}
\Label{ref}%
\newcommand{\bysame}[1]{\underline{\phantom{#1}}}%
\bibitem{baumgartner-taylor:} James E.\ Baumgartner and Alan D.\ Taylor, 
  Saturation properties of Ideals in generic extensions I, Transactions of the 
  American Mathematical Society, Vol.270, No.2, (1982), 557--574.
\bibitem{brickhill-fuchino-sakai} Hazel Brickhill, Saka\'e Fuchino and Hiroshi 
  Sakai, On $\alpha$-stationary subsetes of $\Pkl{}{}$, in preparation. 
\bibitem{corazza} Paul Corazza, Laver Sequences for Extendible and 
  Super-Almost-Huge Cardinals, The Journal of Symbolic Logic, Vol.64, No.3 
  (1999), 963--983. 
\bibitem{cox:} Sean Cox, The diagonal reflection principle, 
Proceedings of the American Mathematical Society, Vol.140, No.8 (2012), 
2893-2902. 
\bibitem{cox2} \bysame{Sean Cox}, Forcing axioms, approachability, and stationary set 
  reflection, preprint.
\bibitem{I} Saka\'e Fuchino, Andr\'e Ottenbreit Maschio Rodrigues and Hiroshi 
Sakai, Strong downward L\"owenheim-Skolem theorems for stationary logics, I, 
to appear in Archive for Mathematical Logic.\ifextended\\
\href{https://fuchino.ddo.jp/papers/SDLS-x.pdf}{https://fuchino.ddo.jp/papers/SDLS-x.pdf}\fi
\bibitem{III} \bysame{Saka\'e Fuchino, eta}, Strong 
downward L\"owenheim-Skolem theorems for stationary logics, III  --- mixed support iteration, sbumitted.
\ifextended\\
\href{https://fuchino.ddo.jp/papers/SDLS-x.pdf}{http://fuchino.ddo.jp/papers/SDLS-III-x.pdf}\fi
\bibitem{IV} \bysame{Saka\'e Fuchino, eta}, Strong 
downward L\"owenheim-Skolem theorems for stationary logics, IV --- more on 
Laver-generically large cardinals and mixed support iterations, in preparation.
\bibitem{millennium-book:} T.\ Jech, Set Theory, The Third Millennium
	Edition, Springer (2001/2006).
\bibitem{higher-inf:} Akihiro Kanamori, The Higher Infinite, 
  Springer--Verlag (1994/2003).
\bibitem{koenig:} Bernhard K\"onig, Generic compactness reformulated, Archive of 
  Mathematical Logic 43, (2004), 311--326. 
\bibitem{MaSaUs} Yo Matsubara, Hiroshi Sakai and Toshimichi Usuba,  On the 
  existence of skinny stationary subsets, Annals of Pure and Applied Logic 170, (2019) 539--557.


\ifextended{\tt
  \bibitem{kunen-book}  Kenneth Kunen, Set Theory, An Introduction to Independence 
	Proofs, Elsevier (1980). 

}\fi 
\end{thebibliography}
\end{document}